\documentclass[a4paper,11pt,french,leqno]{smfart}

\usepackage{amsmath}
\usepackage{amssymb}
\usepackage{amsthm}
\usepackage{amsopn}
\usepackage{amscd}
\usepackage{fullpage}
\usepackage[utf8]{inputenc}
\usepackage{graphics, xcolor}
\usepackage{textcomp} 
\usepackage[all]{xy}
\usepackage{lscape}
\usepackage{url}
\usepackage{verbatim}
\usepackage{mathrsfs}
   \usepackage[francais]{babel}
\usepackage{tikz}
\usepackage{tikz-cd}

\usepackage{bm}

\def\commentaire#1{\ifx\cachercommentaires\undefined \textcolor{red}{#1}\else \fi } 
\newcommand{\scrF}{\mathscr{F}}

\newcommand{\frS}{\mathfrak{S}}

\newcommand{\frU}{\mathfrak{U}}

\newcommand{\frZ}{\mathfrak{Z}}

\newcommand{\fra}{\mathfrak{a}}
\newcommand{\frb}{\mathfrak{b}}

\newcommand{\frg}{\mathfrak{g}}
\newcommand{\frh}{\mathfrak{h}}

\newcommand{\frl}{\mathfrak{l}}
\newcommand{\frmm}{\mathfrak{m}}

\newcommand{\fro}{\mathfrak{o}}
\newcommand{\frp}{\mathfrak{p}}

\newcommand{\frs}{\mathfrak{s}}
\newcommand{\frt}{\mathfrak{t}}

\newcommand{\bbA}{\mathbb{A}}

\newcommand{\bbC}{\mathbb{C}}

\newcommand{\bbN}{\mathbb{N}}

\newcommand{\bbR}{\mathbb{R}}

\newcommand{\bbZ}{\mathbb{Z}}

\newcommand{\caA}{\mathcal{A}}
\newcommand{\caB}{\mathcal{B}}

\newcommand{\caF}{\mathcal{F}}

\newcommand{\caH}{\mathcal{H}}

\newcommand{\caM}{\mathcal{M}}

\newcommand{\caO}{\mathcal{O}}
\newcommand{\caP}{\mathcal{P}}

\newcommand{\caU}{\mathcal{U}}

\newcommand{\GL}{\mathbf{GL}}

\newcommand{\Ind}{\mathrm{Ind}}

\newcommand{\Id}{\mathrm{Id}}

\renewcommand{\Ind}{\mathrm{Ind}}

\newcommand{\SL}{\mathbf{SL}}

\newcommand{\SO}{\mathbf{SO}}
\newcommand{\Sp}{\mathbf{Sp}}
\newcommand{\Or}{\mathbf{O}}

\newcommand{\Mp}{\mathbf{Mp}}

\newcommand{\Tr}{\mathrm{Tr}}

\newcommand{\bil}[2]{\langle  #1,#2 \rangle }

\newcommand{\Triv}{\mathbf{Triv}}

\newcommand{\Std}{\mathbf{Std}}

\theoremstyle{plain}
\newtheorem{thm}{Théorème}[section]
\newtheorem{lemme}[thm]{Lemme}
\newtheorem{cor}[thm]{Corollaire}
\newtheorem{prop}[thm]{Proposition}

\theoremstyle{definition}
 \newtheorem{defi}[thm]{Définition}
 
\newtheorem{rmqs}[thm]{Remarques}
\newtheorem{rmq}[thm]{Remarque}
\newtheorem{exemple}[thm]{Exemple}

\def \dem {\noindent \underline{\sl Démonstration}. }

\begin{document}

\setcounter{MaxMatrixCols}{30}

\numberwithin{equation}{section}

\title{Paquets d'Arthur des groupes classiques  complexes}

 \author{Colette Moeglin}
 \address{CNRS, Institut Mathématique de Jussieu } 
 \email{colette.moeglin@imj-prg.fr}

  \author{David Renard}
 \address{Centre de Mathématiques Laurent Schwartz,  Ecole Polytechnique} 
\email{david.renard@polytechnique.edu}

 \date{\today}

\thanks{Le deuxième auteur a bénéficié d'une aide de  l'agence nationale de la recherche 
ANR-13-BS01-0012 FERPLAY}

\begin{abstract} 
Nous décrivons explicitement les paquets d'Arthur des groupes classiques complexes, ainsi que
 leur paramétrisation interne
par les  caractères du groupe des composantes connexes du centralisateur de leur paramètre. Nous montrons d'abord 
qu'ils sont obtenus par induction parabolique préservant l'irréductibilité à partir des paquets
 unipotents de \og bonne parité \fg. 
Pour ceux-ci, nous montrons qu'ils coïncident avec les paquets définis par Barbasch-Vogan \cite{BV}. Nous utilisons
des résultats profonds de Barbasch entrant dans sa classification du dual unitaire de ces groupes \cite{B1}.

\noindent \textit{\textbf{Abstract}}. --- We describe explicitly Arthur packets for complex classical groups, 
as well as their internal parametrization by the group of characters  of the component group of the 
stabilizer of their parameter.  We first show that they are obtained by parabolic induction 
preserving irreducibility  from unipotent packets of ``good parity''.
 For these, we show that they coincide with the packets defined by Barbasch and  Vogan. 
 We use deep results of Barbasch
 entering his classification of the unitary dual of these groups \cite{B1}.
 \end{abstract}

\maketitle

\tableofcontents

\section{Introduction}

Soit $G$ un groupe classique complexe. Le but de cet article est de décrire de manière aussi explicite que
possible les représentations irréductibles de $G$ (ou plus exactement, les modules de Harish-Chandra)
qui sont composantes locales d'une représentation automorphe de carré intégrable.
Plus précisément, soit $F$ un corps de nombre, et l'on suppose qu'il existe une place archimédienne $v_0$ de $F$ 
telle que $F_{v_0}$ soit isomorphe à $\bbC$. 
Soit $\mathbf{G}$ un groupe algébrique réductif défini sur $F$, que l'on suppose déployé,
et tel que $G$ soit isomorphe à $\mathbf{G}(F_{v_0})$. Les représentations  automorphes de carré
 intégrable de $\mathbf{G}(\bbA_F)$
sont les sous-représentations irréductibles de $\mathbf{G}(\bbA_F)$ dans 
$L^2(\mathbf{G}(F)\backslash \mathbf{G}(\bbA_F))$.
Dans \cite{Art13}, J. Arthur décrit ces représentations. Il montre que leurs composantes locales en une place $v$
se regroupent en \og paquets \fg\, associés à certains homomorphismes $\psi_v :\; W_{F_v} \rightarrow {}^LG$
du groupe de Weil (ou de Weil-Deligne pour les places non-archimédiennes) $W_{F_v}$ vers le $L$-groupe 
 ${}^LG$. Lorsqu'on fixe un tel homomorphisme $\psi_v :\; W_{F_v} \rightarrow {}^LG$, c'est-à-dire 
  un paramètre d'Arthur,  les représentations qui lui sont associées, et qui forment le  \og paquet
  d'Arthur \fg \, noté $\Pi(\psi_v,\mathbf{G}(F_v))$, sont déterminées par des propriétés locales.

 Les homomorphismes $\psi_v$ sont conjecturalement les localisations d'un homomorphisme global
  $\psi$ dont la définition
  précise est encore lointaine, car elle suppose d'avoir montré que la catégorie des 
  représentations automorphes isobariques des groupes $\GL_N$ possède une structure de type tannakien 
  qui ne semble  pas accessible pour l'instant ({\sl cf.} \cite{Lang79}).
   Mais Arthur en a donné un substitut commode en remplaçant les homomorphismes
  du groupe tannakien conjectural vers ${}^LG$ par les classes d'isomorphie de représentations automorphes cuspidales
  des groupes $\GL_N$. On peut localiser ces paramètres globaux et les paquets locaux ne dépendent que de
  cette localisation. D'autre part, une représentation $\pi_v$ est dans le paquet local associé 
  à $\psi_v$ s'il existe
 une représentation automorphe de carré intégrable, associée à  un paramètre global $\psi$ 
  dont la composante locale en $v$ est $\psi_v$, et qui admet $\pi_v$ comme composante locale en $v$.
Le paramètre global se lit sur ce qui se passe aux  places non ramifiées et est donc caractérisé de manière 
extr\^emement simple. Cette présentation sommaire ne donne bien sûr  qu'un aperçu des résultats de \cite{Art13}.
Pour en dire un peu plus, précisons
qu'à chaque caractère $\eta$ du groupe des composantes connexes du 
 centralisateur de $\psi_v$ dans ${}^LG$, Arthur associe une somme directe 
$X_\eta$ (avec éventuellement des multiplicités et qui peut-\^etre nulle) de représentations irréductibles
 de $\mathbf{G}(F_v)$   et il écrit une formule de multiplicité globale qui fait intervenir ces caractères locaux.
Les représentations $X_\eta$ sont uniquement déterminées par des propriétés de transfert endoscopique
 expliquées en \cite{Art13} 2.2.1 (b) et 1.5.1.

Revenons maintenant à notre place complexe $v_0$ telle que $G$ soit isomorphe à $\mathbf{G}(F_{v_0})$, et 
abandonnons l'indice $v_0$ pour les paramètres d'Arthur : $\psi :\, W_\bbC \rightarrow {}^LG$.
Rappelons que $W_\bbC=\bbC^ \times$. Dans ce
travail, nous déterminons donc les représentations associées à un tel paramètre $\psi$, et en s'appuyant 
fortement sur les travaux de Barbasch et Vogan \cite{BV} et de Barbasch \cite{B1}. En fait on démontre que le cas 
général se ramène par une induction irréductible explicite au cas traité en \cite{BV} appelé le cas 
spécial unipotent.  
 Nous montrons au passage que pour $\psi$ fixé, les représentations associées à un 
tel  caractère $\eta$ sont disjointes les unes des autres; de plus elles sont 
irréductibles ou nulles ce qui est très particulier aux places complexes.

Pour décrire un peu plus précisément nos résultats, 
introduisons,  comme ci-dessus, le morphisme de $W_\bbC\times \SL_2({\mathbb C})$ dans $^LG$
qui paramétrise les constructions de \cite{Art13} (cf. définition \ref{ArtPar} dans
le texte).   En voyant $\psi$ comme une représentation de $W_\bbC\times
\SL_2({\mathbb C})$, on vérifie aisément l'existence d'une décomposition
de $\psi$ (cf.(\ref{decomposition}) dans le texte) sous la
forme $\rho \oplus \rho^* \oplus \psi_{u,bp}$,  où $\rho$ est une
représentation dont $\rho^*$ est duale, où $\psi_{u,bp}$ est
caractérisé par le fait que sa restriction  à $W_\bbC$
est triviale et  l'orbite unipotente définie par sa restriction 
 aux éléments unipotents principaux de $\SL_2({\mathbb C})$
a des blocs de Jordan tous de même parité (paire si $^LG$ est un groupe
symplectique et impaire sinon). D'autre part $\psi_{u,bp}$ est maximal 
avec ces propriétés. On note
$m_{u,bp}$ la dimension de la représentation  $\psi_{u,bp}$ et
$G_{u,bp}$ le groupe de même type que $G$ mais de rang $\lfloor m_{u,bp}/2 \rfloor$.
Alors \cite{BV} associe à $G_{u,bp}$ et $\psi_{u,bp}$ un ensemble de
représentations à l'aide de formules de caractères. D'autre part la
classification de Langlands associe à $\rho$ une représentation
irréductible de $\GL_d({\mathbb C})$, notée $\tau$, où $d$ est la
dimension de la représentation $\rho$. Nous montrons que les
représentations associées par Arthur à $\psi$ sont exactement les
représentations induites de $\tau\otimes X_{u,bp}$ pour le parabolique
maximal $\GL_d({\mathbb C})\times G_{u,bp}$ de $G$ ou $X_{u,bp}$
parcourt l'ensemble des représentations associées par \cite{BV} à
$\psi_{u,bp}$. On va même plus loin dans les paramétrisations: \cite{BV}
paramétrise les représentations qu'ils construisent par les caractères
du quotient de Lusztig du groupe des composantes connexes du centralisateur de
$\psi_{u,bp}$ dans $^LG_{u,bp}$; avec le résultat précédent on a donc
une paramétrisation des représentations associées par \cite{Art13} à
$\psi$ et nous montrons que cette paramétrisation coïncide avec celle de
\cite{Art13}. Les propriétés d'irréductibilité sont tirées de \cite{B1}
et l'identification des constructions de \cite{BV} avec celles de
\cite{Art13} résulte des formules de transfert endoscopique.  Pour
éviter de devoir démontrer  la formule de transfert endoscopique tordue
qui n'est pas dans \cite{BV}, on utilise un résultat tiré de
\cite{pourhowe} rappelé dans la remarque  \ref{mult1} ci-dessous.

Le cas des groupes complexes est intéressant car il illustre des phénomènes pouvant paraître surprenants. Par exemple
sur $\bbC$, la notion de stabilité est triviale : toutes les distributions invariantes sur $G$ sont 
stablement invariantes, et c'est le cas en particulier des caractères des représentations irréductibles.
On pourrait en conclure trop rapidement que, à l'instar de ce qui se passe pour les groupes 
$\GL_N$, les paquets d'Arthur sont des singletons.  Mais ceci n'est pas vrai en général, car il faut aussi qu'une
autre propriété fondamentale soit satisfaite, à savoir la compatibilité des constructions d'Arthur 
avec l'induction parabolique. Plus précisément, il s'agit du fait que les
 $X_\eta$ sont caractérisés par des propriétés de transfert endoscopique  
 et que le transfert endoscopique commute à l'induction.
Or certaines induites ne sont pas irréductibles, ceci est manifeste dans \cite{BV} et,
 comme on vient de le voir, se traduit 
par le fait que les paquets ne sont pas des singletons.
 Ces paquets ne sont pas disjoints en général, et l'on peut  même voir apparaître le cas 
où un paquet associé à un paramètre $\psi$ est entièrement contenu dans un paquet associé à un autre paramètre 
$\psi'$. Ce phénomène, dans le cas des groupes complexes, est uniquement lié au fait 
(déjà remarqué par Barbasch et Trapa) que dans les paramétrisations interviennent des orbites non spéciales, 
au sens de Lusztig-Spaltenstein. On renvoie le lecteur à \ref{intersection} pour une 
description du résultat qui est quand même assez technique.

Donnons un aperçu du contenu de cet article. La deuxième  section introduit les notations et quelques
résultats généraux sur les groupes complexes généraux et leurs représentations (c'est-à-dire en fait des modules de
Harish-Chandra). La troisième section rappelle les définitions des paramètres de Langlands et d'Arthur, et la
 paramétrisation de Langlands. On y donne aussi les propriétés conjecturales des paquets
 d'Arthur.
La quatrième section est consacrée aux groupes linéaires, qui apparaissent ici pour deux raisons.
La première mentionnée ci-dessus est que les groupes classiques considérés apparaissent comme groupes
endoscopiques tordus de groupes $\GL_N(\bbC)$. En particulier, si $\Std_G:\, {}^LG\rightarrow \GL_N(\bbC)$
est la représentation standard du groupe dual ${}^LG$, et si $\psi_G$ est un paramètre d'Arthur pour $G$, 
alors $\psi=\Std_G\circ \psi_G$ est un paramètre d'Arthur pour $\GL_N(\bbC)$ et les représentations 
du paquet  $\Pi(\psi_G,G)$ sont reliées à la représentation de $\GL_N(\bbC)$ déterminée par $\psi$ par une  
identité endoscopique (pour $\GL_N(\bbC)$, les paquets d'Arthur sont des singletons parfaitement déterminés).
La décomposition en représentations irréductibles de $\psi$ vue comme représentation de $W_\bbC \times \SL_2(\bbC)$
joue un rôle crucial dans notre étude, en particulier par des arguments de réduction
à certain types de paramètres.
La deuxième façon d'intervenir des groupes généraux linéaires est que les sous-groupes de Levi des groupes
classiques sont isomorphes à des produits d'un groupe classique de même type, et de facteurs $\GL$.

Dans la cinquième section, nous introduisons les groupes classiques considérés et nous  
 rappelons quelques résultats bien connus,  le plus important pour nous étant le lemme \ref{crucial},
donnant des renseignements sur les composantes des induites paraboliques. Ce résultat fait appel à la définition
de certains invariants des représentations que nous appelons de manière un peu abusive \og exposants \fg, 
invariants qui remplaceront les $K$-types dans certains arguments  de Barbasch-Vogan et Barbasch et les rendront
peut-être plus accessibles à certains lecteurs.
La sixième section explique comment on se ramène de l'étude des paquets d'Arthur généraux
à ceux qui  sont unipotents, et de \og bonne parité \fg. Les paquets unipotents sont ceux
dont le paramètre est trivial sur le premier facteur  de $W_\bbC\times \SL_2(\bbC)$.
Ils sont paramétrés par les orbites nilpotentes dans l'algèbre de Lie du groupe dual.
Certains arguments de la réduction utilisent des résultats sur les paquets unipotents de bonne
parité, et la démonstration devra donc attendre la section \ref{DemGrosseReduc}. Ces résultats sont 
tous tirés entièrement des travaux de  Barbasch-Vogan  \cite{BV} et surtout \cite{B1}.
  Néanmoins, comme il s'agit là de résultats importants
 et que l'article de Barbasch est beaucoup plus général que ce dont on a besoin ici, nous en  donnons  une démonstration légèrement différente
 et simplifiée.   La principale différence avec la démonstration de Barbasch
 est que nous allons utiliser des arguments d'irréductibilité basés sur l'analyse des exposants ({\sl cf.} définition 
 \ref{Exposants}) plutôt que sur les $K$-types minimaux. Notre situation est aussi moins générale, 
 ce qui permet certaines simplifications.

 Lorsque le paramètre d'Arthur $\psi_G$ est unipotent, c'est-à-dire trivial sur $W_\bbC$, 
 sa donnée est équivalente via la théorie de Jacobson-Morozov
 à la donnée d'une orbite nilpotente $\caU$ dans l'algèbre de Lie ${}^L\frg$  du groupe dual ${}^LG$. Lorsque
 $\caU$ est   une orbite nilpotente  spéciale paire, Barbasch et Vogan ont proposés
 antérieurement à Arthur, une définition d'un paquet que nous notons $\Pi_{BV}(\caU,G)$
 possédant certaines des propriétés voulues (la principale qui manque pour les identifier
 immédiatement aux paquets définis par Arthur est celle concernant l'endoscopie tordue vers $\GL_N(\bbC)$). 
 Barbasch et Vogan caractérisent les représentations dans $\Pi_{BV}(\caU,G)$ par des conditions portant sur 
  le caractère infinitésimal (il doit être donné par l'élément semi-simple d'un $\frs\frl_2$-triplet
  associé à l'orbite $\caU$) et   le front d'onde (il doit être égal à l'adhérence de l'orbite
  duale de $\caU$  au sens de Lusztig et Spaltenstein ).
  La septième section rappelle la description faite par Barbasch \cite{B1} des paquets de Barbasch-Vogan
 attachés à une telle  orbite  $\caU$ spéciale paire : tout caractère $\eta$ de $A(\caU)=A(\psi)$
 détermine une représentation irréductible $X^{BV}_\eta$ de $G$ si $\eta$ se factorise 
 par le quotient de Lusztig $\bar A(\caU)$ du groupe $A(\caU)$. Si $\eta$ ne se factorise pas de la sorte, on pose 
 $X^{BV}_\eta=0$.

 Dans la huitième section, nous énonçons des résultats  qui réduisent la description de certains
 paquets de Barbasch-Vogan à celle de paquets attachés à des groupes plus petits, ce qui permettra 
 de raisonner par récurrence dans certaines démonstrations.  
Dans la  neuvième  section, nous énonçons des propriétés des exposants des représentations dans 
les paquets de Barbasch-Vogan. Le résultat principal de cet article apparaît  dans la dixième section. 
Il affirme que les paquets unipotents définis par Arthur et  
  ceux définis par Barbasch-Vogan en \cite{BV} coïncident. Il y a deux points clefs dans la 
  démonstration de ce fait. Le premier
  est que l'une et l'autre des constructions donnent  les mêmes formules de transfert pour l'endoscopie
ordinaire. Ces formules sont constitutives de la construction des paquets chez Arthur, et c'est la proposition
12.4 chez Barbasch-Vogan. L'autre point clef est le fait que les représentations $X_\eta$ 
attachées aux caractères $\eta$ comme ci-dessus sont ici  irréductibles ou nulles, par définition chez
 Barbasch-Vogan, et par les résultats de \cite{pourhowe} pour celles définies par Arthur ({\sl cf.} 
 remarque \ref{mult1}).

Dans le onzième section, nous donnons la démonstration des résultats d'irréductibilité d'induites énoncés 
dans la section 6. Comme corollaire, nous en déduisons le fait que lorsque le caractère infinitésimal est régulier,
le paquet d'Arthur est un singleton (il est égal au paquet de Langlands qu'il contient
naturellement).

Dans la douzième section, nous donnons quelques compléments en particulier le calcul du front d'onde 
des représentations étudiées. Le résultat est étonnamment  simple; on introduit comme 
ci-dessus l'orbite nilpotente de ${}^L\frg$ définie par la restriction de $\psi$ à $\SL_2(\bbC)$ 
sans rien supposer sur la restriction de $\psi$ à $W_\bbC$. Alors le front d'onde des représentations associées à 
$\psi$ est exactement la fermeture du dual (au sens de Lusztig-Spaltenstein) de cette orbite.  Nous montrons aussi 
que les représentations dans un paquet d'Arthur sont stables par toute involution 
complexe (il y a une restriction pour les groupes $\SO_{2n}(\bbC)$ où nous nous limitons à un cas particulier 
important). Finalement nous étudions les intersections entre paquets.

Signalons  aussi que nous avons ajouté à la liste des groupes classiques
 étudiés par Arthur les groupes métaplectiques.
En effet,  ils se comportent exactement comme les groupes classiques et nous anticipons
 la généralisation à ces groupes de la théorie d'Arthur, anticipation justifiée par les
progrès faits par Wen Wei Li sur la stabilisation de la formule des
traces. Il se développe en ce moment un programme d'Arthur-Langlands très sophistiqué pour les groupes 
non linéaires tels que le groupe métaplectique (\cite{Weiss}, \cite{GG}),
mais pour  les groupes métaplectiques complexes, beaucoup résultats font partie du folklore depuis 
quelques temps déjà.  En effet le groupe métaplectique $\Mp_{2n}(\bbC)$ est une extension scindée
du groupe symplectique $\Sp_{2n}(\bbC)$, mais son $L$-groupe est $\Sp_{2n}(\bbC)$ plutôt que $\SO_{2n+1}(\bbC)$.
Pour tout ce qui concerne la théorie des représentations spécifiques de cette extension, on se ramène
donc trivialement au groupe symplectique, mais  pour tout ce qui concerne les constructions 
attachées au groupe dual,  les résultats
vont être différents - paramètres et    paquets d'Arthur  par exemple.
Les paquets d'Arthur pour le groupe métaplectique sont définis via la correspondance de Howe entre 
groupes métaplectiques et groupes orthogonaux impairs (\cite{pourhowe}). Ceci est fait en utilisant 
une construction globale et demandant simplement que les représentations
des paquets d'Arthur pour le groupe métaplectique soient la composante locale d'une forme automorphe de carré 
intégrable dont presque partout on a une représentation non ramifiée prescrite. 
Le cas des paramètres unipotents est aussi considéré
de façon implicite (c'est-à-dire sans mentionner le groupe dual) dans \cite{B1} et \cite{barbaschtransparent}, 
où il est remarqué que les résultats de \cite{BV} s'étendent
aux groupes métaplectiques. D'autre part, les résultats d'irréductibilité d'induites paraboliques
utilisés pour la réduction des paramètres généraux aux paramètres unipotents de bonne parité 
dans la section 6 sont aussi valables. Mais la façon 
ad hoc de définir les paquets par la correspondance de Howe, c'est-à-dire sans les caractériser par des propriétés 
de transfert  endoscopiques, empêche de vérifier aisément que la définition se comporte bien par induction.

\

Nous remercions J. Adams et D. Vogan pour l'aide qu'ils nous ont apportée en répondant 
rapidement et précisément à nos questions.
La première auteure remercie aussi l'ESI et en particulier J. Schwermer et S. Kudla pour les excellentes 
conditions de travail fournies au printemps 2015 lors de la période spéciale 
qu'ils y ont organisée et où ce travail a été commencé.

\section{Notations et généralités sur les groupes complexes et leurs représentations}

\subsection{Groupes complexes}\label{grC}

On note $\Gamma=\mathrm{Gal}(\bbC/\bbR)$ le groupe de Galois de $\bbC$ sur $\bbR$ et $\sigma$ 
son élément non trivial.
Soient $\mathbf{G}$ un groupe algébrique connexe réductif défini sur $\bbC$ et  $G=\mathbf{G}(\bbC)$
 le groupe de ses points 
complexes.  On fixe un sous groupe compact maximal $K$ de $G$, et on note $\sigma_c$ l'involution antiholomorphe de 
$G$ dont le groupe des points fixes est $K$.  On fixe aussi une paire de Borel $(B,H)$ de $G$, où l'on suppose 
que le tore maximal $H$ est $\sigma_c$-stable. On note $H=TA$ sa décomposition de Cartan. 
 
Soient $\frg$,  $\frb$, $\frh$, $\frt$ et $\fra$ les algèbres de Lie de $G$, $B$, $H$,  $T$ et $A$ respectivement
(les trois premières sont des algèbres de Lie complexes, les deux dernières réelles).
Leur complexification respective s'écrit :
  \begin{equation}\label{GC10} \frg_\bbC=\frg\times  \frg, \quad    \frb_\bbC=\frb \times \frb, \quad 
   \frh_\bbC=\frh\times  \frh \end{equation}
 \begin{equation}\label{GC11} 
\frt_\bbC= \{ (H,-H))\; \vert \; H \in \frh \}, 
   \quad \fra_\bbC= \{ (H,H)\; \vert \; H \in \frh \} . \end{equation}

Soient $R=R(\frg,\frh)$ le système de racines de $\frh$ dans $\frg$ et $R^+$ le système des racines
positives relativement à $\frb$.
Le système de racine $R_\bbC$ de $\frh_\bbC$ dans $\frg_\bbC$ est alors 
 \begin{equation}\label{GC12} R_\bbC=\{ (\alpha,0), \alpha \in R  \} \coprod \{ (0,\alpha), \alpha \in R  \}. 
\end{equation}
L'action de $\sigma \in \Gamma$ sur $R_\bbC$  est donnée par  $\sigma\cdot ( (\alpha,0))=(0,\alpha)$. 
Ceci montre que la sous-algèbre de Borel
$\frb_c=\frb\times \frb$ de $\frg_\bbC$ est stable sous l'action de $\sigma$. 
Le système de racines positives 
 \begin{equation}\label{GC13} R_\bbC^+=\{ (\alpha,0), \alpha \in R^ +  \} \coprod \{ (0,\alpha), \alpha \in R^+  \} 
\end{equation}
est celui défini par cette sous-algèbre de  Borel. 
On note $W$ le groupe de Weyl de $R$ et l'on identifie celui de $R_\bbC$ à $W\times W$.

La donnée radicielle associé à $G$ est le quadruplet $(X^*(H),\Delta,  X_*(H), \Delta \, \check{}\, )$, où  
$X^*(H)$ est le groupe des caractères algébriques de $H$, $X_*(H)$ celui  des cocaractères et $\Delta   \, \check{}\, $
l'ensemble des coracines simples.

\subsection{Représentations irréductibles des groupes complexes}
Nous reprenons les notations de la section précédente pour un groupe complexe $G=\mathbf{G}(\bbC)$.
Soient $\lambda,\mu \in \frh^*$ tels que $\lambda-\mu$ est un poids d'une représentation holomorphe 
de dimension finie de $G$ (c'est-à-dire $\lambda-\mu\in X^ *(H)$). Définissons un caractère $\bbC_{\lambda,\mu}$ de $H$ par 
\begin{equation}\label{GC21}
\bbC_{\lambda,\mu \vert T}= \bbC_{\lambda-\mu}, \qquad \bbC_{\lambda,\mu \vert A}= \bbC_{\lambda+\mu}.
\end{equation}
On étend $\bbC_{\lambda,\mu}$ en un caractère de $B$ et l'on pose 
\begin{equation}\label{GC22}
X(\lambda,\mu)=\Ind_B^G(  \bbC_{\lambda,\mu}). \end{equation}
L'induction est ici l'induction parabolique infinitésimale ({\sl i.e.} l'induite et l'induisante 
sont des modules de Harish-Chandra), et  normalisée,  (voir la  la section XI.2 de \cite{KV}). 
Le module $X(\lambda,\mu)$ est la série principale de paramètre $(\lambda,\mu)$
 et l'on note 
 \begin{equation}\label{GC23}
\bar X(\lambda,\mu) \end{equation}
son unique sous-quotient irréductible contenant le  $K$-type de poids extrémal $\lambda-\mu$.

On a alors le résultat de classification suivant, dû à Zhelobenko :
\begin{prop}\label{Zhe} Soient  $\lambda,\mu, \lambda',\mu' \in \frh^*$ tels que $\lambda-\mu$ et 
$\lambda'-\mu'$  soient des poids d'une représentation holomorphe 
de dimension finie de $G$. Alors il y équivalence entre
\begin{itemize}
\item[(a)] $X(\lambda,\mu)$ et $X(\lambda',\mu')$ ont mêmes facteurs de composition avec mêmes multiplicités,
\item[(b)] $\bar X(\lambda,\mu)\simeq \bar X(\lambda',\mu')$,
\item[(c)] il existe $w \in W$ tel que  $w\cdot (\lambda,\mu)=(\lambda',\mu')$.
\end{itemize}
De plus, tout $(\frg,K)$-module irréductible est équivalent à un $\bar X(\lambda,\mu)$.
\end{prop}

Le module $\bar X(\lambda,\mu)$ est le sous-quotient de Langlands de $  X(\lambda,\mu)$. Il appara\^it 
comme quotient, ou comme sous-représentation, lorsque $(\lambda,\mu)$ possède les propriétés de positivité 
ou négativité requises (par rapport à $B$) et que nous n'explicitons pas ici. On peut toujours trouver
$w\in W$ tel que $\bar X(\lambda,\mu)\simeq \bar X(w\cdot \lambda,w\cdot \mu)$ soit un quotient (resp. un 
sous-module) de $X(w\cdot \lambda,w\cdot \mu)$.

\subsection{Caractère infinitésimal}

Soit $\frU(\frg_\bbC)$ l'algèbre enveloppante de $\frg_\bbC$. Comme $\frg_\bbC=\frg\times \frg$, on a 
$\frU(\frg_\bbC)=\frU(\frg)\otimes \frU(\frg)$. Notons $\frZ(\frg_\bbC)$ le centre de cette algèbre enveloppante
et $S(\frg_\bbC)$ l'algèbre symétrique sur $\frg_\bbC$.
On a $\frZ(\frg_\bbC)=\frZ(\frg)\otimes \frZ(\frg)$,   $S(\frg_\bbC)= S(\frg)\times S(\frg)$
 et 
 \begin{equation}\label{GC24}
\gamma_{\frg_\bbC} : \;  \frZ(\frg_\bbC)=\frZ(\frg)\otimes \frZ(\frg) \longrightarrow S(\frh_\bbC)^{W\times W}=
 S(\frh)^W \otimes S(\frh)^W \end{equation}
est l'isomorphisme de Harish-Chandra. Via cet isomorphisme, les caractères de  $\frZ(\frg_\bbC)$ sont 
paramétrés par les couples $(\lambda,\mu) \in \frh^*\times \frh^*$, deux caractères, paramétrés respectivement par 
$(\lambda,\mu)$ et $(\lambda',\mu')$ étant égaux si et seulement s'ils sont conjugués par $W\times W$.

Si un   $(\frg,K)$-module $X$  admet un caractère infinitésimal  paramétré par 
 $(\lambda,\mu) \in \frh^*\times \frh^*$, on dira simplement que $X$ a pour caractère infinitésimal $(\lambda,\mu)$.
 
\begin{rmqs} 
\begin{itemize}
\item[(a)]Les choses étant bien faites, le caractère infinitésimal de $X(\lambda,\mu)$ (et donc de 
$\bar X(\lambda,\mu)$) est $(\lambda,\mu)$.
\item[(b)] Il découle de la classification de Zhelobenko que si  le  $(\frg,K)$-module $X$  
 a pour caractère infinitésimal $(\lambda,\mu)$, alors il existe $w\in W$ tel que $\lambda-w\cdot \mu\in X^ *(H)$.
\end{itemize}
 \end{rmqs}

\subsection{Un résultat sur les induites paraboliques}
On continue avec les mêmes notations que dans les sections précédentes.
\begin{lemme} \label{cruc} Soit $P=MN$ un sous-groupe parabolique de $G$ de facteur de Levi $M$ et
 de radical unipotent $N$. On suppose
que $M$ contient le tore maximal $H$. Soient  $ X^M(\lambda,\mu)$ une représentation de la série principale de 
$M$ et  $\bar X^M(\lambda,\mu)$   son sous-quotient de Langlands comme en (\ref{GC23}).

$(i)$ $\Ind_P^G( X^M(\lambda,\mu))=X(\lambda,\mu)$, la série principale de $G$ de paramètre $(\lambda,\mu)$.

$(ii)$ $\Ind_P^G( \bar X^M(\lambda,\mu))$ contient $\bar X(\lambda,\mu)$, comme sous-quotient.
  
$(iii)$ Si $\Ind_P^G( \bar X^M(\lambda,\mu))$ est réductible, alors elle contient un  sous-quotient 
 $\bar X(\lambda,s_\alpha\cdot \mu)$ où  $\alpha$ est une racine  de $H$ dans $N$
telle que  $\bil{\alpha \check{}}{\lambda}$ et  $\bil{\alpha \check{}}{\mu}$ sont des entiers non nuls 
de m\^eme signe et $s_\alpha\in W$ est la reflection par rapport à cette racine.
\end{lemme}
\dem Le premier point est immédiat par transitivité de l'induction parabolique. Le second en découle
car l'induite contient le $K$-type de poids extrémal $\lambda-\mu$. 
Le troisième point est dû à Zhelobenko \cite{Zhe}. \qed

\section{Paramètres de Langlands et d'Arthur}
\subsection{Paramètre de Langlands} \label{LangL} Le groupe de  Weil de  $\bbC$ est $W_\bbC= \bbC^\times$.
\bigskip 

Soit $\mathbf{G}$ un groupe algébrique réductif connexe défini sur $\bbC$,  et l'on adopte les mêmes
notations qu'en  \ref{grC}. Soit ${}^LG$
son  dual de Langlands. Il s'agit ici du groupe complexe connexe déterminé par la donnée radicielle duale, 
c'est-à-dire que l'on suppose fixé une paire de Borel $(\caB,\caH)$ de ${}^LG$ avec les  identifications
$X_*(\caH)=X^*(H)$, $X^*(\caH)=X_*(H)$ .

\begin{defi} \label{DefParLang}
Un paramètre  de Langlands est un morphisme continu: 
\[ \phi: \; W_\bbC \rightarrow  {}^LG  \]
tel que   $\phi$ a  pour image des   éléments  semi-simples de  ${}^LG$.

Le groupe  ${}^LG$ agit par conjugaison sur l'ensemble des paramètres de Langlands, et l'on note 
$\Phi(G)$ l'ensemble de ces classes de conjugaison.
\end{defi}

\bigskip 
Nous commettrons fréquemment l'abus de langage consistant à ne pas distinguer
entre un paramètre de Langlands et l'élément de $\Phi(G)$ qu'il définit.  

\medskip 

Soit $\phi: \; W_\bbC \rightarrow  {}^LG  $ un paramètre de Langlands. A conjugaison près, on peut supposer que 
l'image de $\phi$ est contenue dans le tore $\caH$. On peut donc écrire
 \begin{equation}\label{GC31}\phi(z)=z^\lambda \bar z^\mu, \quad (z\in \bbC^ \times). \end{equation}
où $\lambda,\mu$ sont dans $X_*(\caH)\otimes_\bbZ\bbC =X^*(H)\otimes_\bbZ\bbC=\frh^*$ et $\lambda-\mu \in X^*(H)$.
Notons $\phi=\phi_{\lambda,\mu}$ le paramètre de Langlands défini par (\ref{GC31}).

Le résultat suivant, élémentaire, est le pendant de la proposition \ref{Zhe}.
\begin{prop} Soient  $\lambda,\mu, \lambda',\mu' \in \frh^*$ tels que $\lambda-\mu$ et 
$\lambda'-\mu'$ soient dans  $X^*(H)$. Alors $\phi=\phi_{\lambda,\mu}$ et $\phi=\phi_{\lambda',\mu'}$
sont équivalents si et seulement s'il existe $w \in W$ tel que  $w\cdot (\lambda,\mu)=(\lambda',\mu')$.
De plus tout paramètre de Langlands $\phi:  \; W_\bbC \rightarrow  {}^LG  $
est conjugué à un $\phi_{\lambda,\mu}$.
\end{prop}

On peut donc reformuler la classification des  $(\frg,K)$-modules irréductibles de la proposition \ref{Zhe}.
\begin{cor}
La correspondence $\phi_{\lambda,\mu}\leftrightarrow \bar X(\lambda,\mu)$ induit une bijection entre $\Phi(G)$ et 
l'ensemble des classes d'isomorphie de $(\frg,K)$-modules irréductibles. 
\end{cor}

\begin{rmq}
Les paquets de Langlands pour les groupes complexes sont des singletons. 
\end{rmq}

\subsection{Paramètres d'Arthur}\label{ArtP}
Les notations sont les mêmes que dans le section précédente.
\begin{defi}\label{ArtPar}
Un paramètre d'Arthur pour  $G$ est un morphisme de groupes continu
\begin{equation*}
\psi: \, W_\bbC \times \SL_2(\bbC) \longrightarrow {}^ L G
\end{equation*}
tel que 
\begin{itemize}
\item[(i)] la restriction de $\psi$ à  $W_\bbC$ est un paramètre de Langlands tempéré,
\item[(ii)]  la restriction de  $\psi$ à $\SL_2(\bbC)$ est algébrique.
\end{itemize}

Dans (i), rappelons que le paramètre de Langlands $\psi$ est dit tempéré s'il est d'image bornée.
Avec les notations de (\ref{GC31}), ceci est équivalent au fait que $\lambda+\mu\in X^*(H)\otimes_\bbZ i\bbR$. 
\medskip

Le groupe ${}^ L G$ agit par conjugaison sur l'ensemble des paramètres d'Arthur, et l'on note 
$\Psi(G)$ l'ensemble de ces classes de conjugaison.
\end{defi}

A tout paramètre d'Arthur  $\psi$, on associe un paramètre de Langlands 
\begin{equation}\label{ArtLang}
\phi_\psi: \, W_\bbC \longrightarrow {}^ L G, \qquad z \mapsto \psi(z,
  \left(  \begin{matrix}  (z\bar z)^{\frac{1}{2}} & 0\\ 
0 &  (z\bar z)^{-\frac{1}{2}}   \end{matrix} \right)) . 
\end{equation}

\bigskip 

Soit $\psi: \, W_\bbC \times \SL_2(\bbC) \longrightarrow {}^ L G$ un paramètre d'Arthur.
 A conjugaison près, on peut supposer que  sa restriction à $\bbC^\times$ est donnée par  : 
 \begin{equation}\label{GC41}\psi(z)=z^\lambda \bar z^\mu, \quad (z\in \bbC^ \times). \end{equation}
avec $\lambda-\mu \in X^*(H)$. La condition de tempérance $(i)$ dans la définition des paramètres d'Arthur
nous donne en plus que $\lambda+\mu \in \frt^*$. Ici, on voit $\frt^*$ comme le sous-espace réel de $\frh^*$ 
des formes linéaires qui prennent des valeurs réelles sur $\frt$. De même,  
on voit $\fra^ *$ comme le sous-espace réel de $\frh^*$ 
des formes linéaires qui prennent des valeurs réelles sur $\fra$, de sorte que $\frh^*=\frt^*\oplus \fra^*$
avec $i\fra^ *=\frt^*$.

\bigskip

La restriction de $\psi$ à  $\SL_2(\bbC)$ est un morphisme algébrique. Notons 
 \begin{equation}\label{GC42}\bar \psi :\, \frs\frl_2(\bbC) \longrightarrow  {}^L\frg \end{equation}
 sa différentielle en l'identité ($ {}^L\frg$ est bien s\^ur l'algèbre de Lie de  ${}^ L G$), et 

 \begin{equation}\label{GC43}h_\psi= \bar \psi (  \left(  \begin{matrix}  1 & 0\\ 
0 &  -1   \end{matrix} \right)), \qquad e_\psi= \bar \psi (  \left(  \begin{matrix}  0 & 1\\ 
0 &  0   \end{matrix} \right)). \end{equation}

L'orbite (nilpotente) de $e_\psi$    dans $ {}^L\frg$  sous l'action adjointe de   ${}^ L G$ est notée $\caU_\psi$.
Comme l'image de $\SL_2(\bbC)$ et celle de $\bbC^\times$ par $\psi$ commutent dans ${}^ L G$, on peut aussi supposer 
à conjugaison près que $h_\psi$ est dans l'algèbre de Lie de $\caH$ (que l'on identifie à $\frh^*)$).
On obtient alors que $\phi_\psi$ est donné par 
\begin{equation}\label{PP}
\phi_\psi(z)=z^ {\lambda+\frac{1}{2}h_\psi}\bar z^ {\mu+\frac{1}{2}h_\psi} . 
\end{equation}
Le caractère infinitésimal attaché à $\psi$ est :
\begin{equation}\label{GC44}(\lambda+\frac{1}{2}h_\psi, \mu+\frac{1}{2}h_\psi).\end{equation}
En particulier, les éléments de ce paquet sont de la forme 
\begin{equation}\label{GC45} \bar X ( \lambda+ \frac{1}{2}h_\psi, w\cdot   (\mu+ \frac{1}{2}h_\psi) ) .\end{equation}
pour des éléments $w \in W$ vérifiant 
\begin{equation}\label{GC46}  \lambda+ \frac{1}{2}h_\psi- w\cdot   (\mu+ \frac{1}{2}h_\psi) \in X^ *(H) .\end{equation}

\subsection{Paquets d'Arthur}\label{pacAr}
Dans \cite{Art84}, \cite{Art89}, J. Arthur conjecture l'existence de paquets $\Pi(\psi, G)$  attachés aux 
paramètres $\psi \in \Psi(G)$,   devant posséder certaines propriétés. Parmi les principales,
 citons le fait que les $\Pi(\psi,G)$ sont finis, constitués de (classes d'équivalence de)  représentations
{\sl  unitaires}, ayant toutes le m\^eme caractère infinitésimal, donné par (\ref{GC44}).
 Le paquet d'Arthur $\Pi(\psi,G)$ contient
 le paquet de Langlands $\Pi(\phi_\psi, G)$ (c'est-à-dire la représentation 
 $\bar X( \lambda+ \frac{1}{2}h_\psi,\mu+ \frac{1}{2}h_\psi)$). 
 Ils doivent  satisfaire les identités de caractères  attendues
dans la théorie de l'endoscopie (standard et tordue); c'est ce qui est appelé le {\sl transfert spectral}.
En revanche, ces  paquets ne sont pas disjoints, et ne sont pas des réunions de $L$-paquets. 
Comme nous l'avons déjà remarqué, pour les groupes complexes la notion de conjugaison stable est triviale 
(c'est la conjugaison ordinaire), et ne donne donc aucune contrainte sur les paquets.

Arthur énonce ses conjectures pour des groupes définis
sur un corps local $F$ quelconque (remplacer $W_\bbC$ par le groupe de Weil-Deligne $W_F$). 
Dans \cite{ABV},  pour un corps $F$ local archimédien, 
des paquets  $\Pi^{ABV}(\psi,G)$ possédant les propriétés voulues, à l'exception, malheureusement, de la compatibilité à l'endoscopie tordue
 sont définis par des méthodes géométriques puissantes. C'est une conjecture raisonnable de dire que les constructions de \cite{ABV} coïncident avec celle de \cite{Art13} dans le cas \og particulier\fg\,  de \cite{Art13}.  
Antérieurement, pour $F\simeq \bbC$, Barbasch et Vogan ont défini dans \cite{BV} des paquets attachés aux
 paramètres unipotents (c'est à dire ceux dont la restriction à $\bbC^\times$ est triviale) attachés 
 par la théorie de Jacobson-Morosov à une orbite nilpotente spéciale paire. Ensemblistement les constructions sont, dans ce cas particulier des paramètres unipotents, les mêmes d'après \cite{ABV} chapitre 27. Nous  reviendrons sur la description de ces constructions qui sont fondamentales pour ce que nous faisons ici.

Dans \cite{Art13}, J. Arthur donne une définition des paquets $\Pi(\psi,G)$ lorsque $G$ est un groupe classique.
Donnons quelques précisions au sujet de ces paquets. 
Soit $\psi$ un paramètre d'Arthur pour le groupe $G$. Soit $S_\psi$ le centralisateur de l'image de $\psi$
dans ${}^LG$ et $S_\psi^0$ sa composante connexe neutre.
On pose 
\begin{equation} \label{Apsi} A(\psi)=S_\psi/S_\psi^0Z({}^LG).\end{equation}
Pour les groupes complexes classiques, les groupes $A(\psi)$ sont abéliens, ce sont des produits de facteurs
$\bbZ/2\bbZ$. Nous ne considérerons que ce cas dans la suite. Notons $\widehat{A(\psi)}$ le groupe des caractères 
de $A(\psi)$. 
Arthur  définit une application :
\begin{equation}\label{Xeta}
\eta \in \widehat{A(\psi)} \mapsto X_\eta 
\end{equation}
 où $X_\eta$ est représentation semi-simple de longueur
finie   de $G$ (éventuellement  avec des multiplicités), ou bien  $\{0\}$. Ces représentations sont uniquement définies par les relations suivantes.

Pour tout $s \in A(\psi)$, on considère la représentation virtuelle
\begin{align}\label{distrpasstable}
X_\psi^s= \sum_{\eta \in \widehat{A(\psi)} }  \eta(s_\psi s) \;   X_\eta 
\end{align}
où $s_\psi$ est l'image dans $A(\psi)$ de l'élément $\psi((1,-\Id))$, $(1,-\Id)\in W_\bbC\times \SL_2(\bbC)$.
 Ce sont ces représentations virtuelles qui apparaissent dans les
 identités de transfert endoscopiques qui caractérisent les $X_\eta$.
Lorsque $s=1$, (\ref{distrpasstable}) est la représentation stable attachée au paquet $ \psi$. 
On la note  $X^ {st}_\psi$ : 
\begin{align}\label{distrstable}
X_\psi^{st}= \sum_{\eta \in \widehat{A(\psi)} }  \eta(s_\psi ) \;   X_\eta 
\end{align}

Le paquet  $\Pi(\psi,G)$ est alors l'ensemble des représentations irréductibles de $G$ qui apparaissent dans les 
$X_\eta$, lorsque $\eta$ décrit $\widehat{A(\psi)}$. Et à chaque élément de $\Pi(\psi,G)$ est associée une représentation
du groupe $A(\psi)$:
 
\begin{equation}\label{rhoAX}
\bar X\in \Pi(\psi,G)\mapsto \rho^A_{\bar X}
\end{equation}
qui attache a tout élément d'un paquet   $\Pi(\psi,G)$ une représentation de dimension finie
du groupe $A(\psi)$. Notons ceci $\bar X \mapsto \rho^A_{\bar X}$, de sorte que pour tout $s\in A(\psi)$, 
\begin{equation}\label{rhoeta} \Tr(\rho^A_{\bar X}(s))=\sum_{\eta \in \widehat{A(\psi)} } m(\bar X,X_\eta)\;  \eta(s),
 \end{equation}
 où $m(\bar X,X_\eta)$ désigne la multiplicité de $\bar X$ dans $X_\eta$. On peut alors réécrire les représentations virtuelles 
 $X_\psi^s$ sous la forme 
\begin{align}\label{distrstable2}
X_\psi^s= \sum_{\bar X \in \Pi(\psi,G) }  \Tr(\rho^ A_{\bar X}(s_\psi s)) \;   \bar X. 
\end{align}

\begin{rmq}\label{mult1}
Pour les groupes classiques complexes ou réels, et les paramètres $\psi$ unipotents, 
il est  démontré dans \cite{pourhowe} que les multiplicités $m(\bar X, X_\eta)$ sont 
$0$ ou $1$ et plus généralement que les représentations $\rho^ A_{\bar X}$ sont irréductibles.
Ce qui  est propre au cas des groupes complexes, c'est le fait 
aussi démontré dans {\sl loc. cit.} que pour tout $\eta \in \widehat{A(\psi)}$, 
les  $X_\eta$ définis ci-dessus sont des représentations  irréductibles ou nulles, c'est-à-dire que 
l'application qui à $\bar X\in \Pi(\psi,G)$ associe le caractère $\rho^A_{\bar X}$ est injective. 
Cette dernière propriété n'est pas vraie pour les groupes classiques réels.

\end{rmq}
La première partie de la remarque est un des points clé de la démonstration du théorème \ref{egalitepaq}
 et dans cet article on  généralise toute la remarque à tout paramètre $\psi$.

\section{$\GL_N$}

\subsection{$\GL_N(\bbC)$}\label{GLNPar1}

 Soit $N$ un entier positif.  Dans cette section, on s'intéresse au cas du groupe algébrique  $\mathbf{G}=\GL_N$
 défini sur $\bbC$. Avec les notations de la section \ref{grC}, on prend  $\sigma_c : \, g\mapsto {}^t\bar g^{-1}$ 
 et on fixe  l'épinglage usuel $\mathbf{spl}= (\mathbf{B}_d,\mathbf{H}_d, \{X_\alpha\}_\alpha)$ 
 où $ \mathbf{B}_d$ est le sous-groupe de Borel des matrices triangulaires supérieures,
  $\mathbf{H}_d$ est le tore  diagonal et $X_\alpha$ est un vecteur radiciel pour une racine simple $\alpha$
 du système de racines positives  de $\mathbf{H}_d$ dans $\GL_N$.
 Les sous-groupes de Levi et les sous-groupes paraboliques standard de $G_N$ sont définis relativement à 
$(\mathbf{H}_d,\mathbf{B}_d)$.
\medskip
 
Via l'isomorphisme d'Harish-Chandra, un caractère infinitésimal pour $\GL_N(\bbC)$ est donné
 par un  couple $(\lambda,\mu)$ d'élément de $\frh_d^*$. En identifiant naturellement  $\frh_d$, la sous-algèbre
des matrices diagonales de $\caM_N(\bbC)$  à $\bbC^N$, et de même pour son dual $\frh_d^*$, 
un caractère infinitésimal est alors donné par un élément 
\begin{equation}\label{carinfGL}
(\lambda,\mu)=((\lambda_1,\ldots, \lambda_N), (\mu_1,\ldots,\mu_N)) \in \bbC^N \times \bbC^N,
\end{equation}
 où plutôt
 par une orbite de tels éléments sous l'action du groupe de Weyl, ici identifié au groupe $\frS_N \times \frS_N$. 
 Un tel caractère infinitésimal est entier si  les $\lambda_i-\lambda_j$ et les $\mu_i-\mu_j$  sont entiers, et régulier
 si les $\lambda_i$ sont distincts, ainsi que les $\mu_i$.

\medskip

Soit   $ N_1,\ldots ,N_r \in \bbN^\times $   tels que $\displaystyle \sum_{i=1}^r N_i=N$. 
Le sous-groupe $M=M_{N_1,\ldots ,N_r}$ des matrices diagonales par blocs de taille respective
$N_1,\ldots, N_r$, isomorphe à 
$\GL_{N_1}(\bbC) \times \GL_{N_2}(\bbC)\times ...\times \GL_{N_r}(\bbC)$
est un sous-groupe de Levi standard de $\GL_N(\bbC)$, et le sous-groupe parabolique 
 $P=P_{N_1,\ldots, N_r}$  contenant $M$ et le sous-groupe de  Borel des  matrices triangulaires supérieures 
 est un sous-groupe parabolique standard de radical unipotent  $N=N_{N_1,\ldots, N_r}$.
Pour tout  $1\leq i\leq r$, soit   $X_i$ un module de Harish-Chandra de longueur finie de $\GL_{N_i}(\bbC)$. 
On note  alors  
\begin{equation}\label{prodstarGL
}X_1\star X_2\star \cdots \star X_r \end{equation}
la représentation obtenue par induction parabolique 
(normalisée) à partir de la représentation $X_1\otimes X_2\otimes ...\otimes X_r$ de $M$ 
relativement au sous-groupe parabolique $P$. 

\subsection{Paramètres et paquets d'Arthur pour $\GL_N$}\label{ArtGN}

Un paramètre de Langlands de $\GL_N(\bbC)$ est  un morphisme continu :
\[ \phi: \; W_\bbC \longrightarrow {}^L\GL_N(\bbC)=\GL_N(\bbC),\]
 c'est-à-dire  une représentation de dimension $N$ de $\bbC^\times$.
Le fait que les  $\phi(z)$, $z\in \bbC^\times$  soient semi-simples nous dit que cette 
représentation est complètement réductible. 

Les  représentations irréductibles de  $W_\bbC\simeq \bbC^\times$
sont de dimension $1$ puisque  $\bbC^\times$ est
 abélien. Elles sont paramétrées par les couples 
$(\lambda, \mu)\in \bbC\times \bbC$ avec  $\lambda-\mu\in \bbZ$ de la manière suivante : 
\begin{equation}\label{Chisn}\chi_{\lambda,\mu}(z)=
z^\lambda\bar z^\mu
\end{equation}
Le caractère   $\chi_{\lambda,\mu}(z) $ est unitaire si  $\lambda+\mu \in i\bbR$.

Un paramètre d'Arthur pour $\GL_N(\bbC)$ est  un morphisme continu :
\[ \psi: \; W_\bbC \times \SL_2(\bbC)\longrightarrow {}^L\GL_N(\bbC) =\GL_N(\bbC)\]
vérifiant les propriétés énoncées dans la définition \ref{ArtPar}.
Cette représentation de dimension $N$ de $W_\bbC \times \SL_2(\bbC)$ dans $\bbC^ N$ est complètement réductible. Elle 
s'écrit donc comme une somme directe 
\begin{equation}
\label{opsii} \psi=\oplus_{i=1,\ldots,r}\;  \psi_i, \qquad \psi_i :\; W_\bbC\times \SL_2(\bbC) \rightarrow \GL_{N_i}(\bbC)
  \end{equation}
avec $\psi_i$ irréductible et $\sum_{i=1}^r N_i=N$. Les  représentations irréductibles
de  $W_\bbC\times \SL_2(\bbC) $ sont des produits tensoriels de représentations irréductibles de $W_\bbC$
avec des représentations irréductibles de  $\SL_2(\bbC)$. Les représentations irréductibles
de  $W_\bbC$ sont les caractères $\chi_{\lambda,\mu}$ décrits ci-dessus. 
 Les représentations irréductibles algébriques de $\SL_2(\bbC)$ sont déterminées par leur dimension,
  et l'on note $R_n$ un choix de représentation irréductible 
 de dimension $n$ de $\SL_2(\bbC)$ (ou sa classe d'équivalence). Les représentations irréductibles de 
  $W_\bbC\times \SL_2(\bbC)$ sont donc à équivalence près les 
  \begin{equation}
 \chi_{\lambda,\mu}\boxtimes R_n, \, \lambda-\mu \in \bbZ, n\in \bbN.
  \end{equation}

Celles qui apparaissent dans la décomposition d'un paramètre d'Arthur  ont en plus la propriété d'être à image
bornée,  ce sont donc les $\chi_{\lambda,\mu}\boxtimes R_n$ avec $\lambda+\mu \in i\bbR$.

Comme les paquets de Langlands, les paquets d'Arthur pour $\GL_N(\bbC)$ sont des singletons, et si $\phi_\psi$ 
est le paramètre de Langlands  (\ref{ArtLang}), on a 
donc $\Pi(\psi)=\Pi(\phi_\psi)$.  On note $X^\GL_\psi$  l'unique représentation qu'il contient. 
 Nous allons maintenant déterminer la représentation $X^{\GL}_\psi$ associée à un paramètre d'Arthur 
 $\psi$, en commençant par les $\psi$ irréductibles :
  \begin{equation}\label{PsiSpeh}
\text{si } \psi=\chi_{\lambda,\mu} \boxtimes R_n, \quad  X^{\GL}_\psi= \chi_{\lambda,\mu} \circ  \textstyle \det_n 
  \end{equation} 
  où $\det_n :\, \GL_n(\bbC) \rightarrow \bbC^ \times$ est le déterminant.
 Pour le cas général, on a

\begin{prop}
 Si $\psi=\oplus_{i=1,\ldots,r} \psi_i$
est une décomposition en irréductibles, alors
\begin{equation}\label{PsiSpeh2}
 X_\psi^\GL= \star_{i} \; X_{\psi_i}^\GL. 
\end{equation}
\end{prop}

\begin{rmq} 
Un résultat de Vogan \cite{VogGL} (voir aussi  \cite{Tadic} et  \cite{Baruch}) affirme que cette représentation 
est unitaire et irréductible, en particulier, 
elle ne dépend pas de l'ordre dans lequel on prend le produit.
\end{rmq}

\begin{rmq} Pour palier à l'absence de preuve de la conjecture de Ramanujan, nous sommes aussi obligés de considérer les paramètres qui sont presque unitaires, c'est-à-dire ceux pour lesquels $\Re e(\lambda+\mu)\in ]-1/2,1/2[$. L'extension à ces paramètres est sans difficulté.
\end{rmq}

On note $\theta=\theta_N$ l'automorphisme $g\mapsto {}^tg^{-1}$ de $\GL_N(\bbC)$. 
Les représentations irréductibles de $\GL_N(\bbC)$  auto-duales
(isomorphes à leur contragrégientes) sont celles qui sont stables sous l'action de $\theta_N$. Elles
jouent un rôle fondamental dans les travaux d'Arthur \cite{Art13} de par leur lien avec les représentations
 des groupes classiques. Les paramètres de Langlands ou d'Arthur qui leur sont associés
 sont ceux dont la composition avec $\theta$ (du coté dual donc) leur est conjuguée.

Le lecteur vérifiera facilement que les paramètres d'Arthur $\theta$-stables admettent une décomposition du type
\begin{equation}\label{decirrpsi}
\psi= \bigoplus_i \big( \chi_i\boxtimes R_{N_i} \oplus \chi_i^{-1}\otimes R_{N_i} \big)
 \bigoplus_j (\Triv\boxtimes R_{N_j}) 
\end{equation}
où les caractères $\chi_i$ de $\bbC^\times$ sont tels que $\chi_i\neq \chi_i^{-1}$.

\section{Les  groupes classiques et leurs représentations. Paquets d'Arthur}

\subsection{Les  groupes classiques}\label{grcla}
  
  Les  groupes classiques complexes que nous considérons sont  les groupes de rang $n$ suivants: 
 
 \medskip
 
 \begin{itemize}
 \item[$\mathbf C_n$.] Le groupe symplectique $\Sp_{2n}(\bbC)$.  
  Son dual de Langlands est $\SO_{2n+1}(\bbC)$. 
 
  \item[$\mathbf B_n$.] Le groupe spécial orthogonal impair $\SO_{2n+1}(\bbC)$. 
  Son dual de Langlands   est $\Sp_{2n}(\bbC)$.
   
    \item[$\mathbf D_n$.]  Le groupe spécial orthogonal pair $\SO_{2n}(\bbC)$. 
   Son dual de Langlands  est $\SO_{2n}(\bbC)$.

\medskip 
A cette liste de groupes classiques, nous  ajoutons un cas  un peu moins classique pour lequel nous disposons
de résultats incomplets (ils n'entrent pas dans le cadre des travaux d'Arthur), mais néanmoins intéressants.

\medskip 

   \item[$\mathbf{Mp}_n$.]  Le groupe métaplectique $\Mp_{2n}(\bbC)$. 
   Son dual de Langlands  est $\Sp_{2n}(\bbC)$ ({\sl cf.} \cite{Weiss}).

    \end{itemize}

\begin{rmq}
Le groupe métaplectique  $\Mp_{2n}(\bbC)$ est une extension scindée du groupe
 symplectique $\Sp_{2}(\bbC)$,  mais son  dual de  Langlands est différent de celui-ci. 
 Du point de vue de la théorie des représentations spécifiques,   c'est donc le même groupe, mais pas du point 
 de vue de la théorie d'Arthur-Langlands.
\end{rmq}

Nous notons ces groupes $G$, ou bien $G_n$ lorsqu'on veut garder l'information sur le rang.
On fixe une paire de Borel $(B,H)$ pour chacun de ces groupes comme dans la section \ref{grC}.
Soit  $(\epsilon_1, \ldots , \epsilon_n)$ la base de $\frh^*$ telle que le système de racines de $H$ dans $G$
soit constitué des racines $\pm \epsilon_i\pm\epsilon_j$, $1\leq i<j\leq n$, auxquelles
 on ajoute les racines $\pm 2\epsilon_i$
dans les cas  {$\mathbf C_n$} 
et  $\mathbf{Mp}_n$,  
et $\pm \epsilon_i$ dans le cas {$\mathbf B_n$}, $1\leq i\leq n$. 
Le choix de cette base identifie 
$\frh^*$ à $\bbC^n$ et l'on note $\lambda=(\lambda_1,\ldots, \lambda_n)$ un élément de $\frh^*$.

Soit $\bar X=\bar X(\lambda,\mu)$ une représentation irréductible de $G$.
Les paramètres $\lambda=(\lambda_1,\ldots ,\lambda_n)$ et $\mu=(\mu_1,\ldots , \mu_n)$ 
sont des éléments de $\frh^*$  que l'on a identifié à $\bbC^n$. On a alors pour tout $i=1,\ldots,n$, $\lambda_i-\mu_i\in \bbZ$
et $\chi_{\lambda_i,\mu_i}(z)=z^{\lambda_i}\bar z^{\mu_i}$ est un caractère de $\bbC^\times$.

\subsection{Exposants}

Dans ce paragraphe, on introduit la notion d'exposants (la terminologie est peut-\^etre abusive) d'une 
représentation irréductible d'un groupe classique complexe.

\begin{defi}\label{Exposants}
On appelle exposants de $\bar X=\bar X(\lambda,\mu)$  l'ensemble  (avec multiplicités) des caractères 
$\chi_{\lambda_i,\mu_i}$ de $\bbC^\times$ comme ci-dessus
et on le note $\mathrm{Exp}(\bar X)$.
\end{defi}

\bigskip 

Soit $P=MN$ un sous-groupe parabolique de $G$. On suppose que le facteur de Levi $M$ contient
$H$ et que $N$ est contenu dans $B$ ($P$ est \og standard\fg\, relativement au choix de la paire de Borel
$(B,H)$). Le facteur de Levi $M$ est isomorphe à un produit
\begin{equation}\label{M}  \left( \times_{i=1}^k \GL_{n_i} (\bbC)  \right)\times G' \end{equation}
où $G'$ est un groupe classique de rang $n'$, avec $\sum_{i=1}^k n_i+n'=n$.
En choisissant des réalisations explicites de nos groupes classiques comme sous-groupes de groupes généraux linéaires,
par exemple comme dans \cite{Art13}, \S1.2, 
il est possible de fixer l'isomorphisme entre $M$ et (\ref{M}), de sorte qu'on
va les identifier dans la suite  sans plus de précautions. On a bien conscience que dans le cas où $G$ est 
un groupe orthogonal pair, si $n'=0$ la classe d'isomorphie du sous-groupe de Levi ne détermine une 
classe d'association de parabolique qu'à conjugaison près sous l'action du groupe orthogonal.

Soit $ X'$ une représentation de $G'$, et pour tout $i=1,\ldots, k$, soit $X_i$ une représentation de $\GL_{n_i}
(\bbC)$.
On note alors 
\begin{equation}\label{defstar}  \left(\star_{i=1}^k X_i \right)\star X'=\Ind_P^G(\left(\otimes_{i=1}^k X_i\right) 
\otimes X'). 
\end{equation}

Le résultat suivant sera crucial dans la suite. Nous adoptons le raccourci de langage suivant : on dit qu'un 
réel $a$ est demi-entier s'il appartient à $\frac{1}{2}\bbZ \setminus \bbZ$.

\begin{lemme} \label{crucial} 
Soit $\bar X$ une représentation irréductible d'un groupe classique $G_n$. Soit $\chi=\chi_{a,b}$ un caractère de $\bbC^\times$
avec $a$ et $b$ réels, et supposons que 
\[Y= \left( \chi_{a,b}\circ \textstyle  \det_r \right)\star \bar X \]
soit réductible. Alors l'un des sous-quotients de $ Y$ possède un exposant $\chi_{a',b'}$ avec 
$\vert a'-b'\vert >\vert a-b\vert$.
D'autre part, si   $\frac{r-1}{2}+a$
est demi-entier (resp. entier), alors on a la même conclusion avec de plus $a',b'$ demi-entiers (resp. entiers).   
\end{lemme}

\dem On applique le lemme \ref{cruc} à $G^+=G_{n+r}$,  son sous-groupe de Levi standard $M=\GL_r(\bbC)\times G_n$ 
et la représentation $\bar X^M=(\chi_{a,b}\circ \textstyle  \det_r)\otimes \bar X$ de celui-ci.
Si $\bar X=\bar X(\lambda,\mu)$, avec $\lambda=(\lambda_1,\ldots,\lambda_n)$, $\mu=(\mu_1,\ldots, \mu_n)$, alors 
(en conjuguant au besoin la situation par un élément du groupe de Weyl de $G^+$)
$\bar X^M=\bar X^M(\lambda^+,\mu^+)$, avec 
\[ \lambda^+=(\ a+\frac{r-1}{2}, a+\frac{r-3}{2}, \ldots , a+\frac{1-r}{2}, \lambda_1,\ldots,\lambda_n)\]
\[ \mu^+=( b+\frac{r-1}{2}, b+\frac{r-3}{2}, \ldots , b+\frac{1-r}{2}, \mu_1,\ldots,\mu_n).\]
Le lemme \ref{cruc} nous dit qu'il existe un racine $\alpha$ de $H^+$ dans le radical unipotent du sous-groupe 
parabolique standard $P=MN$ telle que $\bar X (\lambda^+,s_\alpha\cdot \mu^+)$
soit un sous-quotient de $Y$ et $\bil{\alpha\check}{ \lambda^+}$ et $\bil{\alpha\check}{ \mu^+}$ 
sont des entiers non nuls  de même signe. Supposons que $\alpha=\epsilon_j-\epsilon_{i+r}$, avec
 $1\leq j \leq r$ et $1\leq i\leq n$
et que le signe en question soit négatif.   Ceci signifie que $\lambda_i$ et $\mu_i$ sont réels et
\[ \lambda_i>  a+\frac{r-(2j-1)}{2},  \quad  \mu_i>  b+\frac{r-( 2j-1)}{2}.\]
Or $\bar X (\lambda^+,s_\alpha\cdot \mu^+)$ a pour exposants 
$\chi_{a',b'}=\chi_{\lambda_i,  b+\frac{r-(2j-1)}{2} }$ et $ \chi_{a'',b''}=\chi_{ a+\frac{r-(2j-1)}{2}, \mu_i}$, et 
\[a'-b'=\lambda_i-\left( b+\frac{r-(2j-1)}{2}\right)>   a+\frac{r-(2j-1)}{2}-\left( b+\frac{r-(2j-1)}{2}\right)=a-b, \]
\[ a''-b''=a+\frac{r-(2j-1)}{2}- \mu_i < a+\frac{r-(2j-1)}{2}-\left( b+\frac{r-(2j-1)}{2}\right)=a-b.  \]
On a donc soit $\vert a'-b'\vert >\vert a-b\vert$, soit $\vert a''-b''\vert >\vert a-b\vert$.
 Les autres  cas ($\alpha=\epsilon_i+\epsilon_{j}$ ou bien 
 $\bil{\alpha\check}{ \lambda^+}$ et $\bil{\alpha\check}{ \mu^+}$ positifs) se démontrent de la même
manière.

Si $\alpha$ est  une racine multiple de  $\epsilon_{j}$ pour un 
$1\leq j\leq r$ et   $\bil{\alpha\check}{ \lambda^+}$ et $\bil{\alpha\check}{ \mu^+}$ sont des entiers  positifs 
(l'autre cas se  traite de la m\^eme manière), ceci donne  $a+\frac{r-(2j-1)}{2}>0$  et  $b+\frac{r-(2j-1)}{2}>0$. 
Or $\bar X (\lambda^+,s_\alpha\cdot \mu^+)$ a pour exposant 
$\chi_{a',b'}=\chi_{a+\frac{r-(2j-1)}{2},-\left(  b+\frac{r-(2j-1)}{2}\right) }$ 
et 
\[(a'-b')-(a-b)=2b+r-(2j-1)>0, \qquad (a'-b')-(b-a)=2a+r-(2j-1)>0\]
donc $(a'-b')>\vert a-b \vert $.

Si $a+\frac{r-1}{2}$  est  demi-entier (resp. entier), dans le premier cas ci-dessus
$\lambda_i$ et $\mu_i$ sont aussi demi-entiers (resp. entiers) car 
$\bil{\alpha\check}{ \lambda^+}$ et $\bil{\alpha\check}{ \mu^+}$ sont entiers, 
 et il en est de même pour    $a'$, $b'$, $a''$ et $b''$ ainsi que dans le second cas pour $a'$ et  $b'$. \qed

\section{Réduction au cas unipotent de bonne parité}

\subsection{Généralités}
Pour chacun des groupes classiques de \S \ref{grcla},
 on dispose d'une représentation naturelle du $L$-groupe dans un ${}^L \GL_N=\GL_N(\bbC)$ : 
\begin{equation}\label{DefStdG} \Std_G: \; {}^LG \longrightarrow  \GL_N(\bbC). \end{equation}
On a $N=2n+1$, $2n$, $2n$, $2n$ respectivement dans les cas
 {$\mathbf C_n$}, {$\mathbf B_n$}, {$\mathbf D_n$}, {$\mathbf{Mp}_n$}.
 Soit $G$ l'un de ces  groupes classiques  et soit
   \[\psi_G :\; W_\bbC\times \SL_2(\bbC)\longrightarrow {}^L G\] 
     un paramètre d'Arthur  pour $G$.  
     Posons  $\psi=\Std_G \circ \psi_G$. C'est un paramètre d'Arthur $\theta_N$-stable
     de $\GL_N(\bbC)$, il admet donc une décomposition en irréductibles de la forme
      \begin{equation} \label{decPsi}
      \psi= \bigoplus_i \big( \chi_i\boxtimes R_{N_i} \oplus \chi_i^{-1}\boxtimes R_{N_i} \big)
 \bigoplus_j (\Triv\boxtimes R_{N_j}) 
\end{equation}
où les caractères unitaires  $\chi_i$ de $\bbC^\times$ sont tels que $\chi_i\neq \chi_i^{-1}$.

\begin{rmq}\label{parite} Il y a une condition de parité sur les $N_j$ de la deuxième somme :
 si ${}^L G$ est un groupe symplectique, 
les $N_j$ impairs apparaissent avec une multiplicité paire
et  si $G$ est un groupe
orthogonal, les $N_j$ pairs apparaissent avec une multiplicité paire. On a bien sûr
\[ 2\sum_i N_i+\sum_j N_j=N.\]
\end{rmq}
Notons $\psi_{u,bp}$ la somme des sous-représentations intervenant dans la 
décomposition (\ref{decPsi}) de la forme $\Triv\boxtimes R_{N_j}$ avec $N_j$ pair si $^LG$ est un groupe 
symplectique et impair sinon et  $N_{u,bp}$ la dimension de cette 
représentation. D'après ce qui précède $N-N_{u,bp}$ est pair et il existe une représentation $\rho$ de 
$W_\bbC\times \SL_2(\bbC)$  dans $\GL_{N-N_{u,bp}}({\mathbb C})$, $\rho$, non unique telle que
\begin{equation}\label{decomposition} \psi=\rho \oplus \rho^* \oplus \psi_{u,bp}.\end{equation}

Notre but est la description des éléments du paquet $\Pi(\psi_G,G)$ attaché par Arthur à un  paramètre $\psi_G$
dans les cas  {$\mathbf C_n$},  {$\mathbf B_n$} et {$\mathbf{D}_n$}.
 L'importance de la décomposition (\ref{decomposition}) est que nous allons ramener cette description
  à celle de $\Pi(\psi_{u,bp})$ et utiliser  \cite{BV} pour décrire $\Pi(\psi_{u,bp})$.

\begin{rmq}
 Dans le cas {$\mathbf{Mp}_n$}, 
les paquets d'Arthur $\Pi(\psi, \Mp_{2n}(\bbC))$ sont définis via la correspondance de Howe entre 
groupes métaplectiques et groupes orthogonaux impairs (\cite{pourhowe}).
\end{rmq}

Revenons aux cas traités par Arthur, et donc $G$ n'est pas un groupe métaplectique.
Soit $\psi_G$ un paramètre d'Arthur et $\psi=\Std_G\circ \psi_G$. Soit $X^{\GL}_\psi$
la représentation auto-duale irréductible de $\GL_N(\bbC)$  qui est l'unique élément du paquet d'Arthur
$\Pi(\psi,\GL_N(\bbC))$ ({\sl cf. (\ref{PsiSpeh2})}). 
Les éléments du paquets d'Arthur $\Pi(\psi,G)$, et plus précisément, les représentations $X_\eta$
associées à chaque caractère de $A(\psi_G)$ sont caractérisés par les identités de transfert pour 
l'endoscopie ordinaire d'un groupe endoscopique elliptique de $G$ vers $G$, et l'identité de transfert
endoscopique tordu, où $G$ est un groupe endoscopique tordu pour $(\GL_N(\bbC), \theta_N)$
qui stipule que le transfert de la représentation virtuelle stable $X_{\psi_G}^{st}$ (\ref{distrstable})
est la trace tordue de la représentation  $X^{\GL}_\psi$.

\medskip 

\begin{rmq}\label{carinf}
Expliquons comment déterminer le caractère infinitésimal des éléments du paquet $\Pi(\psi,G)$. 
Rappelons qu'un caractère infinitésimal est donné par deux éléments $\lambda$ et $\mu$
de $\frh^*$, c'est-à-dire deux $n$-uplets de nombres complexes, modulo l'action de groupe de Weyl, 
c'est-à-dire à permutation et changements de signes près.
Pour  trouver $\lambda$ (resp. $\mu$), on considère la réunion des  
 $ (a+\frac{r-1}{2}, a+\frac{r-1}{2} \ldots ,a-\frac{r-1}{2})$ (resp.   $(b+\frac{r-1}{2}, b+\frac{r-1}{2} 
\ldots ,b-\frac{r-1}{2})$) pour 
chaque terme de la forme $\chi_{a,b}\boxtimes R_r$ apparaissant dans le paramètre $\psi$ 
(éventuellement avec $a=b=0$ pour les termes $\Triv\boxtimes R_r$). On obtient ainsi un $2n$-uplet
dans les cas {$\mathbf C_n$},  {$\mathbf D_n$} et {$\mathbf{Mp}_n$}, et un $2n+1$-uplet dans le cas {$\mathbf B_n$}, 
$0$ apparaissant avec une multiplicité 
impaire. Dans ce dernier cas, on enlève un $0$, et il reste donc un 
$2n$-uplet ayant la propriété que si $x$ en est  un élément, alors $-x$ aussi, avec la même multiplicité. 
On enlève la moitié des éléments, en groupant les éléments par paires $\{x,-x\}$ et en ne gardant pour 
chaque paire qu'un seul des deux éléments.
\end{rmq}

\subsection{Réduction aux paramètres unipotents de bonne parité}  \label{reducParUnip}

Nous énonçons dans cette section un résultat de réduction qui ne s'applique pas aux groupes métaplectiques.
En revanche, certains résultats intermédiaires importants d'irréductibilité d'induites
sont eux valides  aussi pour les groupes métaplectiques.
Soit $G$ un   groupe classique complexe   et soit
   \[\psi_G :\; W_\bbC\times \SL_2(\bbC)\longrightarrow {}^L G\] 
     un paramètre d'Arthur  pour $G$.        Posons  $\psi=\Std_G \circ \psi_G$ et soit   
              \begin{equation} \label{decPsiunip}
      \psi= \bigoplus_i \big( \chi_i\boxtimes R_{N_i} \oplus 
      \chi_i^{-1}\boxtimes R_{N_i} \big)
 \bigoplus_j (\Triv\boxtimes R_{N_j}) 
\end{equation}    
la décomposition en irréductibles de $\psi$ comme dans la section  précédente.

\bigskip

Supposons que $G$ ne soit pas un groupe métaplectique. 
Considérons une décomposition de $\psi$ de la forme :
\begin{equation}\label{psipsi}
       \psi=\rho\oplus \rho^*\oplus \psi' 
 \end{equation}
où $\rho=  \bigoplus_k  \chi_k\boxtimes R_{N_k}$, et $\rho^*=  \bigoplus_k  \chi_k^{-1}\boxtimes R_{N_k}$.
Ici les $\chi_k$ peuvent être triviaux mais s'ils sont triviaux alors $N_k$ est de mauvaise parité 
c'est-à-dire est pair si $^LG$ est un groupe orthogonal et impair si $^LG$ est un groupe symplectique. 
Le paramètre $\psi'$ se factorise par le $L$-groupe d'un groupe classique 
$G'$ de même type que $G$. Soit $\psi_{G'}$ le paramètre d'Arthur pour le groupe $G'$ tel que $\psi'=\Std_{G'}
\circ\psi_{G'}$. Notons $N_\rho=\sum_k N_k$ la dimension de la représentation $\rho$
de $W_\bbC \times \SL_2(\bbC)$, et soit $ X^{\GL}_\rho$
la représentation de $\GL_{N_\rho}(\bbC)$ de paramètre d'Arthur $\rho$ ({\sl cf. }(\ref{PsiSpeh2})).
 Le groupe $G$ admet un sous-groupe de Levi maximal isomorphe
à $\GL_{N_\rho}(\bbC)\times G'$, et ceci fournit une injection
\begin{equation}\label{injGG}
\iota: \;  \GL_{N_\rho}(\bbC) \times {}^LG'\hookrightarrow {}^LG
\end{equation}
de sorte que $\psi_G=\iota\circ (\rho, \psi_{G'})$.

\begin{rmq}\label{natiso} Les groupes $A(\psi_G)$ et $A(\psi_{G'})$ sont naturellement isomorphes. En effet,
le centralisateur de $ \psi$ est un produit du centralisateur de $\psi'$ et de groupes généraux linéaires
complexes.        
\end{rmq}

\begin{prop} \label{reducPaqArt}
Soit $\eta \in \widehat{A(\psi_G)}$ et soient $X_\eta$ et $X'_\eta$ les représentations semi-simples de $G$ et $G'$
respectivement attachées par Arthur ({\sl cf.} (\ref{Xeta}), où pour $X'_\eta$ on tient compte de la remarque 
ci-dessus). On a alors
\begin{equation}\label{XGLstarr}
X_\eta =   X^{\GL}_\rho \star X'_\eta.
  \end{equation}
\end{prop}

\dem Nous allons démontrer que le terme de droite de (\ref{XGLstarr}) vérifie les identités endoscopiques
qui caractérisent le terme de gauche.
Soit  ${\bf H}=(H,s,\xi: {}^LH\rightarrow {}^LG, \ldots)$ une  donnée endoscopique elliptique de $G$
({\sl cf.} \cite{Art13}) telle que $\psi_G$ se factorise 
par le groupe dual de $H$ et on fixe une telle factorisation $\psi_G=\xi \circ\psi_H$. 
En particulier  l'élément $s\in {}^LG$
 s'identifie à un élément du commutant de $\psi_G$. Il faut alors démontrer qu'il 
  existe une donnée endoscopique elliptique ${\bf H}'=(H',s',\ldots)$ de $G'$,  tel que l'élément $s'$ de 
 cette donnée soit dans le centralisateur de $\psi_{G'}$ et tel que le transfert de la distribution stable 
 associée à ${\bf H}$ et à la factorisation de $\psi_G$ soit l'induite du produit tensoriel 
  des données analogues pour $\psi_{G'}$ et   ${\bf H'}$ et de  la représentation $ X^{GL}_\rho$. 
  Expliquons maintenant comment construire explicitement cette donnée endoscopique ${\bf H}'$.
 Comme on a le droit de le faire, on suppose que $s$ vérifie $s^2=1$ et dans un premier temps on suppose 
 aussi que $s\neq 1$. 
 On décompose alors $\psi$ en $\psi_{+}\oplus \psi_-$ suivant les valeurs propres de $s$. 
 On remarque que l'on a aussi une décomposition analogue pour $\psi'$ et pour $\rho$. On a alors 
$$\psi_{+}= \rho_{+}\oplus \rho_{+}^* \oplus \psi'_{+}$$ 
et une décomposition analogue avec $+$ remplacé par $-$.
C'est ici qu'a servi l'hypothèse sur la parité de $N_k$ si $\chi_k$ est trivial pour que le dual 
de $\rho_+$ apparaisse lui aussi dans  dans l'espace propre de valeur propre $+1$.
Notons $N_{\rho_\pm}$  les dimensions des représentations $\rho_\pm$, 
et $X^\GL_{\rho_\pm}$ la représentation de $\GL_{N_{\rho_\pm}}(\bbC)$ associée à ce paramètre.
On a bien s\^ur  $N_\rho=N_{\rho_+}+N_{\rho_-}$ et $X^\GL_{\rho}= X^\GL_{\rho_+}\star X^\GL_{\rho_-}$.
Ainsi il existe un sous-groupe de Levi 
$$M\simeq \left( \GL_{N_{\rho_+}}(\bbC)\times M^+\right)\times \left(\GL_{N_{\rho_-}} (\bbC)\times M^-\right)$$ 
de  $H$ tel que $\psi_H$ se factorise par le groupe dual de $M$ 
et la distribution stable de $H$ associée à $\psi_H$  ({\sl cf.} (\ref{distrstable}))
 est une induite à partir de ce Levi. Notons $H'$ le facteur    $M^+\times M^-$
 de $M$ : c'est un groupe endoscopique pour $G'$, s'inscrivant dans une donnée endoscopique 
 ${\bf H}'=(H',s', \xi',\ldots)$ de $G'$ et le paramètre d'Arthur $\psi_{G'}$ se factorise
 en $\xi'\circ \psi_{H'}$. L'élément $s'$ est dans le centralisateur de $\psi_{G'}$, on
 peut le prendre tel que ${s'}^2=1$ et $\psi'={\psi'}^+\oplus {\psi'}^-$ est la décomposition de 
 $\psi'$ selon les valeurs propres $\pm 1$ de $s'$. 
 Partons de la représentation stable $X_{\psi_{H'}}^{st}$ 
associée à $\psi_{H'}$ ({\sl cf.} (\ref{distrstable})). On peut d'abord considérer son transfert endoscopique 
vers $G'$, puis induire vers $G$ avec $X^{\GL}_{\rho}$ : 
 $$  X^\GL_{\rho}\star \mathrm{Trans}_{H'}^{G'} (  X_{\psi_{H'}}^{st})  $$
où  désigne le transfert endoscopique (spectral) du groupe endoscopique $H'$ de $G'$ vers $G'$.
Le fait que le transfert commute à l'induction nous dit que l'on obtient le même résultat en induisant 
vers $H$ avec $X^\GL_{\rho_+}$ et $X^\GL_{\rho_-}$ et en prenant ensuite le transfert 
endoscopique de $H$ vers $G$ : 
$$ \mathrm{Trans}_{H}^{G}( X^\GL_{\rho_+}\star   X^\GL_{\rho_-} \star   X_{\psi_{H'}}^{st})= 
\mathrm{Trans}_{H}^{G}(X^\GL_{\rho} \star   X_{\psi_{H'}}^{st}).$$
Ce que l'on obtient est la représentation $X_\psi^s$ de (\ref{distrpasstable}).

Pour la distribution stable (\ref{distrstable}) associée à $\psi$ (le cas $s=1$),
  il faut montrer que le transfert
 pour l'endoscopie tordue de cette représentation vers l'espace tordu associé à $\GL_N({\mathbb C})$
  est aussi une induite. On vérifie aisément l'égalité des traces tordues de 
$X^{\GL}_{\psi}$ et   $X^{\GL}_\rho \times X^{\GL}_{\psi'}\times  X^{\GL}_{\rho^*}$ (il y a même en fait
égalité de   représentations tordues).
Il suffit alors de vérifier que la trace tordue de la représentation de droite est bien
 le transfert pour l'endoscopie tordue de la représentation 
 $  X^{\GL}_\rho \star X_{\psi_{G'}}^{st}$. Cela est aussi dû au fait que le transfert commute à l'induction.

On a ainsi vérifié que toutes les formules de transfert se déduisent pour $\psi$ de leurs analogues pour $\psi'$
 en induisant avec la représentation $X^{\GL}_\rho$. La proposition s'en déduit donc.
 \qed

\medskip

    Revenons à la décomposition (\ref{decPsiunip})     du début de paragraphe.
Alors   $\psi_{u}=\bigoplus_j (\Triv\otimes R_{N_j})$ est un  paramètre
unipotent de $G_{N_u}$ (où $N_u=\sum_j N_j$), c'est-à-dire que la restriction de 
$\psi_{u}$     à $W_\bbC={\mathbb C}^\times$ est triviale. Soit $G_u$ le groupe classique de m\^eme type que celui
de $G$ admettant une représentation standard dans $\GL_{N_u}(\bbC)$. Alors 
$\psi_{u}$ se factorise en $\psi_{u}= \Std_{G_u} \circ \psi_{G_u}$, où 
\begin{equation}
\psi_{G_u} : \; W_\bbC \times \SL_2(\bbC) \longrightarrow {}^ L G_u 
 \end{equation}
 est un paramètre d'Arthur unipotent de $G_u$.
        
Le groupe $G$ admet un sous-groupe parabolique $P=MN$ tel que $M$ soit isomorphe
à $ \left(\times_i  \GL_{N_i}(\bbC) \right) \times G_u$. Le paramètre $\psi_G$ se factorise par un paramètre
$\psi_M : \; W_\bbC \times \SL_2(\bbC) \longrightarrow {}^ L M $.
Si l'on identifie $M$ et $ \left( \times_i \GL_{N_i}( \bbC) \right)\times G_u$
comme expliqué ci-dessus,    les $\bar X^M \in \Pi(\psi_M, M)$ sont de la forme
\begin{equation}
\bar X^M=  \left( \bigotimes_i \chi_{i} \circ \textstyle \det_{N_i} \right) \otimes    \bar X_u 
 \end{equation}
où $\bar X_u$ décrit le paquet unipotent $\Pi(\psi_{G_u},G_u)$ de $G_u$.

Or, dans cette situation, nous avons le résultat suivant dû à Barbasch (\cite{B1}, Thm. 14.1) (que l'on redémontrera) : 

\begin{thm}\label{Ba141}
Avec les notations qui précèdent, pour toute représentation $\bar X_u$ dans  $\Pi(\psi_{G_u},G_u)$, 
la représentation induite
\[ \bar X= \Ind_P^G ( \left( \bigotimes_i \chi_{i}\circ \textstyle \det_{N_i} \right) 
\otimes \bar X_u  )=
\star_i \left(\chi_{i}\circ \textstyle \det_{N_i} \right)\star \bar X_u  \] 
 est irréductible.
\end{thm}
\begin{rmq}\label{nonramanujan}
Le théorème précédent reste vrai si les caractères $\chi_i$ sont seulement presque unitaires c'est-à-dire vérifient
 $\Re e\, (\lambda_i+\mu_i)\in ]-1/2,1/2[$.
\end{rmq}
On passe du théorème à la remarque on utilisant la description des séries complémentaires données en (\cite{B1} \S 12) 
mais cela peut se démontrer de façon totalement élémentaire.       

\
       
Lorsque $G$ n'est pas métaplectique, on peut appliquer la proposition \ref{reducPaqArt}.
Nous montrerons plus loin que pour tout $\eta \in \widehat{A(\psi_u)}$, la représentation 
$X_\eta^{G_u}$ de $G_u$ associée par Arthur est irréductible ou nulle
et le paquet d'Arthur $\Pi(\psi_u,G_u)$ est constitué des $X_\eta^{G_u}$ non nuls lorsque 
 $\eta$ décrit  $\widehat{A(\psi_u)}$. Nous en déduisons le :

\begin{cor} \label{CorBa141}  Avec les notations qui précèdent, la représentation 
$X_\eta$ de $G$ associée par Arthur à un caractère $\eta \in \widehat{A(\psi)}$,
est la représentation  irréductible 
$X_\eta=\left(\star_i \chi_{i} \circ \textstyle \det_{N_i} \right) \star   X_\eta^ {G_u}$
(on rappelle ici la remarque \ref{natiso},   et  bien évidemment $X_\eta=0$   si   $ X_\eta^ {G_u}=0$).

En particulier le paquet $ \Pi(\psi_G, G)$ est constitué d'induites irréductibles : 
  \[ \Pi(\psi_G,G)= \Big\{    
\left(\star_i \chi_{i} \circ \textstyle \det_{N_i} \right) \star  \bar X_u  ,\;  
\bar X_u \in  \Pi(\psi_{G_u}, G_u) \Big\}.\]     
\end{cor}

\bigskip 

Nous nous sommes donc ramenés à la détermination des paquets d'Arthur unipotents. Supposons  maintenant 
que $\psi_G=\psi_{G_u}$ et $\psi= \Std_G \circ \psi_G =\psi_u=\bigoplus_j (\Triv\boxtimes R_{N_j})$.

\begin{defi} \label{bonneparite}Nous définissons la {\sl bonne parité}  $\epsilon_G$ 
pour le groupe $G$ comme étant $1 \mod 2$ dans les cas {$\mathbf C_n$}
et {$\mathbf D_n$} (le groupe dual est un groupe orthogonal), et   $0 \mod 2$ dans les cas {$\mathbf B_n$} 
et {$\mathbf{Mp}_n$}(le groupe dual est 
un groupe symplectique).  
\end{defi}

Le paramètre $\psi$ se décompose en une partie de bonne parité, et une partie de mauvaise parité : 
\[\psi=\psi_{bp}\oplus\psi_{mp},  \quad  \psi_{bp}= 
\bigoplus_{\substack{j\\ N_j \, \mathrm{mod}\,  2=\epsilon_G}  }   (\Triv\boxtimes R_{N_j}), \quad  \psi_{mp}= 
\bigoplus_{\substack{j\\ N_j+1 \, \mathrm{mod}\,  2=\epsilon_G}  }   (\Triv\boxtimes R_{N_j}) .
\]
Les parties de mauvaise parité apparaissent avec une multiplicité paire ({\sl cf.} remarque \ref{parite}),
 et l'on va donc plutôt écrire : 
\[ \psi_{mp}=  \bigoplus_{\substack{ k\\ N_k+1 \, \mathrm{mod}\,  2=\epsilon_G}  }   (\Triv\boxtimes R_{N_k}\oplus 
 \Triv\boxtimes R_{N_k}) .
\]

On a $N=2 N_{mp}+N_{bp}= 2 \sum_k N_k+\sum_jN_j$ et   $\displaystyle\psi_{bp}=
\bigoplus_{\substack{j\\ N_j \, \mathrm{mod}\,  2=\epsilon_G}} (\Triv\otimes R_{N_j})$ 
est un  paramètre unipotent de $G_{N_{bp}}$. Soit $G_{bp}$ le groupe classique de m\^eme type que celui
de $G$ admettant une représentation standard dans $\GL_{N_{bp}}(\bbC)$. Alors 
$\psi_{bp}$ se factorise en $\psi_{bp}= \Std_{G_{bp}} \circ \psi_{G_{bp}}$, où 
\begin{equation}
\psi_{G_{bp}} : \; W_\bbC \times \SL_2(\bbC) \longrightarrow {}^ L G_{bp} 
 \end{equation}
 est un paramètre d'Arthur unipotent de bonne parité de  $G_{bp}$.
 
 Le groupe $G$ admet un sous-groupe parabolique $P=MN$ tel que $M$ soit isomorphe
à $\left(\times _k \GL_{N_k}(\bbC)\right)\times G_{bp}$. Le paramètre $\psi$ se factorise par un paramètre
$\psi_M : \; W_\bbC \times \SL_2(\bbC) \longrightarrow {}^ L M $.
Si on identifie $M$ et $\left(\times _k \GL_{N_k}(\bbC)\right)\times G_{bp}$, 
 les $\bar X^M \in \Pi(\psi_M, M)$ sont de la forme
\begin{equation}
\bar X^M=  \left( \bigotimes_k \Triv_{N_k}  \right) \otimes \bar X_{bp}
 \end{equation}
où $\bar X_{bp}$ décrit le paquet unipotent $\Pi(\psi_{G_{bp}}, G_{bp})$ de $G_{bp}$.
Là encore, nous avons le résultat d'irréductibilité suivant, dû à Barbasch (\cite{B1}, \S 6.6) : 

\begin{thm}\label{redbp}
Avec les notations qui précèdent, pour toute représentation $\bar X_{bp}$ dans  $\Pi(\psi_{G_{bp}}, G_{bp})$
\[ \bar X=  
\left( \star_k \Triv_{N_k}  \right)\star \bar X_{bp}. \] 
 est irréductible.
\end{thm}
        
    De m\^eme que ci-dessus, lorsque $G$ n'est pas métaplectique, on peut appliquer la proposition \ref{reducPaqArt} et obtenir le
\begin{cor} \label{CorBa142}  Avec les notations qui précèdent, la représentation 
$X_\eta$ de $G$ associée par Arthur à un caractère $\eta \in \widehat{A(\psi_G)}$,
est la représentation  irréductible 
\[ X_\eta=\left(\star_k \Triv_{k} \circ \textstyle \det_{N_k} \right) \star   X_\eta^{G_{bp}} \]
(on rappelle ici la remarque \ref{natiso},   et  bien évidemment $X_\eta=0$   si   $ X_\eta^ {G_{bp}}=0$).

En particulier le paquet $ \Pi(\psi_G, G)$ est constitué d'induites irréductibles : 
  \[ \Pi(\psi_G,G)= \Big\{    
\left(\star_k \Triv_{k} \circ \textstyle \det_{N_k} \right) \star  \bar X_{bp}  , \; 
\bar X_{bp }\in  \Pi(\psi_{G_{bp}}, G_{bp}) \Big\}.\]     
\end{cor}

Revenons maintenant à un paramètre $\psi_G$  quelconque, avec 
{\small    \begin{align*}
&\psi=\Std_G \circ \psi_G = 
&\bigoplus_i \big( \chi_{i}\boxtimes R_{N_i} \oplus 
\chi_{i}^{-1}\boxtimes R_{N_i} \big) 
\bigoplus_{k}   (\Triv\otimes R_{N_k} \oplus \Triv\boxtimes R_{N_k}  ) \bigoplus_j (\Triv\boxtimes R_{N_j})  
\end{align*}}    
où la première somme est la partie non unipotente, la seconde, la partie unipotente de mauvaise parité, et 
la troisième, la partie unipotente de bonne parité. On retrouve (\ref{decomposition}).
  Soit $G_{bp}$ le groupe classique de m\^eme type que celui
de $G$ admettant une représentation standard dans $\GL_{N_{bp}}(\bbC)$, 
où $N_{bp}=\sum_j N_j$ est la dimension de la représentation $\psi_{u,bp}$ de (\ref{decomposition}).
 
 Le groupe $G$ admet un sous-groupe parabolique $P=MN$ tel que $M$ soit isomorphe
à $ \left(\times _i \GL_{N_i}(\bbC)\right)\times  \left(\times_k \GL_{N_k}(\bbC)\right) \times G_{bp} $. 
Le paramètre $\psi$ se factorise par un  paramètre d'Arthur 
$\psi_M $ pour $M$. 
Si on identifie $M$ et $ \left(\times _i \GL_{N_i}(\bbC)\right)\times  \left(\times_k \GL_{N_k}(\bbC)
\right)\times G_{bp}$,
  les $\bar X^M  \in \Pi(\psi_M, M)$ sont de la forme
\begin{equation}
\bar X^M = 
\left( \bigotimes_i \chi_{i} \circ \textstyle \det_{N_i} \right) 
\otimes \left( \bigotimes_k \Triv_{N_k}  \right) \otimes \bar X_{bp}
 \end{equation}
où $\bar X_{bp}$ décrit le paquet unipotent $\Pi(\psi_{{u,bp}}, G_{bp})$ de $G_{bp}$

Les deux réductions effectuées ci-dessus peuvent se résumer en une seule grâce 
à la transitivité du foncteur d'induction parabolique :

\begin{thm}\label{Grossereduc}
Avec les notations qui précèdent, pour toute représentation $\bar X_{bp}$ dans  $\Pi(\psi_{G_{bp}}, G_{bp})$
\[ \bar X=  
\left(\star_i \chi_{i}\circ \textstyle \det_{N_i} \right) \star
 \left(  \star_k \Triv_{N_k}  \right)\star \bar X_{bp}.
\] 
 est irréductible et  si $G$ n'est pas un groupe métaplectique,  la représentation 
$X_\eta$ de $G$ associée par Arthur à un caractère $\eta \in \widehat{A(\psi_G)}$,
est la représentation  irréductible 
\[ X_\eta=\left( (\star_i \chi_{i}\circ \textstyle \det_{N_i})\star (\star_k \Triv_{k} \circ \textstyle \det_{N_k}) 
\right) \star   X_\eta^ {G_{bp}} \]
(on rappelle ici la remarque \ref{natiso},   et  bien évidemment $X_\eta=0$   si   $ X_\eta^ {G_{bp}}=0$).
 
En particulier  le paquet $ \Pi(\psi,G)$ est constitué d'induites irréductibles : 
  \[   \Pi(\psi,G)=\Big\{ \bar X= 
\left(\star_i \chi_{i}\circ \textstyle \det_{N_i} \right)\star 
\left( \star_k \Triv_{N_k}  \right)\star \bar X_{bp} , \; \bar X_{bp}  \in  
  \Pi(\psi_{{u,bp}}, G_{bp}) \Big\}.\]     
\end{thm}

\section{Description des paquets  unipotents (Barbasch-Vogan)}\label{descrip}
Soit $G$ l'un des groupes classiques du paragraphe \ref{grcla}, de rang $n$.
Soit $\psi_G$ un paramètre d'Arthur unipotent, c'est-à-dire que 
$\psi_G$ est trivial sur $W_\bbC={\mathbb C}^\times$. On peut donc le voir comme un morphisme de $\SL_2(\bbC)$
dans ${}^LG$. Soit 
\begin{equation}\label{diffpsi}
\bar \psi_G  :\, \frs\frl_2(\bbC)\rightarrow \mathrm{Lie}({}^ L G)= {}^ L\frg \end{equation}
la  différentielle de $\psi_G$ en l'identité
qui envoie le  $\frs\frl(2,\bbC)$-triplet usuel  $\{ e,h,f \}$ sur un  $\frs\frl(2)$-triplet  dans
 ${}^ L \frg$ (les images de  $e$ et $f$ sont des éléments  nilpotents  dans  ${}^L \frg$).
Notons $\caU=\caU_{\psi_{G}}$ l'orbite nilpotente  ${}^LG\cdot e$.
La théorie de Jacobson-Morosov établit que la correspondence $\psi_G \leftrightarrow \caU_{\psi_{G}}$
entre classes de conjugaison de paramètres d'Arthur unipotents et orbites nilpotentes dans ${}^ L\frg$ 
est bijective.

\bigskip

{\bf Orbites nilpotentes dans les algèbres de Lie simples classiques}.
Rappelons la classification des orbites nilpotentes dans les algèbres de Lie simples classiques ({\sl cf. } 
\cite{CMcG}).
La représentation standard $\Std_G: \, {}^LG \rightarrow G_N$ donne par différentiation un morphisme 
d'algèbres de Lie injectif :
\[\overline \Std_G  :\,  {}^ L\frg \longrightarrow \frg\frl_N(\bbC) \]
et l'image d'une orbite nilpotente $\caU$ de ${}^L\frg$ est une orbite nilpotente
de $\frg\frl_N(\bbC)$. Une telle orbite est caractérisée par une partition de $N$, qui donne la taille 
des blocs de Jordan d'un élément de l'orbite. 

\medskip 

On note $\caP(N)$ l'ensemble des partitions de $N$. Les  éléments de $\caP(N)$ sont des suites
$\mathbf{d}=[d_1,  d_2,  \ldots ,  d_k ]$ avec $d_1\geq d_2\geq \ldots \geq d_k> 0$ et 
$\sum_{i=1}^k d_i=N$.
On note $r_{\mathbf{d}}(i)$ la multiplicité de l'entier  strictement positif $i$ dans $\mathbf{d}$.
On définit  $\caP_{1}(N)$ comme le sous-ensemble de   $\caP(N)$ des partitions
$\mathbf{d}=[d_1,  d_2,  \ldots ,  d_k ]$ de $N$ telles que la multiplicité 
$r_{\mathbf{d}}(i)$ de tout entier $i$ pair soit paire.
De même, on définit   $\caP_{-1}(N)$ par la condition que la multiplicité de tout $i$
impair soit paire (et donc $N$ est pair).

\medskip 

Reprenons nos groupes classiques du paragraphe \ref{grcla}.
 Dans le cas {$\mathbf C_n$}, ${}^L\frg=\frs\fro_{2n+1}$, et une 
orbite   $\caU$ de ${}^L\frg$ donne via $\overline \Std_G $ une partition
 $\mathbf{d}=\mathbf{d}_{\caU}=[d_1, d_2, \ldots,  d_k ]$ de $N=2n+1$.
Cette partition est dans $\caP_1(2n+1)$ et la correspondance $\caU \mapsto  \mathbf{d}_{\caU}$ 
est une bijection entre l'ensemble des orbites nilpotentes de  $\frs\fro_{2n+1}$ et $\caP_1(2n+1)$.

De m\^eme, dans les cas {$\mathbf B_n$} et {$\mathbf {Mp}_n$}, ${}^L\frg=\frs\frp_{2n}$, et l'on a une bijection 
$\caU \mapsto  \mathbf{d}_{\caU}$   entre l'ensemble des orbites nilpotentes de 
 $\frs\frp_{2n}$ et $\caP_{-1}(2n)$. Dans le cas {$\mathbf D_n$}, ${}^L\frg=\frs\fro_{2n}$,
mais la correspondance   $\caU \mapsto  \mathbf{d}_{\caU}$ ne  donne plus  une bijection 
entre les orbites nilpotentes de  $\frs\fro_{2n}$ et $\caP_1(2n)$. Il faut remplacer l'action adjointe
du groupe $\SO_{2n}(\bbC)$ par le groupe orthogonal $\Or(2n)$ pour obtenir une bijection. En fait, une orbite
 nilpotente ${}\caO_{ \mathbf{d}_{\caU}}$ de  $\frs\fro_{2n}$ pour l'action de  $\Or(2n)$
  est une orbite pour  $\SO_{2n}(\bbC)$, sauf lorsque la partition  $ \mathbf{d}_{\caU}$ est \og très paire \fg, 
c'est-à-dire que tous les $d_i$ sont pairs (avec des multiplicités paires), auquel cas 
l'orbite ${}\caO_{ \mathbf{d}_{\caU}}$ se scinde en deux orbites 
${}\caO_{ \mathbf{d}_{\caU}}^I$, ${}\caO_{ \mathbf{d}_{\caU}}^{II}$ pour $\SO_{2n}(\bbC)$.
 
\medskip

 Revenons aux  paramètres d'Arthur unipotents $\psi_G$ du début du paragraphe, auxquels nous avons 
 attaché une orbite nilpotente $\caU$ de ${}^ L\frg$ et une partition
$\mathbf{d}=\mathbf{d_{\caU}}$. 
Il est donc équivalent de se donner une classe de conjugaison de paramètres d'Arthur unipotents, une orbite nilpotente
dans ${}^L\frg$, ou la partition qui lui correspond. Posons comme précédemment $\psi=\Std\circ \psi_G$, et soit 
$\textstyle \psi=\oplus_{j=1}^r \Triv \boxtimes R_{N_j}$ la décomposition de $\psi$ en irréductibles.
La partition associée à $\psi_G$ est alors $\mathbf{d}=[N_1,\ldots, N_r]$.

\begin{rmq}\label{AA}
Si l'orbite nilpotente $\caU$ est attachée au paramètre unipotent $\psi_G$, alors 
on a une identification $A(\psi_G)=A(\caU)$, où $A(\caU)$ est le groupe des composantes
connexes du centralisateur dans ${}^LG$ d'un élément  $e$ de $\caU$ ({\sl cf.} \cite{BV}, Prop. 2.4).
 Il sera commode quand $\caU$ est donnée 
de noter $\psi_{\caU}$ le morphisme correspondant et à l'inverse quand $\psi$ est donné de noter 
$\caU_\psi$ l'orbite correspondante.
\end{rmq}

\begin{defi}
 Si le paramètre $\psi$ est de bonne parité, 
c'est-à-dire que dans la partition $\mathbf{d}$, il n'y a que des $d_i$ vérifiant $d_i\mod 2=\epsilon_G$
({\sl cf.} définition \ref{bonneparite}), on dit aussi que $\caU$ et $\mathbf{d_{\caU}}$ sont de bonne parité.
\end{defi}

Dans les travaux de Barbasch-Vogan, une notion importante est celle d'orbite spéciale, qui est basée sur la dualité
de Lusztig-Spaltelstein ({\sl cf.} \cite{BV}, Appendix). Cette dualité est une application $D_\frg$
 de l'ensemble des orbites nilpotentes d'une algèbre de Lie réductive $\frg$ vers l'ensemble des orbites
  nilpotentes de ${}^L\frg$ vérifiant $D_\frg\circ D_{{}^L\frg}\circ D_\frg=D_\frg$.
Les orbites dans l'image de $D$ sont appelées orbites spéciales. Lorsqu'on restreint $D$ aux orbites
 spéciales, on obtient  une bijection échangeant orbites nilpotentes spéciales  de  $\frg$ et de ${}^L\frg$.
Les r\^oles de $\frg$ et ${}^L\frg$ sont ici bien évidemment totalement symétriques. Dans cet article,
axé sur le point de vue d'Arthur, ce sont les orbites nilpotentes de ${}^L\frg$ qui jouent le r\^ole
majeur, car elles donnent les paramètres, alors que dans \cite{BV}, l'accent est mis sur les orbites
 nilpotentes de $\frg$ qui donnent les fronts d'onde des représentations.
C'est la dualité de Lusztig-Spaltenstein qui fait le lien entre ces deux points de vue ({\sl cf.} \cite{BV} 
repris ici dans un cadre un peu plus général en \ref{frontdonde}).

 On trouve dans \cite{CMcG} une procédure expliquant comment calculer $D_\frg$ en terme de partitions, que nous 
 rappelons brièvement. Pour ceci, il nous faut une procédure permettant d'obtenir, à partir
 d'une partition $\mathbf{d}$ de $2n+1$ ou $2n$, une partition dans $\caP_{1}(2n+1)$, $\caP_{-1}(2n)$, $\caP_{1}(2n)$,
 appelées respectivement $B$-collapse, $C$-collapse et $D$-collapse de $\mathbf{d}$.
Si $ \mathbf{d}$ est une partition de $2n+1$ qui n'est pas dans $\caP_{1}(2n+1)$, alors un entier pair $r$
apparaît dans $ \mathbf{d}$ avec une multiplicité impaire. Considérons le plus grand de ces entiers,  retirons
1 au dernier bloc de 
cette taille, et rajoutons 1 au bloc de taille maximale strictement inférieure
à $r-1$. On répète ceci jusqu'à obtenir une partition dans $\caP_{1}(2n+1)$. La procédure est similaire pour 
le $C$-collapse (resp. $D$-collapse) d'une partition de $2n$, où l'on se débarrasse des blocs impairs (resp. pairs)
de multiplicités impaires en commençant par le plus grand.

 {\bf Cas $\mathbf{C}_n$}. Soit $\caO$ une orbite nilpotente de $\frs\frp_{2n}$ et $\mathbf{d}\in \caP_{-1}(2n)$
 la partition associée.  On ajoute un bloc $1$ pour avoir une partition de $2n+1$, on prend la partition transposée 
 (on échange lignes et colonnes du tableau de Young) et  on prend  le $B$-collapse de celle-ci, pour 
 obtenir une partition dans  $\caP_{1}(2n+1)$,  ce qui nous donne une 
 orbite nilpotente de $\frs\fro_{2n+1}$.

{\bf Cas $\mathbf{B}_n$}. Soit $\caO$ une orbite nilpotente de $\frs\fro_{2n+1}$ et $\mathbf{d}\in  \caP_{1}(2n+1)$
 la partition associée.  On prend la partition transposée  et on
 enlève   $1$ au plus petit bloc  pour avoir une partition de $2n$.  On prend le $C$-collapse de cette partition
 pour   obtenir une partition dans  $\caP_{-1}(2n)$,  ce qui nous donne une 
 orbite nilpotente de $\frs\frp_{2n}$.

{\bf Cas $\mathbf{D}_n$}.  On prend le $D$-collapse de la partition transposée en éliminant les blocs pairs ayant une multiplicité impaire.

{\bf Cas $\mathbf{Mp}_n$}.
Pour le groupe métaplectique, il y a aussi une dualité (étudiée en \cite{duke})
 de l'ensemble des  orbites nilpotentes de 
$\frs\frp_{2n} (\bbC)$  dans lui-même, définie comme suit.
Soit ${\mathcal U}$ une orbite nilpotente de $\frs\frp_{2n} (\bbC)$ et $\mathbf{d}_{\caU}\in \caP_{-1}(2n)$
la partition qui lui correspond. On ajoute un bloc 1 à $\mathbf{d}_{\caU}$,
on prend la partition conjuguée, $\mathbf{d}_{\caU}^t$ et on applique l'algorithme pour calculer 
le \og $C$-collapse \fg, 
en évacuant les blocs impairs de multiplicité impaires du plus grand au plus petit. Mais ici
comme on a une partition de $2n+1$, on se retrouve à la fin avec un seule bloc de taille impaire 
avec une multiplicité impaire. 
On enlève 1 au dernier de ces blocs. Notons $\mathbf{d}_{\caO}$ la partition de $\caP_{-1}(2n)$ obtenue et 
$\caO$ l'orbite de $\frs\frp_{2n} (\bbC)$ qui lui correspond.
Par exemple, l'orbite $\caU$ de partition $[2,1,1]$ a pour duale l'orbite $\caO$ de partition $[4]$,
 et réciproquement. De même, l'orbite $\caU$ de partition $[2,2,2]$ a pour duale l'orbite $\caO$ de partition $[4, 2]$, et réciproquement.

\begin{defi}
 On dit que l'orbite nilpotente  $\caU$ de  $\frs\frp_{2n} (\bbC)$ est antispéciale si elle est égale à sa biduale.  
 Par exemple, les orbites
 de partitions  $[2,1,1]$, $[4]$, $[2,2,2]$ et  $[4, 2]$ sont antispéciales.
 Celle de partition $[3,3]$ a pour duale $[2,2,2]$ et n'est donc pas antispéciale.
\end{defi}

\begin{rmq}
Une orbite nilpotente  de $\frs\frp_{2n} (\bbC)$ est à la fois spéciale et antispéciale 
si et seulement elle est de bonne parité 
({\sl ie.} tous ses blocs sont pairs).
\end{rmq}

\begin{rmq}
Pour les groupes classiques, une orbite est {\sl paire} (au sens usuel, {\sl cf.} \cite{BV} ou \cite{CMcG})
 si tous les $d_i$ sont de même parité. Une orbite nilpotente  $\caU$ de bonne parité est spéciale
 (antispéciale dans le cas $\mathbf{Mp}_n$)  et paire. 
Dans le cas {$\mathbf C_n$}, cette parité ne peut être que la bonne parité (impaire), 
car il y a un nombre impair de blocs
impairs. 
\end{rmq}

Dans le cas où $\caU$ est une orbite spéciale paire, et $G$ un groupe classique non métaplectique
une définition du  paquet d'Arthur  associé  a été donné par Barbasch et Vogan
 (\cite{BV}). Notons le  $\Pi_{BV}(\caU,G)$. Cette définition est étendue par Barbasch \cite{B1} au cas des groupes
 métaplectiques.
 
 \begin{defi} [Barbasch-Vogan] Soit $G$ un groupe classique et $\caU$ une orbite
 nilpotente spéciale (antispéciale dans  le cas $\mathbf{Mp}_n$) paire de ${}^L\frg$. 
Soit $\caO$ l'orbite duale de $\caU$. C'est une orbite nilpotente spéciale   (antispéciale dans 
le cas $\mathbf{Mp}_n$)  dans $\frg$.
Nous avons vu en (\ref{GC44}) que le caractère infinitésimal doit être 
 $\left( \frac{1}{2}h,\frac{1}{2}h\right)$.
Les éléments du paquet  $\Pi_{BV}(\caU,G)$ sont  les 
  modules de Harish-Chandra irréductibles $\bar X$ ayant ce caractère infinitésimal, 
et dont le front d'onde  est $\mathrm{WF}(X)=\overline{\caO}$.
\end{defi}

On suppose dans la suite que   $\caU$  est spéciale  (antispéciale dans 
le cas $\mathbf{Mp}_n$) et paire.

\begin{rmq}  
Pour les groupes classiques, on calcule facilement les coordonnées du caractère   
 infinitésimal de $\Pi_{BV}(\caU,G)$ à partir de l'orbite $\caU$, ou plut\^ot
 de la partition $\mathbf{d}_{\caU}$. Si $\mathbf{d}_{\caU}=[d_1,\ldots, d_k]$,
 on considère pour chaque $d_i$ la suite 
 \[ \left(\frac{d_i-1}{2}, \frac{d_i-3}{2}, \ldots,\frac{3-d_i}{2},  \frac{1-d_i}{2} \right). \] 
 On concatène toute ces suites en réordonnant les éléments dans l'ordre décroissant.
On ne garde  ensuite que les termes strictement positifs, et la moitié des termes nuls 
(s'il y en a un nombre impair $2\ell+1$, on en garde $\ell$). 
 
 Par exemple, dans le cas {$\mathbf C_n$}, considérons le paramètre unipotent associé à 
 $\mathbf{d}_{\caU}=[5,3,3,1,1]$. Le caractère infinitésimal est $(2,1,1,1,0,0)$.
Dans le cas {$\mathbf B_n$},   si  $\mathbf{d}_{\caU}=[6,6,4,2]$, 
le caractère infinitésimal est $(\frac{5}{2},\frac{5}{2},\frac{3}{2},\frac{3}{2},
\frac{3}{2},\frac{1}{2},\frac{1}{2},\frac{1}{2},\frac{1}{2})$.

Si l'orbite $\caU$ est paire, le caractère infinitésimal associé est entier.
Dans le cas {$\mathbf C_n$}, nous avons vu qu'une orbite paire est de bonne parité (tous les $d_i$ sont impairs), 
et ainsi les coordonnées du caractère infinitésimal sont des entiers.
 Dans le cas {$\mathbf B_n$}, une orbite paire est soit de bonne parité (paire), auquel cas
les coordonnées du caractère infinitésimal sont des demi-entiers (dans $\frac{1}{2}\bbZ\setminus \bbZ$), 
soit totalement de mauvaise parité, auquel cas les coordonnées du caractère infinitésimal sont des entiers. 
  Dans le cas {$\mathbf D_n$}, une orbite paire est soit de bonne parité (impaire), auquel cas
les coordonnées du caractère infinitésimal sont des entiers, soit totalement de mauvaise parité, 
auquel cas les coordonnées du caractère infinitésimal sont des demi-entiers.
 Dans le cas {$\mathbf{Mp}_n$}, une orbite paire antispéciale est de bonne parité (paire), auquel cas
les coordonnées du caractère infinitésimal sont des demi-entiers.
 \end{rmq}

 Nous allons donner  la  description des éléments du  paquet $\Pi_{BV}(\caU,G)$ en suivant \cite{BV} et \cite{B1}.
Ils sont paramétrés par les caractères    d'un certain quotient $\bar A(\caU)$  de $A(\caU)$. 

Pour les groupes classiques non métaplectiques, le quotient $\bar A(\caU)$ est le quotient de Lusztig.
Il  est isomorphe à $(\bbZ/2\bbZ)^m$ pour un certain entier $m$, ceci sera rendu explicite plus loin.
On a donc une bijection 
\begin{equation}\label{BVpar1} 
\widehat{\bar A(\caU)} \longrightarrow \Pi_{BV}(\caU,G), \qquad \eta\mapsto  X^{BV}_{\eta}
\end{equation}

\begin{rmq} En fait, Barbasch et Vogan \cite{BV} paramètrent les représentations dans un paquet $\Pi_{BV}(\caU,G)$
par les caractères du  quotient de Lusztig attaché à l'orbite $\caO$, duale de $\caU$. Ils remarquent que 
ces quotients sont isomorphes: $\bar A(\caO)\simeq \bar A(\caU)$, mais bien entendu, pour passer
de l'une à l'autre des paramétrisations, il faut dire quel est cet isomorphisme. Ceci est assez subtil, et nous 
l'expliquons plus loin ({\sl cf.} remarque  \ref{AOAU}).
 On se sert ici de la section 6 de \cite{B1}, où Barbasch donne la paramétrisation
des  représentations dans $\Pi_{BV}(\caU,G)$ par $\widehat{\bar A(\caU)}$.
\end{rmq}

Pour  \^etre en phase avec la paramétrisation d'Arthur, nous réinterprétons (\ref{BVpar1}) : à un caractère
$\eta$ de $A(\caU)$ est attachée une représentation du paquet  $\Pi_{BV}(\caU,G)$ ou bien $0$ 
de la manière suivante. Si $\eta$ est un caractère de $ A(\caU)$ qui ne se factorise pas par 
$\bar A(\caU)$, on pose alors $ X^ {BV}_\eta=0$. 
\begin{equation}\label{BVpar2} 
\widehat{A(\caU)} \longrightarrow \Pi_{BV}(\caU,G) \coprod \{0\}, \qquad \eta\mapsto  X^{BV}_{\eta}
\end{equation}

\bigskip 
Pour les groupes métaplectiques, les représentations dans le paquet attaché à une orbite nilpotente $\caU$
de l'algèbre de Lie ${}^L\frg=\frs\frp_{2n}(\bbC)$  sont paramétrés par les caractères d'un autre groupe 
quotient de $A( \caU)$ que nous allons noter $\overline{A}^{Mp}({\mathcal U})$, pour le distinguer
du groupe   $\overline{A}({\mathcal U})$ qui apparaît lui lorsque ${}^L\frg=\frs\frp_{2n}(\bbC)$ est 
vue comme l'algèbre de Lie du dual de Langlands de $\SO_{2n+1}$.

Soit ${\mathcal U}$  une orbite nilpotente antispéciale  paire de $\frs\frp_{2n}({\mathbb C})$.
 On note $\mathbf{d}_{\caU}=[d_1,  \cdots , d_t]$ la partition (dans $\caP_{-1}(2n)$) définissant ${\mathcal U}$.
  Tous les blocs sont  alors pairs, et l'on pose  $n_i=d_i/2$. Alors  $\times_{i\in [1,t]}
  \Sp_{2n_i}(\bbC)$ est un sous-groupe de $\Sp_{2n}({\mathbb C})$.
   Le centre de ce sous-groupe s'envoie surjectivement sur le groupe des 
composantes connexes du centralisateur d'un élément de ${\mathcal U}$ inclus dans ce sous-groupe. 
Pour tout $i\in [1,t]$, on note $\epsilon_i$ 
l'image de l'élément non trivial du centre de $\Sp_{2n_i}({\mathbb C})$ dans ce groupe de composantes connexes. 
En particulier $\epsilon_i=\epsilon_j$ si $n_i=n_j$. On  définit alors $\overline{A}^{Mp}({\mathcal U})$ comme le 
quotient du groupe des composantes connexes par le sous-groupe engendré par les 
éléments $\epsilon_{i}\epsilon_{i+1}$ où $i$ parcourt l'ensemble des entiers pairs dans $[1,t]$ et 
$\epsilon_{t+1}=1$ si $t$ est pair.

 \bigskip

 A la  partition $\mathbf{d}_{\caU}$ , on attache un symbole ({\sl cf.} \cite{CMcG} \S 10
 pour les groupes classiques non métaplectiques. Pour ces derniers, voir ci-dessous). 
Dans les cas  {$\mathbf C_n$}, et  {$\mathbf B_n$}   ce symbole est de la forme
\begin{equation} \label{symb1}
 \begin{pmatrix} x_0 &  & x_2 & & \ldots &&&& x_{2k}\\
&x_1 && x_3 && \ldots  && x_{2k-1}
\end{pmatrix} 
\end{equation}
 et dans le cas  {$\mathbf D_n$}, ce symbole est de la forme
\begin{equation} \label{symbC}
 \begin{pmatrix} x_0  & x_2 & & \ldots &&&& x_{2k}\\
x_1 & x_3 && \ldots  &&&& x_{2k+1}
\end{pmatrix} 
\end{equation}

 Un entier $i$ donné n'apparaît dans le symbole qu'au plus deux fois. Comme on est parti d'une orbite spéciale, 
 on a de plus :
 \[x_0\leq x_1\leq x_2 \ldots \leq x_i\leq x_{i+1}\leq \ldots   \]
 
 Voici comment on détermine $\bar A(\caU)$ : on enlève du symbole tous les $x_i$ apparaissant avec multiplicité
 $2$ (l'un apparaît sur la ligne du haut, l'autre sur celle du bas).
  Il reste un symbole de même type sans multiplicité.  Soit $m$ la  cardinalité  de la ligne inférieure.
  Alors   $\bar A(\caU)\simeq (\bbZ/2\bbZ)^m$.  
 
Donnons maintenant une description des $2^m$ éléments de $\Pi_{BV}(\caU,G)$, selon les cas, en illustrant
 ce qui se passe sur des exemples. On suit bien évidemment les descriptions de \cite{B1}.
Nous allons pour cela définir à chaque fois un sous-groupe parabolique de $G$, et $2^m$ représentations
 irréductibles de celui-ci. Les éléments du paquets seront les $2^m$  sous-quotients irréductibles
contenant les $K$-types minimaux des induites de $P$ à $G$ de ces représentations.
 
\bigskip

{\bf Cas $\mathbf{C}_n$}.  Nous illustrons ce cas par l'exemple 
 $\mathbf{d}_{\caU}=[9,5,5,5,3,1,1]$. Le rang est 14.
 
 On repart du symbole attaché à $\mathbf{d}_{\caU}$ :
 \begin{equation} \label{SymbBC}
 \begin{pmatrix} x_0 &  & x_2 & & \ldots &&&& x_{2k}\\
&x_1 && x_3 && \ldots  && x_{2k-1}
\end{pmatrix} 
\end{equation}

Dans notre exemple, on obtient 
\begin{equation} \label{SymbEx}
 \begin{pmatrix} 0&  & 2   && 4 && 7\\
&1 && 4&& 5   
\end{pmatrix} 
\end{equation}

\bigskip

On définit  une sous-algèbre $\frmm$ de $\frg=\frs \frp_{2n}$ de la manière suivante
\begin{equation}\label{defmcasA}  \frmm \simeq \frmm^0  \times \caA^1 \times \ldots \times \caA^k    \end{equation}
 où $\frmm^0 \simeq \frs \frp_{2x_0}$, $\caA^j\simeq \frg \frl (x_{2k-2j+2}+x_{2k-2j+1}-2k+2j-1 )$, $1\leq j\leq k$.

Dans notre exemple, on obtient 
\[  \frmm \simeq \frg \frl(7) \times \frg\frl(5) \times \frg\frl(2) \]

On écrit le caractère infinitésimal attaché à  $\caU$ en coordonnées :
\[ \lambda=(\lambda_1, \lambda_1, \ldots, \lambda_1,\lambda_2, \ldots, \lambda_2,   
 \ldots ,\lambda_s, \ldots, \lambda_s) \]
 avec $\lambda_i>\lambda_{i+1}$, $\lambda_s \geq 0$.
Les $\lambda_i$ sont des entiers dans ce cas. 
 
 Dans notre exemple,
 \[ \lambda=(4,3,2,2,2,2,1,1,1,1,1,0,0,0) .\]
  
 On extrait alors de $\lambda$ la suite  strictement décroissante maximale,
 $(\lambda_1,\lambda_2, \ldots ,\lambda_s)$ et on la complète en extrayant des opposés des éléments non 
 nuls du reste (que l'on retire aussi de $\lambda$)  une suite strictement décroissante maximale.
 Dans notre exemple, on extrait d'abord  ainsi $(4,3,2,1)$, il reste $(2,2,2,1,1,1,1,0,0)$ et donc on 
 complète en $(4,3,2,1,0,-1,-2)$. On note $F^1$ la représentation de dimension finie 
 holomorphe de $\caA^1$ de caractère infinitésimal donné par la suite ainsi obtenue. Dans notre exemple
 $F^1$ est  la représentation de dimension finie 
 holomorphe de $\frg \frl(7)$ de caractère infinitésimal $(4,3,2,1,0,-1,-2)$.
 
 On définit par le même procédé, à partir de ce qui reste de $\lambda$, une représentation $F^2$ de $\caA^2$, 
 et ainsi de suite jusqu'à  $F^k$, représentation de $\caA^k$. Dans notre exemple, après extraction du 
 caractère infinitésimal de $F^1$, il nous reste $(2,2,1,1,1,0,0)$ et le caractère infinitésimal de la représentation 
 holomorphe de dimension finie $F^2$ de $\frg\frl(5)$  est   donc $(2,1,0,-1,-2)$.
 Il nous reste  ensuite $(1,0)$ et le caractère infinitésimal de la représentation 
 holomorphe de dimension finie $F^3$ de $\frg\frl(2)$  est   donc $(1,0)$.

 Il ne reste après ceci de $\lambda$ qu'une suite d'entiers décroissante $\lambda^0$ avec  $x_0$  éléments, que que l'on voit comme un élément du dual de la sous-algèbre de Cartan de  $\frs \frp_{2x_0}$ (il ne reste plus rien dans notre exemple, puisque $x_0=0$).
  
   \medskip

Lorsque un entier $a$ apparaît avec multiplicité $2$ dans le symbole, on vérifie, en utilisant 
le fait que tous les blocs de $\mathbf{d}_{\caU}$ sont pairs, que c'est sous la forme suivante : 
\begin{equation*}
 \begin{pmatrix}
  x_0 &&  \ldots &&  x_{2j}=a && \ldots && x_{2k}\\
& \ldots  && x_{2j-1}=a  && \ldots  &&  x_{2k-1} &
\end{pmatrix} 
\end{equation*}

Lorsque l'on enlève ces paires $(x_{2j-1},x_{2j})$ avec $x_{2j-1}=x_{2j}$ du symbole, il nous reste donc 
$m$ paires   $(x_{2j-1},x_{2j})$ avec $x_{2j-1}>x_{2j}$. 
Pour une telle paire, on définit la représentation de dimension finie 
 holomorphe $\tilde F^j$ de $\caA^j$ de la manière suivante. On reprend $F^j$ et son caractère infinitésimal,
 et dans celui-ci, on change le signe de la plus petite coordonnée strictement positive dont l'opposé n'apparaît
  pas, et on réordonne dans l'ordre décroissant. Ceci nous donne le caractère infinitésimal de $\tilde F^j$.
   Dans notre exemple, on a $m=2$, $\tilde F^1$ a pour caractère infinitésimal $(4, 2,1,0,-1,-2,-3)$ et  
 $\tilde F^3$ a pour caractère infinitésimal $(0,-1)$.

On note alors, pour tout $1\leq j\leq k$, tel que   $x_{2j-1}=x_{2j}$,
\begin{equation}\label{caF} \scrF_1^j=F^j \otimes F^j,  \end{equation}
et pour tout $1\leq j\leq k$  tel que   $x_{2j-1}<x_{2j}$,  
\begin{equation}\label{caF2} \scrF^j_{1}= F^ j \otimes F^j,\qquad \scrF^j_{-1}= F^j \otimes \tilde F^j.  \end{equation}
Posons :
\begin{equation} \label{epsj} 
\hat{\boldsymbol{\epsilon}}=( \hat \epsilon_1, \ldots, \hat \epsilon_k )\qquad \text{ où  } \hat \epsilon_j=\pm 1, \; 
\text{ et } \hat \epsilon_j=1  \text{ si   } x_{2j-1}=x_{2j}, 
\end{equation}
de sorte que 
\begin{equation} \label{cafe} \scrF_{\hat{\boldsymbol{\epsilon}}}=\otimes_j\scrF^j_{\hat\epsilon_j}
\end{equation}
est une représentation de dimension finie de $\prod_{j=1}^k \caA^j$, que l'on relève, avec la même notation, 
en une représentation de dimension finie du produit des groupes généraux linéaires d'algèbre de Lie  $\caA^j$.  

Soit $P=MN$ un sous-groupe parabolique de $G$ dont le facteur de Levi a pour algèbre de Lie $\frmm$.
C'est un produit  des groupes généraux linéaires d'algèbre de Lie  $\caA^j$, $j=1, \ldots ,k$, et 
d'un facteur $G(x_0)=\Sp_{2x_0}(\bbC)$. Considérons la représentation sphérique 
$\bar X(\lambda^0,\lambda^0)$ de ce groupe\footnote{C'est la représentation triviale de $\Sp_{2x_0}(\bbC)$.},
 et définissons
  \begin{equation} \label{defpaq} \bar X_{\hat{\boldsymbol{\epsilon}}} : \text{sous-quotient irréductible
   contenant le $K$-type minimal de } \Ind_P^ G \left(\bar X(\lambda^0,\lambda^0) \otimes 
   \scrF_{\hat{\boldsymbol{\epsilon}}} \right).
\end{equation}

Le paquet associé à $\caU$ par Barbasch et Vogan est alors 

 \begin{equation} \label{defpaq2} \Pi_{BV}(\caU,G)=\{  \bar X_{\hat{\boldsymbol{\epsilon}}} \}.
 \end{equation}

\begin{rmq}\label{parpaq}
Le cardinal de ce paquet est le nombre de choix possibles pour les signes $\hat{\epsilon_j}$, 
c'est-à-dire $2^m$ avec les notations qui précèdent. En fait, chaque
 $\hat{\boldsymbol{\epsilon}}=( \hat \epsilon_1, \ldots, \hat \epsilon_k )$ définit de manière évidente
un caractère de $\bar A(\caU)$, et les éléments du paquet sont donc paramétrés
par ces caractères.
\end{rmq}

\medskip

Revenons à notre exemple. On obtient donc dans ce cas un paquet à 4 éléments, 
les sous-quotients contenant les
$K$-types minimaux des induites de 
$\GL_7(\bbC)\times \GL_5(\bbC)\times \GL_2(\bbC)$ des caractères 
$\scrF_{1,1,1}$, $\scrF_{1,1,-1}$, $\scrF_{-1,1,1}$, $\scrF_{-1,1,-1}$.
Renommons-les, en les indexant par les caractères du groupes $\bar A(\caU)=\left(\bbZ/2\bbZ\right)^2$ :
\begin{align*} \bar X_{1,1}&= \bar X((4,3,2,1,0,-1,-2,2,1,0,-1,-2,1,0) , (4,3,2,1,0,-1,-2,2,1,0,-1,-2,1,0))\\
\bar X_{1,-1}&= \bar X((4,3,2,1,0,-1,-2,2,1,0,-1,-2,1,0) , (4,3,2,1,0,-1,-2,2,1,0,-1,-2,0,-1))\\
\bar X_{-1,1}&= \bar X((4,3,2,1,0,-1,-2,2,1,0,-1,-2,1,0) , (4,2,1,0,-1,-2,-3,2,1,0,-1,-2,1,0))\\
\bar X_{-1,-1}&= \bar X((4,3,2,1,0,-1,-2,2,1,0,-1,-2,1,0) , (4,2,1,0,-1,-2,-3,2,1,0,-1,-2,0,-1))\\
\end{align*}

\begin{exemple} \label{castrian}
Traitons l'exemple  des \og cas triangulaires \fg\, ({\sl cf.} \cite{BV}, \S 9).
Il s'agit des cas où la partition $d_{\caU}$ est de la forme   
$$d_{\caU}=[2m+1, 2m-1,2m-1, \ldots, 3,3, 1,1].$$
L'orbite duale $\caO$ est aussi triangulaire, avec  $d_{\caO}=[2m, 2m, \ldots, 4,4, 2,2]$, et le symbole
(le même pour $\caU$ et $\caO$, ce qui est très particulier) est
\begin{equation*}
 \begin{pmatrix}
  0 &&  2 &&  \ldots  && \ldots && 2m\\
& 1  && 3  && \ldots  &&  2m-1 &
\end{pmatrix} .
\end{equation*}
On a alors $\frmm= \prod_{j=1}^m \frg\frl(2j)$, $\bar A(\caU)\simeq \left(\bbZ/2\bbZ \right)^m$
et \[ \scrF^j_1=\chi_{\frac{1}{2},\frac{1}{2} } \circ \textstyle \det_{2j},\quad 
 \scrF^j_{-1}=\chi_{\frac{1}{2},-\frac{1}{2} } \circ \textstyle \det_{2j}.\]
Les induites (\ref{defpaq}) sont dans ce cas irréductibles (\cite{BV}, Prop. 9.11) 
(rappelons que pour $\alpha,\beta \in \bbC$ avec $\alpha-\beta \in \bbZ$, 
on a noté $\chi_{\alpha,\beta}$ le caractère $z\mapsto z^{\alpha}\bar z^{\beta}$ de $\bbC^\times$).
\end{exemple}
 
\begin{exemple}
$d_{\caU}=[7, 3,3]$,  $\lambda=(3,2,1,1,1,0)$, le symbole est 
$\begin{pmatrix}
  1 &&  4 \\
& 2  &
\end{pmatrix} $, $\bar A({}^L\caA)\simeq Z/2\bbZ$,
$\frmm=\frs\frp_2(\bbC)\times \frg\frl_5(\bbC)$, $F^1$ a pour caractère infinitésimal $(3,2,1,0,-1)$
et $\tilde F^1$ a pour caractère infinitésimal $(3,1,0,-1,-2)$, $\lambda^0=(1)$, $ \bar X(\lambda^0,\lambda^0)
=\Triv_{\Sp_2(\bbC)}$.
On a deux éléments dans le paquet, qui sont respectivement les sous-quotients irréductibles contenant les $K$-types
minimaux de 
\[ \Ind_P^G\left( \Triv_{\Sp_2(\bbC)}  \otimes \scrF^1_{\pm 1} \right).   \]
Où encore, en les  indexant par les caractères du groupes $\bar A(\caU)=\bbZ/2\bbZ$ :
\begin{align*} \bar X_{1}&= \bar X((3,2,1,0,-1,1) , (3,2,1,0,-1,1))\\
\bar X_{-1}&= \bar X((3,2,1,0,-1,1) , (3,1,0,-1,-2,1)).\\
\end{align*}
\end{exemple}

\medskip 
{\bf Cas $\mathbf{B}_n$, bonne parité}.
 Nous supposons que $\caU$ est de bonne parité, {\sl i.e.} les blocs sont tous pairs.
  Nous illustrons ce cas par l'exemple  $\mathbf{d}_{\caU}=[8,4,4,4,2]$, de rang $11$.
 Le symbole est encore de la forme (\ref{SymbBC}) et dans notre exemple, il est donné par 
 $\begin{pmatrix}
  1 &&  3&&6 \\
& 2  &&3 &
\end{pmatrix} $
 
On définit  une sous-algèbre $\frmm$ de $\frg=\frs \fro_{2n+1}$ de la manière suivante
\begin{equation}\label{defm}  \frmm \simeq \frmm^0  \times \caA^1 \times \ldots \times \caA^k    \end{equation}
 où $\frmm^0 \simeq \frs \fro_{2(x_{2k}-k)+1}$, 
 $\caA^j\simeq \frg \frl (x_{2k-2j+1}+x_{2k-2j}-2(k-j) )$, $1\leq j\leq k$.
Dans notre exemple,  $\frmm^0 \simeq \frs \fro_9(\bbC)\times \frg\frl_4(\bbC)\times  \frg\frl_3(\bbC)$.

Comme dans le cas {$\mathbf C_n$},  on écrit le caractère infinitésimal attaché à  $\caU$ en coordonnées:
\[ \lambda=(\lambda_1, \lambda_1, \ldots, \lambda_1,\lambda_2, \ldots, \lambda_2,   
 \ldots ,\lambda_s, \ldots, \lambda_s) \]
 avec $\lambda_i>\lambda_{i+1}$, $\lambda_s \geq 0$.
Les $\lambda_i$ sont des demi-entiers dans ce cas. 
 
Dans notre exemple, on obtient 
 \[ \lambda=\left(\frac{7}{2},\frac{5}{2},\frac{3}{2},\frac{3}{2},\frac{3}{2},\frac{3}{2},
 \frac{1}{2},\frac{1}{2},\frac{1}{2},\frac{1}{2} ,  \frac{1}{2}  \right). \]
  
 On extrait alors de $\lambda$ la suite  strictement décroissante maximale,
 $\lambda^0=(\lambda_1,\lambda_2, \ldots ,\lambda_s)$.
  Dans notre exemple, $\lambda^0=\left(\frac{7}{2},\frac{5}{2},\frac{3}{2},\frac{1}{2}\right)$.
   On définit les $F^j$, représentations de dimension
  finie holomorphe de $\caA^j$ à partir de ce qui reste de $\lambda$ comme dans le cas {\bf  A}.
   Dans notre exemple $F^1$ est  la représentation de dimension finie 
 holomorphe de $\frg \frl(4)$ de caractère infinitésimal $\left(\frac{3}{2},\frac{1}{2},  -\frac{1}{2},
  -\frac{3}{ 2}\right)$ (c'est la représentation triviale),
 et $F^2$  est  le caractère  holomorphe de $\frg \frl(3)$  de caractère infinitésimal 
  $\left(\frac{3}{2},\frac{1}{2},   -\frac{1}{2}\right)$.

   \medskip

Lorsque un entier $a$ apparaît avec multiplicité $2$ dans le symbole, on vérifie, en utilisant 
le fait que tous les blocs de $\mathbf{d}_{\caU}$ sont impairs, que c'est sous la forme suivante : 
\begin{equation*}
 \begin{pmatrix}
  x_0 &&  \ldots &  x_{2j}=a &&& \ldots && x_{2k}\\
& \ldots  &&& x_{2j+1}=a  & \ldots  &&  x_{2k-1} &
\end{pmatrix} .
\end{equation*}

Lorsque l'on enlève ces paires $(x_{2j},x_{2j+1})$ avec $x_{2j}=x_{2j+1}$ du symbole, il nous reste donc 
$m$ paires   $(x_{2j},x_{2j+1})$ avec $x_{2j+1}>x_{2j}$. 
Pour une telle paire, on définit la représentation de dimension finie 
 holomorphe $\tilde F^j$ de $\caA^j$ de la m\^eme manière que dans le cas {$\mathbf C_n$}.
   Dans notre exemple, on a $m=1$, $\tilde F^2$ a pour caractère infinitésimal 
   $\left(\frac{1}{2},-\frac{1}{2},-\frac{3}{2}\right)$.

On définit  alors, pour tout $1\leq j\leq k$ tel que   $x_{2j}=x_{2j+1}$, les représentations de dimension finie
$ \scrF_1^j=F^j \otimes F^j$, 
 pour tout $1\leq j\leq k$  tel que   $x_{2j}<x_{2j+1}$,  
$\scrF^j_{1}= F^ j \otimes F^j$, $ \scrF^j_{-1}= F^j \otimes \tilde F^j$, et   
$\scrF_{\hat{\boldsymbol{\epsilon}}}=\otimes_j\scrF^j_{\hat\epsilon_j}$ comme en (\ref{caF}) (\ref{caF2}) et 
(\ref{cafe}).

Soit $P=MN$ un sous-groupe parabolique de $G$ dont le facteur de Levi a pour algèbre de Lie $\frmm$.
C'est un produit  des groupes généraux linéaires d'algèbre de Lie  $\caA^j$, $j=1, \ldots ,k$, et 
d'un facteur $G(x_{2k}-k)=\SO_{2(x_{2k}-k)+1}(\bbC)$. Considérons la représentation sphérique 
$\bar X(\lambda^0,\lambda^0)$ de ce groupe\footnote{C'est la représentation triviale.},
 et définissons
$ \bar X_{\hat{\boldsymbol{\epsilon}}}$  comme en (\ref{defpaq}) 
Le paquet associé à $\psi_G$ par Barbasch et Vogan est alors 
\begin{equation} \label{defpaq2casB} \Pi_{BV}(\caU,G)=\{  \bar X_{\hat{\boldsymbol{\epsilon}}} \}.
 \end{equation}
et la remarque \ref{parpaq} est encore valide.

\medskip

Revenons à notre exemple, où l'on a  
\[   \bar X(\lambda^0,\lambda^0)=\Triv_{\SO_9(\bbC)}, \;      \scrF^1_1=\Triv_{\GL_4}, 
 \; \scrF^2_1=\chi_{\frac{1}{2},\frac{1}{2} }\circ \textstyle \det_2, 
\; \scrF^2_{-1}=\chi_{\frac{1}{2},-\frac{1}{2} }\circ \textstyle \det_2. \]
On obtient donc dans ce cas un paquet à 2 éléments : 
\[ \bar X_1=\bar X(  (\frac{7}{2},\frac{5}{2},\frac{3}{2},\frac{1}{2},\frac{3}{2},\frac{1}{2},  
-\frac{1}{2}, -\frac{3}{2},\frac{3}{2},\frac{1}{2},  -
\frac{1}{2} ),(\frac{7}{2},\frac{5}{2},\frac{3}{2},\frac{1}{2},\frac{3}{2},\frac{1}{2},  -\frac{1}{2}, 
-\frac{3}{2},\frac{3}{2},\frac{1}{2},  -\frac{1}{2} )   ),  \] 
\[ \bar X_{-1}=\bar X(  (\frac{7}{2},\frac{5}{2},\frac{3}{2},\frac{1}{2},\frac{3}{2},\frac{1}{2},  
-\frac{1}{2}, -\frac{3}{2},\frac{3}{2},\frac{1}{2},  -
\frac{1}{2} ),(\frac{7}{2},\frac{5}{2},\frac{3}{2},\frac{1}{2},\frac{3}{2},\frac{1}{2},  -\frac{1}{2}, 
-\frac{3}{2},\frac{1}{2},-\frac{1}{2},  -\frac{3}{2} )   ).  \]

 \begin{exemple} \label{castrianB}
Le cas triangulaire est ici celui  où la partition $d_{\caU}$ est de la forme   
$$d_{\caU}=[2m, 2m,2m-2,2m-2 \ldots, 4,4, 2,2].$$
Le symbole est
\begin{equation*}
 \begin{pmatrix}
  0 &&  2 &&  \ldots  && \ldots && 2m\\
& 1  && 3  && \ldots  &&  2m-1 &
\end{pmatrix} .
\end{equation*}
On a alors $\frmm=\frs\fro_{2m+1}(\bbC)\times  \prod_{j=1}^m \frg\frl(2j-1)$, $\bar A(\caU)\simeq \left(\bbZ/2\bbZ \right)^m$
$  \bar X(\lambda^0,\lambda^0)=\Triv_{\SO_{2m+1}(\bbC)}$, 
et \[ \scrF^j_1=\chi_{\frac{1}{2},\frac{1}{2} } \circ \textstyle \det_{2j-1},\quad 
 \scrF^j_{-1}=\chi_{\frac{1}{2},-\frac{1}{2} } \circ \textstyle \det_{2j-1}.\]
Les induites (\ref{defpaq}) sont dans ce cas irréductibles (\cite{BV}, Prop. 9.11).
\end{exemple}

\begin{exemple}
$d_{\caU}=[6, 4,2]$, $d_{\caO}=[3, 3, 3,1,1,1]$, $\lambda=\left(\frac{5}{2},\frac{3}{2}, \frac{3}{2}, 
 \frac{1}{2},\frac{1}{2}, \frac{1}{2},   \right)$, le symbole est 
$\begin{pmatrix}
  1 &&  4 \\
& 2  &
\end{pmatrix} $, $\bar A(\caU)\simeq Z/2\bbZ$,
$\frmm=\frs\fro_7(\bbC)\times \frg\frl_3(\bbC)$, $\lambda^0=\left(\frac{5}{2},\frac{3}{2}, \frac{1}{2}\right)$,
 $F^1$ a pour caractère infinitésimal $\left(\frac{3}{2},\frac{1}{2}, -\frac{1}{2}\right)$
et $\tilde F^1$ a pour caractère infinitésimal $ \left(\frac{1}{2},-\frac{1}{2}, -\frac{3}{2}\right)$, 
$\bar X(\lambda^0,\lambda^0)=\Triv_{\SO_7(\bbC)}$ et  
on a deux éléments dans le paquet, qui sont respectivement les sous-quotient irréductibles contenant les $K$-types
minimaux de 
\[ \Ind_P^G\left(  \Triv_{\SO_7(\bbC)} \otimes \scrF^1_{\pm 1} \right),   \]
où encore, 
\[ \bar X_1=\bar X\left(  \left(\frac{5}{2},\frac{3}{2},\frac{1}{2},\frac{3}{2},\frac{1}{2},-\frac{1}{2}\right),
\left( \frac{5}{2},\frac{3}{2},\frac{1}{2},\frac{3}{2},\frac{1}{2},-\frac{1}{2} \right) \right) ,  \] 
\[ \bar X_{-1}=\bar X\left ( \left(\frac{5}{2},\frac{3}{2},\frac{1}{2},\frac{3}{2},\frac{1}{2},-\frac{1}{2}\right), 
\left(\frac{5}{2},\frac{3}{2},\frac{1}{2},\frac{1}{2},-\frac{1}{2},-\frac{3}{2}\right) \right).  \]

\end{exemple}

\medskip
{\bf Cas $\mathbf{B}_n$, mauvaise parité}.
  Nous supposons que $\caU$ est de mauvaise parité, {\sl i.e.} les blocs sont tous 
impairs, et ont une multiplicité paire : 
$\mathbf{d}_{\caU}=[2d_1+1,2d_1+1,\ldots, 2d_t+1,2d_t+1]$.
Le symbole est encore de la forme \begin{equation*}
 \begin{pmatrix}
  0 &&  d_t+1 && \ldots && d_1+t\\
& d_t+1  && \ldots && d_1+t   &
\end{pmatrix} 
\end{equation*}
 On définit  une sous-algèbre $\frmm$ de $\frg=\frs \fro_{2n+1}$ de la manière suivante
\begin{equation}\label{defmcasBmp}  \frmm \simeq  \caA^1 \times \ldots \times \caA^k    \end{equation}
 où $\caA^j=\frg\frl(2d_j+1)$, $1\leq j\leq k$
 comme en (\ref{defmcasA}).
 Soit $P=MN$ un sous-groupe  parabolique de $G$ dont le facteur de Levi a pour algèbre de Lie 
$\frmm$, c'est-à-dire
$M\simeq \GL_{2d_1+1}(\bbC)\times \cdots \times \GL_{2d_t+1}(\bbC)$.
Le paquet $\Pi_{BV}(\caU,G)$ est un singleton,  
l'induite de  $P$ à $G$ de la représentation triviale (qui est irréductible).

\bigskip {\bf Cas $\mathbf{D}_n$, bonne parité}. Nous supposons que $\caU$ est de bonne parité, {\sl i.e.} 
les blocs sont tous impairs.
  Nous illustrons ce cas par l'exemple  $\mathbf{d}_{\caU}=[11,7,7,7,3,3]$, de rang 19.
 En général, le symbole est de la forme
 \begin{equation} \label{SymbD}
 \begin{pmatrix} x_0 &  & x_2 & & \ldots &&&& x_{2k}\\
x_1 && x_3 && \ldots  &&  && x_{2k+1}
\end{pmatrix} 
\end{equation}

Dans notre exemple, on obtient 
\begin{equation} \label{SymbExD}
 \begin{pmatrix} 1&  & 4   && 5 \\
2 && 5 && 8   
\end{pmatrix} 
\end{equation}

On définit  une sous-algèbre $\frmm$ de $\frg=\frs \fro_{2n+1}$ de la manière suivante
\begin{equation}\label{defmC}  \frmm \simeq \frmm^0  \times \caA^1 \times \ldots \times \caA^k    \end{equation}
 où $\frmm^0 \simeq \frs \fro_{2(x_{2k+1}+x_0-k)}$, 
 $\caA^j\simeq \frg \frl (x_{2k-2j+2}+x_{2k-2j+1}-2(k-j)-1 )$, $1\leq j\leq k$.

On écrit le caractère infinitésimal attaché à  $\caU$ en coordonnées:
\[ \lambda=(\lambda_1, \lambda_1, \ldots, \lambda_1,\lambda_2, \ldots, \lambda_2,   
 \ldots ,\lambda_s, \ldots, \lambda_s) \]
 avec $\lambda_i>\lambda_{i+1}$, $\lambda_s \geq 0$.
Les $\lambda_i$ sont des entiers dans ce cas.

Dans notre exemple, on obtient 
 \[ \lambda=\left(5,4,3,3,3,3,2,2,2,2,1,1,1,1,1,1,0,0,0   \right). \]
  
 On extrait alors de $\lambda$ la suite  strictement décroissante maximale,
 $\tilde \lambda^0=(\lambda_1,\lambda_2, \ldots ,\lambda_s)$.
  Dans notre exemple, $\tilde \lambda^0=\left(5,4,3,2,1,0\right)$. On définit les $F^j$, représentations de dimension
  finie holomorphe de $\caA^j$ à partir de ce qui reste de $\lambda$ comme dans le cas {\bf  A}.
   Dans notre exemple $F^1$ est  la représentation de dimension finie 
 holomorphe de $\frg \frl(7)$ de caractère infinitésimal $\left(3,2,1,0,-1,-2,-3\right)$ (la triviale)
et $F^2$  est  le caractère  holomorphe de $\frg \frl(5)$  de caractère infinitésimal $ \left(3,2,1,0, -1 \right)$.

Lorsque un entier $a$ apparaît avec multiplicité $2$ dans le symbole, on vérifie, en utilisant 
le fait que tous les blocs de $\mathbf{d}_{\caU}$ sont impairs, que c'est sous la forme suivante : 
\begin{equation*}
 \begin{pmatrix}
  x_0 &&  \ldots &&&  x_{2j}=a && \ldots & x_{2k}\\
x_1 && \ldots  && x_{2j-1}=a  & \ldots  &&&  x_{2k+1} &
\end{pmatrix} 
\end{equation*}

Lorsque l'on enlève ces paires $(x_{2j-1},x_{2j})$ avec $x_{2j-1}=x_{2j}$ du symbole, il nous reste donc 
$m$ paires   $(x_{2j}-1,x_{2j})$ avec $x_{2j-1}<x_{2j}$. 
Pour une telle paire, on définit représentation de dimension finie 
 holomorphe $\tilde F^j$ de $\caA^j$ de la m\^eme manière que dans le cas {$\mathbf C_n$}.
   Dans notre exemple, on a $m=1$, $\tilde F^2$ a pour caractère infinitésimal 
   $\left(3,1,0,-1,-2 \right)$.

On définit  alors, pour tout $1\leq j\leq k$ tel que   $x_{2j-1}=x_{2j}$, les représentations de dimension finie
$ \scrF_1^j=F^j \otimes F^j$, 
 pour tout $1\leq j\leq k$  tel que   $x_{2j-1}<x_{2j}$,  
$\scrF^j_{1}= F^ j \otimes F^j$, $ \scrF^j_{-1}= F^j \otimes \tilde F^j$, et   
$\scrF_{\hat{\boldsymbol{\epsilon}}}=\otimes_j\scrF^j_{\hat\epsilon_j}$ comme en (\ref{caF}) (\ref{caF2}) et 
(\ref{cafe}). 

On complète $\tilde \lambda^0$ pour former $\lambda^0$  en ajoutant les opposés de ce qui reste des termes de 
$\lambda$ (dans notre exemple, il ne reste que $(1)$, donc   $\lambda^0=\left( 5,4,3,2,1,0,-1\right)$.

Soit $P=MN$ un sous-groupe parabolique de $G$ dont le facteur de Levi a pour algèbre de Lie $\frmm$.
C'est un produit  des groupes généraux linéaires d'algèbre de Lie  $\caA^j$, $j=1, \ldots ,k$, et 
d'un facteur $G(x_{2k+1}+x_0-k)=\SO_{2(x_{2k+1}+x_0-k)}(\bbC)$. Considérons la représentation sphérique 
$\bar X(\lambda^0,\lambda^0)$ de ce groupe, et définissons
$ \bar X_{\hat{\boldsymbol{\epsilon}}}$  comme en (\ref{defpaq}) 
Le paquet associé à $\caU$ par Barbasch et Vogan est alors 
\begin{equation} \label{defpaq2casC} \Pi_{BV}(\caU,G)=\{  \bar X_{\hat{\boldsymbol{\epsilon}}} \}.
 \end{equation}
et la remarque \ref{parpaq} est encore valide.

\medskip

Revenons à notre exemple, où le paquet compte donc 2 éléments :
\[ \bar X_{1}= \bar X(\lambda,\lambda), \quad \bar X_{-1}= \bar X( \lambda,\mu),\]
avec $\lambda= (5,4,3,2,1,0,-1,3,2,1,0,-1,-2,-3,3,2,1,0,-1)$ 

\noindent et  $ \mu=(5,4,3,2,1,0,-1,3,2,1,0,-1,-2,-3,3,1,0,-1,-2)$.

\begin{exemple} 
$d_{\caU}=[9, 5,5,3]$,  $\lambda=\left(4,3,2,2,2,1,1,1,1,0,0  \right)$, 
le symbole est 
$\begin{pmatrix}
  1 &&  3 \\
3&& 6 
\end{pmatrix} $, $\bar A(\caU)$ est trivial,
$\frmm=\frs\fro_{12}(\bbC)\times \frg\frl_5(\bbC)$, $\tilde \lambda^0=\left(4,3,2,1,0 \right)$,
 $F^1$ a pour caractère infinitésimal $(2,1,0,-1,-2 )$ 
 (c'est la représentation triviale de $\GL_5(\bbC)$),    $\lambda^0=\left(4,3,2,1,0,-1 \right)$,
on a un seul  élément dans le paquet, le sous-quotient irréductible contenant le $K$-type
minimal de $ \Ind_P^G\left(   \bar X(\lambda^0,\lambda^0) \otimes \Triv_{\GL_5(\bbC)} \right)$, c'est-à-dire :
\[\bar X=\bar X( (4,3,2,1,0,2,1,0,-1,-2),  (4,3,2,1,0,2,1,0,-1,-2)).  \]
\end{exemple}

 \begin{exemple} 
$d_{\caU}=[11,9, 5,5,3,1]$, $\lambda=\left(5,4,4,3,3,2,2,2,2,1,1,1,1,1,0,0,0  \right)$, 
le symbole est 
$\begin{pmatrix}
  0 &&  3&&6 \\
2&& 4 &&8 
\end{pmatrix} $, $\bar A(\caU)\simeq \left( \bbZ/2\bbZ\right)^2$,
$\frmm=\frs\fro_{12}(\bbC)\times \frg\frl_5(\bbC)$, $\tilde \lambda^0=\left(5,4,3,2,1,0 \right)$,
 $F^1$ a pour caractère infinitésimal $(4,3,2,1,0,-1,-2 )$, 
  $\tilde F^1$ a pour caractère infinitésimal $(4,2,1,0,-1,-2,-3 )$,  
 $F^2$ a pour caractère infinitésimal $(2,1,0,-1,)$, 
  $\tilde F^2$ a pour caractère infinitésimal $(1,0,-1,-2)$, 
     $\lambda^0=\left(5,4,3,2,1,0\right)=\tilde \lambda^0$.
     Il y a  a  quatre   éléments dans le paquet.
{\small
\begin{align*} \bar X_{1,1}&=\bar X((5,4,3,2,1,0 , 4,3,2,1,0,-1,-2 ,2,1,0,-1) ,
(5,4,3,2,1,0 , 4,3,2,1,0,-1,-2 ,2,1,0,-1)  )\\
\bar X_{1,-1}&=\bar X((5,4,3,2,1,0 , 4,3,2,1,0,-1,-2 ,2,1,0) ,(5,4,3,2,1,0 , 4,3,2,1,0,-1,-2 ,1,0,-1,-2)  )\\
 \bar X_{-1,1}&=\bar X((5,4,3,2,1,0 , 4,3,2,1,0,-1,-2 ,2,1,0,-1) ,(5,4,3,2,1,0 , 4,2,1,0,-1,-2,-3 ,2,1,0,-1)  )\\
\bar X_{-1,1}&=\bar X((5,4,3,2,1,0 , 4,3,2,1,0,-1,-2 ,2,1,0,-1) ,(5,4,3,2,1,0 , 4,2,1,0,-1,-2,-3 ,1,0,-1,-2)  ).\\
\end{align*} }
\end{exemple}

\medskip 
{\bf Cas $\mathbf{D}_n$, mauvaise parité}.
  Nous supposons que $\caU$ est de mauvaise parité, {\sl i.e.} les blocs sont tous 
pairs, et ont une multiplicité paire.
$\mathbf{d}_{\caU}=[2d_1,2d_1,\ldots, 2d_t,2d_t]$.
 Le symbole est  de la forme \begin{equation*}
 \begin{pmatrix}
   d_t && \ldots && d_1+t-1\\
 d_t  && \ldots && d_1+t -1  &
\end{pmatrix} 
\end{equation*}
 On définit  une sous-algèbre $\frmm$ de $\frg=\frs \fro_{2n}$ de la manière suivante
\begin{equation}\label{defmcasDmp}  \frmm \simeq  \caA^1 \times \ldots \times \caA^t    \end{equation}
 où $\caA^j=\frg\frl(2d_j)$, $1\leq j\leq t$.
 Soit $P=MN$ un sous-groupe parabolique de $G$ dont le facteur de Levi a pour algèbre de Lie $\frmm$, c'est-à-dire
$M\simeq \GL_{2d_1}(\bbC)\times \cdots \times \GL_{2d_t}(\bbC)$.
Le paquet $\Pi_{BV}(\caU,G)$ est un singleton,  
l'induite de  $P$ à $G$ de la représentation triviale (qui est irréductible).

\bigskip

{\bf Cas $\mathbf{Mp}_n$}.
 Nous supposons que $\caU$ est de bonne parité, {\sl i.e.} les blocs sont tous pairs.
  Nous illustrons ce cas par l'exemple  $\mathbf{d}_{\caU}=[12, 8,4,4,4,2]$, de rang $17$.
La règle de calcul pour le symbole change : posons $\mathbf{d}_\caU=[d_1, \ldots d_t]$. Quitte 
à ajouter un bloc $(0)$ on suppose que $t$ est impair. On forme la suite strictement croissante
\[(d_t+1, d_{t-1}+2, \ldots , d_1+t)   \]
 Notons $(2x_0+1, \ldots , 2x_{2k}+1)$ les  termes impairs de cette suite et  $(2x_1, \ldots , 2x_{2k-1})$ 
  les termes pairs.  Le symbole est alors
    \begin{equation} \label{symbMp}
 \begin{pmatrix} x_0 &  & x_2 & & \ldots &&&& x_{2k}\\
&x_1 && x_3 && \ldots  && x_{2k-1}
\end{pmatrix} \end{equation}
  et dans notre exemple, il est donné par 
 $\begin{pmatrix}
  0 &&  3&&4&&9 \\
& 2  &&4 &&7&
\end{pmatrix} $

A partir de ce symbole, on retrouve $\bar A^{Mp}(\caU)$ par une recette très similaire à celle du cas $\mathbf{C}_n$:
chaque couple $(x_{2j-1},  x_{2j})$ avec $x_{2j-1}<  x_{2j}$ contribue d'un facteur $\bbZ/2\bbZ$, mais en plus, 
si $x_0=d_1/2\neq 0$, on a un facteur $\bbZ/2\bbZ$ supplémentaire.
Dans notre exemple, c'est donc $(\bbZ/2\bbZ)^2$.
 La règle qui donne le caractère infinitésimal est la même que dans les autres cas, dans l'exemple on obtient
$$\left(\frac{11}{2},\frac{9}{2},\frac{7}{2},\frac{7}{2},\frac{5}{2},\frac{5}{2},\frac{3}{2} ,\frac{3}{2} ,
\frac{3}{2}  ,\frac{3}{2}  ,\frac{3}{2}, \frac{3}{2}, 
  ,\frac{1}{2}  ,\frac{1}{2}  ,\frac{1}{2}  ,\frac{1}{2} ,\frac{1}{2}\right).$$
Le reste se calcule à partir de ce caractère infinitésimal en suivant aussi la recette du cas $\mathbf{C}_n$, 
mais si $x_0=d_1/2$ n'est pas nul, 
il faut remplacer le facteur $\Sp_{2x_0}(\bbC)$ par un facteur $\Mp_{2x_0}(\bbC)$ , 
et la représentation triviale de celui-là par les deux représentations métaplectiques.
Dans notre exemple, on commence par extraire
 $\left(\frac{11}{2},\frac{9}{2},\frac{7}{2},\frac{5}{2},\frac{3}{2} ,\frac{1}{2} ,   \right)$ 
 puis on remonte en prenant les opposés, ce qui donne une représentation holomorphe 
 $F^1$ de caractère infinitésimal  $\left(\frac{11}{2},\frac{9}{2},\frac{7}{2},\frac{5}{2},\frac{3}{2} ,
 \frac{1}{2} , - \frac{1}{2},  -\frac{3}{2},  - \frac{5}{2}, - \frac{7}{2} \right)$ d'un facteur 
 $\GL_{10}(\bbC)$, une représentation $\tilde F^1$ de caractère infinitésimal 
  $\left(\frac{11}{2},\frac{7}{2},\frac{5}{2},\frac{3}{2} ,
 \frac{1}{2} , - \frac{1}{2},  -\frac{3}{2},  - \frac{5}{2}, - \frac{7}{2} , -\frac{9}{2}\right)$ 
et les deux représentations $\scrF_1^1=F^1\otimes F^1$ et $\scrF_{-1}^1=F^1\otimes \tilde F^1$.
Ensuite, on extrait  $\left(\frac{3}{2} ,\frac{1}{2}    \right)$ 
 puis on remonte en prenant les opposés, ce qui donne une représentation holomorphe 
 $F^2$ de caractère infinitésimal 
  $\left(\frac{3}{2} ,\frac{1}{2} , - \frac{1}{2},  -\frac{3}{2},  \right)$ d'un facteur 
 $\GL_{4}(\bbC)$, et la représentation $\scrF^2= F^ 2\otimes F^2=\Triv_{\GL_4(\bbC)}$.
   Enfin, on extrait $\left(\frac{3}{2} ,\frac{1}{2}    \right)$ 
 puis on remonte en prenant les opposés, ce qui donne une représentation holomorphe 
 $F^3$ de caractère infinitésimal 
  $\left(\frac{3}{2} ,\frac{1}{2} , - \frac{1}{2}  \right)$ de 
 $\GL_{3}(\bbC)$,  une représentation holomorphe 
 $\tilde F^3$ de caractère infinitésimal 
  $\left(\frac{1}{2} , - \frac{1}{2} - \frac{3}{2}, \right)$ et 
les deux représentations $\scrF_1^3=F^3\otimes F^3$ et $\scrF_{-1}^3=F^3\otimes \tilde F^3$.
On obtient $4$ représentations dans notre paquet, les sous-quotients de Langlands des induites
d'un sous-groupe parabolique $P$ de facteur de Levi $M\simeq \GL_{10}(\bbC)\times   \GL_4(\bbC)\times \GL_3(\bbC)$ 
\[ \Ind_P^G(\scrF^1_{\pm}\otimes \scrF^2\otimes \scrF^ 3_{\pm} ). \]

Une autre manière de comprendre  cette recette est de voir que l'on retrouve bien la construction des paquets
par la correspondence de Howe ({\sl cf.} \cite{pourhowe}). En effet, partons d'une orbite 
$\caU$ de partition $\mathbf{d}_\caU$ et considérons l'orbite $\caU^-$ obtenue en enlevant le plus grand 
bloc $d_1$ (pair). Calculons le paquet de $\SO_{2n-d_1+1}(\bbC)$ correspondant à cette 
orbite $\caU^-$. Dans notre cas,  $\mathbf{d}_{\caU^ -}=[8,4,4,4,2]$, et c'est l'exemple du cas $\mathbf{B}_n$.
On avait un paquet  à 2 éléments, car $\bar A(\caU^-)=\bbZ/2\bbZ$.
Si $d_1 >d_2$ chacune des  représentations du paquet de $\SO_{2n-d_1}(\bbC)$ 
 s'étend  de deux manières différente en une représentation de 
$\Or_{2n-d_1+1}(\bbC)$, et  ces deux  représentations ont une image par la correspondance de Howe entre 
$\Or_{2n-d_1+1}(\bbC)$ et $\Mp_{2n}(\bbC)$. On voit en comparant les symboles
que  l'on a bien ajouté un facteur $\bbZ/2\bbZ$ en passant de $\caU^-$ à $\caU$. 
C'est le cas dans notre exemple.
En revanche, si $d_1=d_2$, une seule des deux représentations étendues à $\Or_{2n-d_1+1}(\bbC)$
possède  une image par la correspondance de Howe, et il n'y a pas de facteur $\bbZ/2\bbZ$ supplémentaire.
Par exemple, pour l'orbite $[8,8,4,4,4,2]$, on obtient le symbole
$\begin{pmatrix}
  0 &&  3&&4&&7 \\
& 2  &&4 &&7&
\end{pmatrix} $
et  $\bar A^{Mp}(\caU)\simeq \bbZ/2\bbZ$.

Dans les exemples donnés ci-dessus, nous sommes partis d'une partition ayant un nombre pair de termes, il a donc fallu
ajouter un bloc $0$ pour calculer le symbole, d'où un $x_0=0$. Voyons ce qui se passe lorsque 
la partition de départ à un nombre impair de termes. Disons $\mathbf{d}_{\caU}=[12, 8, 4,4,2]$.
Le symbole est alors 
$\begin{pmatrix}
  1 &&  3&&9 \\
& 3  &&6 &
\end{pmatrix} $, et $\bar A^{Mp}(\caU)=(\bbZ/2\bbZ)^2$.
Le $1$ dans ce symbole contribue à un facteur $\bbZ/2\bbZ$, et dans la construction des représentations
du paquets comme sous-quotient de Langlands d'induites paraboliques, il donne un facteur $\Mp_2(\bbC)$ 
au sous-groupe de Levi du parabolique. Sur ce facteur, nous mettons les deux représentations
métaplectiques $\bar X((1/2),(1/2))$ et  $\bar X((1/2),(-1/2))$.
Si on applique la recette avec la correspondence de Howe, on part de l'orbite  
$\mathbf{d}_{\caU}=[ 8, 4,4,2]$ qui a pour symbole 
$\begin{pmatrix}
  0 &&  3&&6 \\
& 1  &&3 &
\end{pmatrix} $, et $\bar A(\caU^-)=\bbZ/2\bbZ$ comme il se doit.

\section{Lemmes de réduction}\label{LemmeRed}

Dans cette section, on étudie l'effet de l'induction parabolique. Plus
précisément, on fixe un parabolique maximal de $G$ dont le Levi est
isomorphe à un groupe $\GL_a(\bbC)\times G^-$  où $G^-$ est de même type que $G$. On 
considère les induites du caractère trivial de  $\GL_a(\bbC)\times G^-$   avec une
représentation spéciale unipotente de $G^-$. Ces induites ont déjà été
considérées dans \cite{BV} et on traduit les résultats de {\sl loc. cit.} en
termes combinatoires.

Commençons par introduire quelques notations commodes.
   On note $\{0\}_\frg$ l'orbite adjointe nilpotente $\{0\}$ dans l'algèbre de Lie $\frg$.
Si $\frg=\frg\frl_n(\bbC)$, on la note plus simplement $\{0\}_n$.

Soient $G$  un groupe   classique complexe, et  $P=MN$  un sous-groupe 
parabolique de  $G$ de facteur de Levi $M$ isomorphe à $\left(\times_i \GL_{n_i}(\bbC) \right)\times G^-$, 
où $G^-$ est un groupe    classique de même type que $G$. Soit $\caO^-$ une orbite nilpotente dans $\frg^-$. 
On note  $ \left( \star_i \{0\}_{n_i} \right)\star \caO^-$ l'orbite induite (voir \cite{LSp} ou \cite{BV},
(4.13))  de $\frmm$ à $\frg$ de $\left( \times_i  \{0\}_{n_i} \right)\times \caO^-$.

Nous allons utiliser le résultat suivant, qui est une version plus précise de la proposition 12.5 de \cite{BV}
et de la proposition 6.6 de \cite{B1}.

\begin{prop}\label{PropB166}
Soient $G$  un groupe   classique complexe de rang $n$, 
 $\caU$  une orbite nilpotente de bonne parité de ${}^L\frg$,  $P=MN$  un sous-groupe parabolique maximal
  de $G$ de facteur de Levi  $M$ isomorphe à $ \GL_{n_1}(\bbC)\times G^-$, 
  où $G^-$ est un groupe classique de même type que 
$G$ de rang $n-n_1$, $\caU^-$ une orbite  nilpotente spéciale paire de ${}^L\frg^-$.
Supposons que,  $\{0\}_{n_1}\times \caU^-$ soit contenue dans l'intersection de  $\caU$ et de $\frmm$.
Alors si l'on note $\caO$ et $\caO^-$ les orbites duales respectivement de  $\caU$ et 
 $\caU^-$, on a $\caO=\{0\}_{n_1}\star \caO^-$. Supposons de plus qu'il existe une représentation
 de dimension finie $\scrF$ de $\frg \frl_{n_1}(\bbC)$ tel que pour tout élément $\bar X^- \in \Pi_{BV}(\caU^-,G^-)$, 
 \[ \Ind_P^G (\scrF \otimes \bar X^-) \]
ait pour caractère infinitésimal celui des éléments de  $\Pi_{BV}(\caU,G)$.
Alors tout les sous-quotients irréductibles de $\Ind_P^G ( \scrF \otimes \bar X^-) $ sont dans 
 $\Pi_{BV}(\caU,G)$. 
 
 Dans le cas où $\caF$ est la représentation triviale $\Triv_{n_1}$,  qui est le cas qui nous intéresse, la représentation induite est semi-simple et nous pouvons décrire ses composantes irréductibles comme suit. 
 Rappelons que les éléments de  $\Pi_{BV}(\caU^-,G^-)$ sont paramétrés par les éléments de 
$\widehat{\bar A(\caU^-)}$, et de même pour  $\Pi_{BV}(\caU,G)$ et $\widehat{\bar A(\caU)} $.
On a une inclusion naturelle de $A(\caU^-)$ dans  $A(\caU)$.
Un caractère de  $\bar A(\caU)$ peut être vu comme un caractère de $A(\caU)$
trivial sur le noyau de la projection $A(\caU) \rightarrow \bar A(\caU)$, et de même 
pour  $\bar A(\caU^-)$ et  $ A(\caU^-)$. Quels que soient $\eta^-\in \widehat{ \bar A(\caU^-)} $
et  $\eta\in \widehat{\bar A(\caU)} $, notons $m(\eta^ -, \eta_{\vert A(\caU^-)})$ la multiplicité de 
$\eta^-$ dans la restriction de $\eta$ à $A(\caU^-)$  (ce sont des caractères, donc cette multiplicité
vaut $0$ ou $1$, et ceci exactement lorsque  la restriction de $\eta$ à $A(\caU^-)$ est égale à $\eta^-$).
Pour tout  $\eta \in \widehat{\bar A(\caU)} $ (resp.  $\eta^- \in \widehat{\bar A(\caU^-)} $),  
soit $\bar X_{\eta}$ (resp.  $\bar X_{\eta^-}^-$) l'élément de 
$\Pi_{BV}(\caU^-, G^-)$ (resp.  $\Pi_{BV}(\caU, G)$) correspondant. 
On a alors 
\begin{equation} \label{B166}
\Triv_{n_1}\star \bar X^-_{\eta^-} = \Ind_P^G (\Triv_{n_1} \otimes \bar X^-_{\eta^-})
=\bigoplus_{\eta\in  \widehat{\bar A(\caU^-)}} 
m(\eta^ -, \eta_{\vert A(\caU)}) \; \bar X_\eta.
\end{equation}
\end{prop}

\dem Par hypothèse,  pour tout $\bar X^-\in \Pi_{BV}(\caU^-,G^-) $, $X:= \Ind_P^G (\scrF\otimes \bar X^-)$, 
ainsi que tous ses sous-quotients irréductibles ont bien le caractère infinitésimal requis. Soit $\bar Y$ un  
de ses  sous-quotients irréductibles. Son front d'onde $\mathrm{WF}(\bar Y)$ est contenu dans celui de 
$X$, qui est $\overline{\caO}=\overline{\{0\}_{n_1} \star \caO^- }$ d'après \cite{BV}, (1.9).
Or la dimension du front d'onde de $\bar Y$ est au moins celle de $\caO$ d'après   \cite{BV}, Cor. 5.19, et comme
ce front d'onde est l'adhérence d'une orbite nilpotente, cette orbite ne peut être que $\caO$. On sait que 
$\{0\}_{n_1} \star \caO^-$ est l'orbite duale de ${\mathcal U}$ (cf. par exemple \cite{spaltenstein} 11.7 page 217, 
ou mieux un calcul à la main) et par hypothèse ${\mathcal U}$ est une orbite paire.  
Ceci montre que $\bar Y$ satisfait aux conditions de \cite{BV} pour être dans $\Pi_{BV}({\caU},G)$.

La formule (\ref{B166}) se démontre alors à partir des formules de caractères pour les représentations
spéciales unipotentes établies dans \cite{BV} comme dans la proposition 12.5 de cet article.\qed

\subsection{Effet de l'ajout de deux  blocs de même taille  de bonne parité, cas {$\mathbf C_n$}}
 \label{ajout}
 
On suppose que $G$ est un groupe symplectique de rang $n$. Soient $\caU$ une orbite nilpotente de bonne parité de 
${}^L\frg=\frs\fro_{2n+1}$  et  $\caO$ l'orbite duale dans $\frg$.
Soient $d_{\caU}$ et $d_\caO$ les partitions de $2n+1$ et $2n$ respectivement associées à ces orbites.

Les blocs de $d_{\caU}$ sont donc impairs, et il y en a un nombre impair on peut donc écrire 
\[ d_{\caU}={\scriptstyle [2 (x_{2p}-p)+1,  2(x_{2p-1}-p)+1,   \ldots,  2(x_{2p-2j+2}-p+j)-1, 
 2(x_{2p-2j+1}-p+j)-1    , \ldots, 
2x_2-1, 2x_1-1 , 2x_0+1] } \]
de sorte que le symbole est 
 \begin{equation}\label{Symbo}
 \begin{pmatrix}
  x_0 &&  \ldots &&  x_{2j} && \ldots && x_{2p}\\
& \ldots  && x_{2j-1}  && \ldots  &&  x_{2p-1} &
\end{pmatrix} 
\end{equation}
 
 \medskip 
 
 Ajoutons 2 blocs de taille $2M+1$ à $d_{\caU}$, pour former une nouvelle partition $d_{\caU^+}$,
 correspondant à une orbite nilpotente  $\caU^+$ dans ${}^L\frg^+=\frs\fro_{2(n+2M+1)+1}$.
 Soient $\caO^+$ sa duale dans $\frg^+=\frs \frp_{2(n+2M+1)}$. 
 On a $\caO^+=\{0\}_{2M+1} \star \caO$. On obtient $d_{\caO^+}$ en ajoutant 2 aux $2M+1$ premiers blocs
 de $d_\caO$ et en prenant le \og $C$-collapse\fg \,  de la partition obtenue. Ceci est mentionné dans \cite{B1},
 en haut de la page 174, et l'on s'en convainc en faisant l'exercice combinatoire.

 \medskip 

Si $M < x_0$, le symbole devient :
 \begin{equation}\label{newsymb1}
 \begin{pmatrix}
  M &&x_0+1 &&  \ldots &&  \ldots && x_{2p}+1 \\                       
 &M+1 &&\ldots    && \ldots  &&  x_{2p-1}+1 &
\end{pmatrix} 
\end{equation}

S'il existe $j$ tel que $x_{2j-1}<M+j<x_{2j}$, le symbole devient :
  \begin{equation}\label{newsymb2}
 \begin{pmatrix}
 x_0&&   \ldots  && M+j  &&x_{2j}+1&&\ldots    &x_{2p}+1 \\
&\ldots&&x_{2j-1}&&M+j+1&&\ldots &\quad x_{2p-1}+1
\end{pmatrix} 
\end{equation}
Si  $x_{2j}\leq M+j\leq x_{2j+1}-1$  (on pose ici $x_{2p+1}=+\infty$ par convention) le symbole devient 
 \begin{equation}\label{newsymb3}
 \begin{pmatrix}
  x_0 &&  \ldots &&  x_{2j} &&M+j+1 &&  x_{2j+2}+1&& \ldots \\
& \ldots  && x_{2j-1}  &&M+j+1 &&   x_{2j+1}+1 &&\ldots  
\end{pmatrix} 
\end{equation}

\medskip
\medskip
\begin{prop}\label{pac1pac2}
Le paquet  $\Pi_{BV}({\caU^+}, G^+)$ est l'ensemble des facteurs de composition
 des $\Triv_{2M+1}  \star \bar X  $, où $\bar X$ décrit $\Pi_{BV}({\caU},G)$.
 Les représentations induites $ \Triv_{2M+1} \star \bar X $ sont soit irréductibles et l'on pose alors
  $\bar X^+_0= \Triv_{2M+1} \star \bar X $, soit de longueur $2$ et l'on pose
  $\bar X\star \Triv_{2M+1}=\bar X^+_1 \oplus \bar X^+_{-1}$.
  \end{prop}

\dem Il est clair que pour tout $\bar X\in\Pi_{BV}({\caU},G) $, $ \Triv_{2M+1} \star \bar X $, 
ainsi que tous ses sous-quotients irréductibles ont bien le caractère infinitésimal requis.
On peut donc appliquer la proposition \ref{PropB166}.
La formule (\ref{B166}) montre que 
$\Triv_{2M+1} \star \bar X $ est de longueur 
\begin{equation}\label{rap} \vert \bar A(\caU^+)\vert\; /\; \vert  \bar A(\caU)\vert= 
\vert \Pi_{BV}({\caU^+},G^+)\vert\;  /\;  \vert \Pi_{BV}({\caU},G)\vert= 1 \text{ ou }  2.\end{equation}
D'autre part, si $\bar X=\bar X(\lambda,\mu)$, alors $\bar X\star \Triv_{2M+1}$ contient 
le sous-quotient de Langlands $\bar X(\lambda^+,\mu^+)$, où $\lambda^+$ et $\mu^ +$ sont obtenus 
respectivement à partir de  $\lambda$ et $\mu$ en ajoutant les coordonnées $(M,M-1,\ldots, -M)$. Comme 
$\bar X\star \Triv_{2M+1}$ est un sous-quotient de $X(\lambda^+,\mu^ +)$ et que 
$\bar X(\lambda^+,\mu^+)$ a multiplicité 1 dans $X(\lambda^+,\mu^ +)$, on voit que  
 $\bar X(\lambda^+,\mu^+)$ apparaît avec multiplicité 1 dans $\bar X\star \Triv_{2M+1}$.
Cette représentation étant de longueur au plus deux, elle est sans multiplicité.

 Dans le  cas (\ref{newsymb3}),  (\ref{rap}) vaut 1 et l'on en déduit que pour tout   
 $\bar X\in \Pi_{BV}(\caU,G)$,  $\Triv_{2M+1} \star \bar X $ est irréductible. On note 
$\bar X_0^+$ cette représentation. 

 Dans les  cas (\ref{newsymb1}) et (\ref{newsymb2}),  (\ref{rap}) vaut 2.
 En effet, dans le cas  (\ref{newsymb1})  $\bar A(\caU^+)=\bbZ/2\bbZ\times   \bar A(\caU^+)$, le facteur 
 $\bbZ/2\bbZ$ supplémentaire venant de la paire $(M+1,x_0+1)$ dans le symbole, et 
  dans le cas (\ref{newsymb2}),  la contribution d'un facteur  $\bbZ/2\bbZ$ dû à la paire $(x_{2j-1},x_{2j})$
est remplacée par  un facteur $(\bbZ/2\bbZ)^2$ dû aux paires $(x_{2j-1},M+j)$ et $(M+j+1, x_{2j}+1)$.
Ainsi $ \Triv_{2M+1} \star \bar X $ est de longueur 2, 
et l'on pose $ \Triv_{2M+1} \star \bar X =\bar X_1^+\oplus \bar X_{-1}^+$. 
La description explicite des éléments des paquets  $\Pi_{BV}({\caU^+},G^+)$ et $\Pi_{BV}({\caU},G)$ 
faite dans la section  \ref{descrip} particulièrement en (\ref{defpaq2}) nous dit exactement ce que sont 
les $\bar X_0^+$,   $\bar X_1^+$ et   $\bar X_{-1}^+$ ci-dessus.
Si $\bar X$ est paramétré par $\hat{\boldsymbol{\epsilon}}$ (que l'on identifie à un caractère de $\bar A(\caU)$,
alors dans le  cas (\ref{newsymb3}), $\bar X_0^+$ est paramétré par le même $\hat{\boldsymbol{\epsilon}}$, 
et dans les cas (\ref{newsymb1}) et (\ref{newsymb2}),  $\bar X_1^+$ et   $\bar X_{-1}^+$ sont paramétrés
respectivement par  $\hat{\boldsymbol{\epsilon}}^\pm$, obtenus en mettant $\pm 1$ sur le facteur 
 $\bbZ/2\bbZ$ que l'on a ajouté.
 \qed

 \subsection{Effet de l'ajout de deux  blocs même taille  de bonne parité, cas {$\mathbf B_n$}}
 \label{ajoutB}
 
On suppose que $G$ est un groupe orthogonal impair de rang $n$. 
Soient $\caU$ une orbite nilpotente de bonne parité de 
${}^L\frg=\frs\frp_{2n}$  $\Pi_{BV}(\caU,G)$ et  $\caO$ l'orbite duale dans $\frg$.
Soient $d_{\caU}$ et $d_\caO$ les partitions de $2n$ et $2n+1$ respectivement associées à ces orbites.

Les blocs de $d_{\caU}$ sont donc pairs, et quitte à rajouter le bloc $0$, 
 il y en a un nombre impair on peut donc écrire 
\[ d_{\caU}={\scriptstyle [2 (x_{2p}-p),  2(x_{2p-1}-p+1),   \ldots,  2(x_{2p-2j}-p+j), 
 2(x_{2p-2j-1}-p+j+1)    , \ldots, 
2x_2-2, 2x_1 , 2x_0] } \]
de sorte que le symbole est 
 \begin{equation}\label{SymboB}
 \begin{pmatrix}
  x_0 &&  \ldots  && x_{2j} && \ldots && &&x_{2p}\\
&&& \ldots  && x_{2j+1}  && \ldots  &&  x_{2p-1} &
\end{pmatrix} 
\end{equation}
 
 \medskip 
 
 Ajoutons 2 blocs de taille $2M$ à $d_{\caU}$, pour former une nouvelle partition $d_{\caU^+}$,
 correspondant à une orbite nilpotente  $\caU^+$ dans ${}^L\frg^+=\frs\frp_{2(n+2M)}$.
 Soient $\caO^+$ sa duale dans $\frg^+=\frs \fro_{2(n+2M)+1}$. 
 On a $\caO^+=\{0\}_{2M}\star \caO$. On obtient $d_{\caO^+}$ en ajoutant 2 aux $2M$ premiers blocs
 de $d_\caO$ et en prenant le \og $B$-collapse\fg \,  de la partition obtenue. 

 \medskip 

Si $M > x_{2p}-p$, le symbole devient :
 \begin{equation}\label{newsymb1B}
 \begin{pmatrix} 
 x_0 &&  \ldots  && x_{2j} && \ldots && &&x_{2p}&&M+p+1\\
&&& \ldots  && x_{2j+1}  && \ldots  &&  x_{2p-1} &&M+p&
\end{pmatrix} 
\end{equation}

S'il existe $j$ tel que $x_{2j}<M+j<x_{2j+1}$, le symbole devient :
  \begin{equation}\label{newsymb2B}
 \begin{pmatrix}
  x_0 &&  \ldots  && x_{2j} && M+j+1&&\ldots && &&x_{2p}+1\\
&&& \ldots  && M+j&&x_{2j+1}+1  && \ldots  &&  x_{2p-1}+1 &
\end{pmatrix} 
\end{equation}
Si  $x_{2j-1}+1\leq M+j\leq x_{2j}$  (on pose ici $x_{-1}=-\infty$ par convention) le symbole devient 
 \begin{equation}\label{newsymb3B}
 \begin{pmatrix}
   x_0 &&  \ldots  && &&M+j&& x_{2j}+1 && \ldots &&x_{2p}+1\\
&&& \ldots && x_{2j-1}&&M+j  && \ldots  &&  x_{2p-1}+1 &
\end{pmatrix} 
\end{equation}

\medskip

Comme dans le cas {$\mathbf C_n$}, 
 \begin{equation}\label{rap2} \vert \bar A(\caU^+)\vert\; /\; \vert  \bar A(\caU)\vert= 
\vert \Pi_{BV}({\caU^+},G^+)\vert\;  /\;  \vert \Pi_{BV}({\caU},G)\vert= 1 \text{ ou }  2, \end{equation}
et plus précisément, ceci vaut 1 dans le cas  (\ref{newsymb3B}), et 2 dans les cas  (\ref{newsymb1B}), 
et  (\ref{newsymb2B}). La proposition \ref{pac1pac2} est alors valide dans ce contexte : 
\begin{prop}\label{pac1pac2B}
Le paquet  $\Pi_{BV}({\caU^+},G^+)$ est l'ensemble des facteurs de composition
 des $\Triv_{2M}\star \bar X$, où $\bar X$ décrit $\Pi_{BV}({\caU},G)$.
 Les représentations induites $\Triv_{2M}\star \bar X$ sont  irréductibles 
 dans le cas  (\ref{newsymb3B}) et l'on pose alors
  $\bar X^+_0=\Triv_{2M}\star \bar X$, et de longueur $2$ dans les cas  (\ref{newsymb1B}), 
et  (\ref{newsymb2B}), et l'on pose alors
  $ \Triv_{2M+1} \star \bar X=\bar X^+_1 \oplus \bar X^+_{-1}$.
  \end{prop}

\subsection{Effet de l'ajout de deux  blocs même taille  de bonne parité, cas {$\mathbf D_n$}}
 \label{ajoutC}
 
On suppose que $G$ est un groupe orthogonal pair de rang $n$. Soient $\caU$ une orbite nilpotente de bonne parité de 
${}^L\frg=\frs\fro_{2n}$  et  $\caO$ l'orbite duale dans $\frg$.
Soient $d_{\caU}$ et $d_\caO$ les partitions de  $2n$ respectivement associées à ces orbites.

Les blocs de $d_{\caU}$ sont donc impairs, et il y en a un nombre pair on peut donc écrire 
\[ d_{\caU}={\scriptstyle [ 2(x_{2p+1}-p)-1,2(x_{2p}-p)+1,    \ldots,  2(x_{2p-2j+1}-p+j)-1, 
 2(x_{2p-2j}-p+j)+1    , \ldots, 
2x_2-1, 2x_1-1 , 2x_0+1] } \]
de sorte que le symbole est 
 \begin{equation}\label{SymboC}
 \begin{pmatrix}
  x_0 &&  \ldots &&  x_{2j} && \ldots && x_{2p}\\
x_1&& \ldots  && x_{2j+1}  && \ldots  &&  x_{2p+1} 
\end{pmatrix} 
\end{equation}
 
 \medskip 
 
 Ajoutons 2 blocs de taille $2M+1$ à $d_{\caU}$, pour former une nouvelle partition $d_{\caU^+}$,
 correspondant à une orbite nilpotente  $\caU^+$ dans ${}^L\frg^+=\frs\fro_{2(n+2M+1)}$.
 Soient $\caO^+$ sa duale dans $\frg^+=\frs \fro_{2(n+2M+1)}$. 
 On a $\caO^+=\{0\}_{2M+1} \star \caO$. On obtient $d_{\caO^+}$ en ajoutant 2 aux $2M+1$ premiers blocs
 de $d_\caO$ et en prenant le \og $D$-collapse\fg \,  de la partition obtenue. 

 \medskip

S'il existe $j$ tel que $x_{2j-1}<M+j<x_{2j}$, le symbole devient :
  \begin{equation}\label{newsymb2C}
 \begin{pmatrix}
 x_0&&   \ldots  && x_{2j-2}&& M+j  &&x_{2j}+1&&\ldots    &x_{2p}+1 \\
x_1&&\ldots&&x_{2j-1}&&M+j+1&& x_{2j+1}+1&&  \ldots &&\quad x_{2p+1}+1
\end{pmatrix} 
\end{equation}
Si  $x_{2j}\leq M+j\leq x_{2j+1}-1$  (on pose ici $x_{-1}=-\infty$ et $x_{2p+1}=+\infty$ par convention) 
le symbole devient 
 \begin{equation}\label{newsymb3C}
 \begin{pmatrix}
  x_0 &&  \ldots &&  x_{2j} &&M+j +1&&  x_{2j+2}+1&& \ldots && x_{2p}+1\\
x_1& & \ldots  && M+j+1  && x_{2j+1}+1 && \ldots && && x_{2p+1}+1
\end{pmatrix} 
\end{equation}

Comme dans le cas {$\mathbf C_n$}, 
 \begin{equation}\label{rap3} \vert \bar A(\caU^+)\vert\; /\; \vert  \bar A(\caU)\vert= 
\vert \Pi_{BV}({\caU^+}, G^+)\vert\;  /\;  \vert \Pi_{BV}({\caU}, G)\vert= 1 \text{ ou }  2, \end{equation}
et plus précisément, ceci vaut 1 dans le cas  (\ref{newsymb3C}), et 2 dans le cas 
  (\ref{newsymb2B}). La proposition \ref{pac1pac2} est alors valide dans ce contexte.
\medskip

 \subsection{Ajout de blocs de mauvaise parité\label{mauvaiseparite}}

 Dans les sections qui suivent, on se préoccupe de l'ajout de deux blocs de mauvaise parité. 
 La situation générale est la suivante: soit $a$ un entier de mauvaise parité et soit 
 ${\mathcal U}^-$ une orbite nilpotente de ${}^L\frg^-$  de bonne parité pour un groupe $G^-$ de même type 
 que $G$ mais de rang $a$ de moins. On suppose que l'orbite ${\mathcal U}'$ de $\frg$ obtenue 
 en ajoutant deux fois $a$ à la partition définissant ${\mathcal U}^-$ n'est pas spéciale.
  On note alors ${\mathcal U}$ la plus petite orbite spéciale contenant ${\mathcal U}'$. 
  Elle s'obtient en ajoutant à la partition définissant ${\mathcal U}^-$, les entiers $a+1$ et $a-1$.
   On note ${\mathcal O}^-$ l'orbite duale de ${\mathcal U} ^-$ et ${\mathcal O}$ l'orbite induite 
   $\{0\}_{a} \star \caO^-$. Alors  ${\mathcal U}$ est l'orbite duale de ${\mathcal O}$. 
    On peut déduire cette assertion de \cite{spaltenstein} III.11.7,   avec  une
    petite  difficulté car la dualité considérée par Spaltenstein n'est pas exactement celle considérée ici.
     Spaltenstein a défini sa dualité sans passer au groupe dual expliquant en \cite{spaltenstein} III.10.3
      le passage de l'une des dualités à l'autre. On préfère redonner une démonstration dans les cas importants 
      pour nous en laissant au lecteur le soin de généraliser.

  \subsubsection{Un lemme d'induction, cas {$\mathbf C_n$}}
 
Dans ce paragraphe, on suppose que $G$ est un groupe symplectique de rang $n$.
 Soient $\caU$ une orbite nilpotente de bonne parité de 
${}^L\frg=\frs\fro_{2n+1}$  et  $\caO$ l'orbite duale dans $\frg$.
Soient $d_{\caU}$ et $d_\caO$ les partitions de $2n+1$ et $2n$ respectivement associées à ces orbites.

Ecrivons
$d_{\caU}=[2d_1+1, 2d_2+1, 2d_3+1, \ldots , 2d_t+1]$. 

\begin{lemme}\label{lemmeinduc} On suppose que $d_1=d_2+1$. Soit $\caU^-$  l'orbite correspondant à la partition 
$d_{\caU^-}= [2d_3+1, \ldots , 2d_t+1]$
obtenue à partir de  $d_{\caU}$ en supprimant les deux premiers blocs.
  Alors $\caO=\{0\}_{2d_1}\star \caO^-$. De plus $\caO$ est aussi l'orbite duale de l'orbite ${\mathcal U}'$ 
  dont la partition est obtenue à partir de celle de ${\mathcal U}^-$ en ajoutant deux fois $2d_1$.
\end{lemme}

\dem On  note $\caO^-$ l'orbite duale de $\caU^-$. Rappelons ({\sl cf.} \cite{CMcG}) que l'on obtient 
les partitions correspondantes $\mathbf{d}_{\caO}$ et $\mathbf{d}_{\caO^-}$ à partir respectivement de
  $ \mathbf{d}_{\caU}$ et $\mathbf{d}_{\caU^-}$ en prenant les partitions transposées, 
  ${}^t \mathbf{d}_{\caU}$ et  ${}^t\mathbf{d}_{\caU^-}$ et en calculant leur \og $C$-collapse \fg.
Le plus grand bloc de ${}^t \mathbf{d}_{\caU}$ est $t$ et celui de ${}^t\mathbf{d}_{\caU^-}$ est $t-2$, et ces blocs
sont impairs, de multiplicité $2d_t+1$. En revanche, les multiplicités des autres blocs sont pairs,
 car elles sont de la forme
$2d_j+1-(2d_{j+1}+1)=2(d_j-d_{j+1})$ pour des $j$ tels que $d_j>d_{j+1}$.
Notons 
\[  {}^t\mathbf{d}_{\caU^-}=[a_1=t-2, \ldots , a_{2d_3+1}].\]
On a alors 
\[  {}^t\mathbf{d}_{\caU}=[a_1+2=t,  \ldots , a_{2d_3+1}+2, 2, \ldots, 2, 1, 1]\]
où la multiplicité de $2$ est  $2(d_2-d_3)$ et celle de $1$ est $2$.
Pour obtenir le \og $C$-collapse \fg\,  de  ${}^t\mathbf{d}_{\caU^-}$, on procède de la manière suivante. 
On note $i_0$ le 
plus grand indice tel que $a_{i_0}=a_1$ et $i'_0$ le plus grand indice tel 
que $a_{i'_0}<a_1$ et $a_{i'_0}$ est impair.
On définit ensuite $i_1$ comme grand indice tel que $a_{i'_0}=a_{i_1}$ et $i'_1$ le plus grand indice, s'il existe, 
 tel que $a_{i'_1}<a_{i'_0}$ et $a_{i'_1}$ est impair, et ainsi de suite jusqu'à ce que l'on arrive à un indice
 tel que $i_l$ soit défini, mais pas $i'_l$.  Le \og $C$-collapse \fg\,  de  ${}^t\mathbf{d}_{\caU^-}$
 est alors obtenu en remplaçant, pour $s=0, \ldots, l$,  les $a_{i_s}$ par $a_{i_s}-1$
  et pour $s=0, \ldots, l-1$,les 
$a_{i'_s}$ par $a_{i'_s}+1$. On obtient   $\mathbf{d}_{\caO^-}$, qui a $2d_3+1$ termes si $a_{2d_3+1}>1$ et 
$2d_3$ termes si $a_{2d_3+1}=1$.

Pour obtenir   $\mathbf{d}_{\caO}$, on procède de même, à partir de  ${}^t\mathbf{d}_{\caU}$ pour calculer 
son  \og $C$-collapse \fg.
 On remarque que les suites d'indices qui apparaissent coïncident jusqu'au rang $l$, 
 mais ici $i'_l$ est défini, et vaut 
 $2d_1$,  $i_{l+1}=2d_1+1   $ et $i'_{l+1}$ n'est pas défini. La partition $\mathbf{d}_{\caO}$
 se déduit de $\mathbf{d}_{\caO^-}$ en ajoutant 2 à tous les blocs, puis en ajoutant $2(d_2-d_3)+1$ 
 (resp. $2(d_2-d_3)+2$)
 blocs de longueur 2 si $a_{2d_3+1}>1$ (resp.   si $a_{2d_3+1}=1$). Ceci est bien l'orbite induite 
 $\caO^-\star\{0\}_{2d_1}$. Cela démontre la première assertion du lemme. Pour la deuxième, on procède de la même 
 façon en partant de ${\mathcal U}'$ au lieu de ${\mathcal U}$. Quand on calcule ${}^t{\mathbf{d} }_{{\caU}'}$, on 
  obtient la même partition que ${}^t{\mathbf{d}_{\caU}}$ sauf que le nombre $2$ intervient $2(d_2-d_3+1)$ et le 
  nombre $1$ n'intervient pas. Le \og  $C$-collapse \fg \,  de cette partition est le même que celui de
   ${}^t{\mathbf{d}}_{\caU}$.\qed

\begin{cor} \label{corinduc}
Soient $\Pi_{BV}(\caU,G)$ et     $\Pi_{BV}(\caU^-, G^-)$  les paquets associés par Barbasch-Vogan 
aux orbites $\caU$ et  $\caU^-$ respectivement pour les groupes  $G=\Sp_{2n}(\bbC)$ et $G^-=\Sp_{2(n-2d_1)}
(\bbC)$. Supposons que $\bar A(\caU^-)=(\bbZ/2\bbZ)^m$, de sorte que le cardinal de $\Pi_{BV}(\caU^-, G^-)$
soit $2^m$.  Considérons   les représentations induites de $G^-\times \GL_{2d_1}(\bbC)$ à $G$ :
\begin{equation}\label{EqindLI}\textstyle 
X_1=(\chi_{\frac{1}{2}, \frac{1}{2}}\circ \det_{2d_1})\star \bar X^ - , \qquad X_{-1}=
 (\chi_{\frac{1}{2},- 
\frac{1}{2}}\circ \det_{2d_1})\star \bar X^ -\end{equation} 
lorsque $\bar X^-$ décrit $\Pi_{BV}(\caU^-, G^-)$, et leurs sous-quotients de Langlands respectifs
$\bar X_1$ et $\bar X_{-1}$. Alors les $2^{m+1}$ représentations $\bar X_1$ et $\bar X_{-1}$ sont   
non équivalentes deux à deux,  et constituent donc les  $2^{m+1}$  représentations unipotentes du paquet
 $\Pi_{BV}(\caU,G)$. De plus les  $X_{-1}$ sont irréductibles, donc  $X_{-1}=\bar X_{-1}$. 
\end{cor}

\dem Le symbole de  $\caU$ s'obtient à partir de celui de $\caU^-$
en ajoutant $ \begin{matrix}
 {\qquad } d_1+t   \\
d_2+t
\end{matrix} $ à la droite de celui-ci. On voit donc que $\bar A(\caU)$ posède un facteur $\bbZ/2\bbZ$
supplémentaire, et le cardinal de $\Pi_{BV}(\caU,G)$ est $2^{m+1}$.

Les induites (\ref{EqindLI}),  ainsi que leurs sous-quotients irréductibles, ont bien le même  caractère infinitésimal
que celui des représentations du paquet $\Pi_{BV}(\caU,G)$.
On peut alors appliquer la première partie de   la proposition \ref{pac1pac2} qui nous dit que les sous-quotients
irréductibles de $X_1$ et $X_{-1}$ sont dans le paquet $\Pi_{BV}(\caU,G)$.
 Supposons que $\bar X^-=\bar X(\lambda^-, \mu^-)$, avec 
\[\lambda^-=(\lambda^-_1,\ldots, \lambda^-_{n-2d_1}), \quad  \mu^-=(\mu^-_1,\ldots, \mu^-_{n-2d_1}).\] 
Posons 
\begin{align*}
\lambda&= ( d_1, d_1-1, \ldots, -d_1+1,\lambda^-_1,\ldots, \lambda^-_{n-2d_1}),\\
 \mu&= ( d_1, d_1-1, \ldots, -d_1+1, \mu^-_1,\ldots, \mu^-_{n-2d_1}), \\
 \mu'&= ( d_1-3, d_1-1, \ldots, -d_1-1, \mu^-_1,\ldots, \mu^-_{n-2d_1}).\end{align*}
Alors $\bar X_1=\bar X(\lambda,\mu)$ 
et   $\bar X_{-1}=\bar X(\lambda,\mu')$.

Il est clair que les $2^{m+1}$ sous-quotients irréductibles ainsi obtenus en faisant varier  $\bar X$ dans 
$\Pi_{BV}(\caU,G)$ sont non équivalents deux à deux et donc décrivent entièrement
$\Pi_{BV}(\caU,G)$. Ainsi, tous les sous-quotients irréductibles des $X_1$ et $X_{-1}$ se trouvent parmi les 
 $\bar X_1$ et $\bar X_{-1}$

D'autre part, si $X_1$ (resp. $X_{-1}$) possède un autre sous-quotient irréductible que $\bar X_1$
(resp. $\bar X_{-1}$), on voit facilement avec le lemme \ref{crucial} que celui-ci est $\bar X_{-1}$
(resp. n'existe pas). 
Ceci montre en particulier  que $X_{-1}=\bar X_{-1}$. \qed

\begin{rmq}\label{rmqinduc} 
 Le lemme se généralise en remplaçant les blocs $d_1$ et $d_2=d_1-2$ par n'importe quelle 
paire $(d_{2j-1}, d_{2j})$ avec $d_{2j}=d_{2j-1}-2$ et c'est exactement ce qui est annoncé au début du paragraphe.
\end{rmq}

  \subsubsection{Un lemme d'induction, cas {$\mathbf B_n$}}
 On suppose que $G$ est un groupe orthogonal impair de rang $n$.
  Soient $\caU$ une orbite nilpotente de bonne parité de 
${}^L\frg=\frs\frp_{2n}$   et  $\caO$ l'orbite duale dans $\frg$.
Soient $d_{\caU}$ et $d_\caO$ les partitions de $2n$ et $2n+1$ respectivement associées à ces orbites.

Ecrivons
$d_{\caU}=[2d_1, 2d_2, 2d_3, \ldots , 2d_t]$. 

\begin{lemme}\label{lemmeinducB} On suppose que $d_2=d_3+1$. Soit $\caU^-$  l'orbite correspondant à la partition 
$d_{\caU^-}= [2d_1,2d_4, \ldots , 2d_t]$
obtenue à partir de  $d_{\caU}$ en supprimant les deuxième et troisième  blocs.
  Alors $\caO=\{0\}_{2d_2-1}\star \caO^-$ et c'est aussi l'orbite duale de l'orbite ${\mathcal U}'$ obtenue en ajoutant à la partition définissant ${\mathcal U}^-$ deux fois le nombre $2d_2-1$.
\end{lemme}

 \dem On  note $\caO^-$ l'orbite duale de $\caU^-$. 
Les partitions  $\mathbf{d}_{\caO}$ et $\mathbf{d}_{\caO^-}$ sont obtenues à partir respectivement de
  $ \mathbf{d}_{\caU}$ et $\mathbf{d}_{\caU^-}$ en ajoutant un bloc $1$, en prenant les partitions transposées, 
 et en calculant leur \og $B$-collapse \fg.
Notons $\mathbf{p}=[p_1, \ldots,p_{2d_1}]$ et $\mathbf{p}^-=[p_1^-, \ldots,p^-_{2d_1}]$   les partitions transposées 
obtenues. On a   $p_1=t+1$ et  $p_1^-=t-1$, $p^-_2=t-2$,  $p^-_{2d_3+1}=\cdots =p^-_{2d_1}=1$ et 
$p_1=p^-_1+2$, $p_2=p^-_2+2$, ..., $p_{2d_3}=p^-_{2d_3}+2$, $p_{2d_3+1}=p^-_{2d_3+1}+1=2$, ... ,
$p_{2d_2}=p^-_{2d_2}+1=2$, $p_{2d_2+1}=1$,...,  $p_{2d_1}=1$. 
Comme $p_1^-=t-1$, $p^-_2=t-2$, l'un des deux est pair, notons le $p^-_{i_0}$. Pour obtenir $\caO^-$ en calculant le 
\og $B$-collapse \fg, on prend celui qui est pair, on lui enlève 1, et on ajoute 1 au plus grand bloc pair
strictement plus petit, notons le $p^-_{i'_0}$. S'il n'existe pas de tel bloc pair, on ajoute juste un bloc 1, et
la procédure se termine. Sinon, comme il  y a un nombre pair de $p^-_j$ qui vaut
 $p^-_{i_0}$ (c'est la différence entre
les tailles de deux blocs de la partition $\mathbf{d}_{\caO^-}$), soit $i_1$ le plus grand indice tel que 
$p^-_{i_1}=p^-_{i'_0}$.  On enlève alors $1$ à $p^-_{i_1}$, et on ajoute 1 au plus grand bloc pair
strictement plus petit (s'il existe, sinon on ajoute un bloc 1), notons le $p^-_{i'_1}$.
On continue cette procédure jusqu'à ce que l'on ait défini $i_0, \ldots ,i_\ell, i'_0,\ldots,i'_{\ell-1}$, mais que 
$i'_{\ell}$ ne soit  pas défini, et l'on ajoute alors le bloc 1. Pour $\caO$, on fait de m\^eme, mais la procédure
se termine différemment car $p^-_{2d_3+1}=p^-_{2d_3+2}$ sont égaux à 2. Le premier est remplacé par 3 et le second 
par 1 et on ajoute un 1 final. On voit que l'on obtient $\mathbf{d}_{\caO}$ à partir de 
$\mathbf{d}_{\caO}$ en ajoutant 2 au $2d_2-1$ premiers blocs. Quand on travaille avec ${\mathcal U}'$ et la partition 
duale de la partition qui définit ${\mathcal U}'$ on obtient $\mathbf{p}$ sauf que une occurence de 2 est remplacée 
par deux occurences de 1. Ensuite on calcule comme ci-dessus le \og $B$-collapse \fg. \qed

\begin{cor} \label{corinducB}
Soient $\Pi_{BV}(\caU,G)$ et  et   $\Pi_{BV}(\caU^-,G^-)$  les paquets associés par Barbasch-Vogan 
aux orbites $\caU$ et  $\caU^-$ respectivement pour les groupes  $G=\SO_{2n+1}(\bbC)$ et $G^-=\SO_{2(n-2d_2)+3}
(\bbC)$. Supposons que $\bar A(\caU^-)=(\bbZ/2\bbZ)^m$, de sorte que le cardinal de $\Pi_{BV}(\caU^-,G^-)$
soit $2^m$.  Considérons   les représentations induites de $G^-\times \GL_{2d_2-1}(\bbC)$ à $G$ :
\begin{equation}\label{EqindLIB}\textstyle 
X_1=(\chi_{\frac{1}{2}, \frac{1}{2}}\circ \det_{2d_2-1})\star \bar X^-, \qquad X_{-1}
= (\chi_{\frac{1}{2},- 
\frac{1}{2}}\circ \det_{2d_2-1})\star \bar X^-\end{equation} 
lorsque $\bar X^-$ décrit $\Pi_{BV}(\caU^-,G^-)$, et leurs sous-quotients de Langlands respectifs
$\bar X_1$ et $\bar X_{-1}$. Alors les $2^{m+1}$ représentations $\bar X_1$ et $\bar X_{-1}$ sont   
non équivalentes deux à deux,  et constituent donc les  $2^{m+1}$  représentations unipotentes du paquet
 $\Pi_{BV}(\caU,G)$. De plus les  $X_{-1}$ sont irréductibles, donc  $X_{-1}=\bar X_{-1}$. 
\end{cor}

La démonstration est la même que dans le cas {$\mathbf C_n$}.

\begin{rmq}\label{rmqinducB} 
 Le lemme se généralise en remplaçant les blocs $d_2$ et $d_3=d_2-2$ par n'importe quelle 
paire $(d_{2j},d_{2j+1})$ avec $d_{2j+1}=d_{2j}-2$.
\end{rmq}

 \subsubsection{Un lemme d'induction, cas {$\mathbf D_n$}}

On suppose que $G$ est un groupe orthogonal pair de rang $n$.
 Soient $\caU$ une orbite nilpotente de bonne parité de 
${}^L\frg=\frs\fro_{2n}$   et  $\caO$ l'orbite duale dans $\frg$.
Soient $d_{\caU}$ et $d_\caO$ les partitions de $2n$  respectivement associées à ces orbites.
Ecrivons
$d_{\caU}=[2d_1+1, 2d_2+1, 2d_3+1, \ldots , 2d_t+1]$. 

\begin{lemme}\label{lemmeinducC} On suppose que $d_2=d_3+1$. Soit $\caU^-$  l'orbite correspondant à la partition 
$d_{\caU^-}= [2d_1+1,2d_4+1, \ldots , 2d_t+1]$
obtenue à partir de  $d_{\caU}$ en supprimant les deuxième et troisième  blocs.
  Alors $\caO= \{ 0 \}_{2d_2} \star \caO^-$ et c'est aussi l'orbite duale de l'orbite ${\mathcal U}'$ qui s'obtient à partir de ${\mathcal U}^-$ en ajoutant les deux blocs $2d_2$.
\end{lemme}

\dem  Maintenant, $t$ est pair. Les partitions  $\mathbf{d}_{\caO}$ et $\mathbf{d}_{\caO^-}$ sont obtenues à
 partir respectivement de  $ \mathbf{d}_{\caU}$ et $\mathbf{d}_{\caU^-}$ en prenant les partitions
  transposées,  et en calculant leur \og $D$-collapse \fg. Prenons des notations analogues à celles du cas 
  des groupes orthogonaux impairs. On a $p^-_{1}=t-2$ qui est pair et le cardinal de $\{ j\, \vert\,
  p^-_j=t-2 \}$ est $2d_t+1$ donc impair. Les autres blocs de 
$\mathbf{p}^-$ sont de multiplicité paires et $p^-_{2d_3+3}=\ldots =p^-_{2d_1+1}=1$.
Ensuite, tout se passe  exactement comme pour les groupes orthogonaux impairs.
\qed
\begin{cor} \label{corinducC}
Soient $\Pi_{BV}(\caU,G)$ et  et   $\Pi_{BV}(\caU^-,G^-)$  les paquets associés par Barbasch-Vogan 
aux orbites $\caU$ et  $\caU^-$ respectivement pour les groupes  $G=\SO_{2n}(\bbC)$ et $G^-=\SO_{2(n-2d_2)}
(\bbC)$. Supposons que $\bar A(\caU^-)=(\bbZ/2\bbZ)^m$, de sorte que le cardinal de $\Pi_{BV}(\caU^-,G^-)$
soit $2^m$.  Considérons   les représentations induites de $G^-\times \GL_{2d_2}(\bbC)$ à $G$ :
\begin{equation}\label{EqindLIC}\textstyle 
X_1=(\chi_{\frac{1}{2}, \frac{1}{2}}\circ \det_{2d_2})\star \bar X^-, \qquad X_{-1}
=(\chi_{\frac{1}{2},- \frac{1}{2}}\circ \det_{2d_2})\star \bar X^-\end{equation} 
lorsque $\bar X^-$ décrit $\Pi_{BV}(\caU^-,G^-)$, et leurs sous-quotient de Langlands respectifs
$\bar X_1$ et $\bar X_{-1}$. Alors les $2^{m+1}$ représentations $\bar X_1$ et $\bar X_{-1}$ sont   
non équivalentes deux à deux,  et constituent donc les  $2^{m+1}$  représentations unipotentes du paquet
 $\Pi_{BV}(\caU,G)$. De plus les  $X_{-1}$ sont irréductibles, donc  $X_{-1}=\bar X_{-1}$. 
\end{cor}

La démonstration est la même que dans les autres cas.

\begin{rmq}\label{rmqinducC} 
 Le lemme se généralise en remplaçant les blocs $d_2$ et $d_3=d_2-2$ par n'importe quelle 
paire $(d_{2j}, d_{2j+1})$ avec $d_{2j+1}=d_{2j}-2$.
\end{rmq}

\subsection{ Un autre lemme de  réduction} \label{plusdeux}
On suppose que $G$ est un groupe classique de rang $n$ et que $\caU$ est une orbite nilpotente de bonne parité
de ${}^L\frg$. On note $\mathbf{d}_\caU=[d_1,\ldots ,d_k]$ la partition associée, et l'on suppose que
$d_1\geq d_2+2$.
 On note maintenant $\caU^-$ l'orbite obtenue en remplaçant $d_1$ par $d_1-2$ dans 
  $\mathbf{d}_{\caU}$, c'est-à-dire que
 $\mathbf{d}_{\caU^-}=[d_1-2,\ldots d_t]$. C'est une orbite dans l'algèbre de Lie du groupe dual 
 de $G^-$, qui est un groupe de
 même type que $G$, et de rang $n-1$. L'orbite $\caU$ est induite, on a $\caU=\{0\}_1 \star \caU^-$.
 Notons $\caO$ et $\caO^-$ les orbites duales respectives de $\caU$ et $\caU^-$.

\begin{lemme}\label{lemmeplusdeux}
On suppose que $G$ est un groupe orthogonal, ou bien que $d_1>d_2+2$. 
Soit $M$ le sous-groupe de Levi standard de $G$ isomorphe à  $ \GL_1(\bbC)\times G^-$  
(et identifié à ce dernier). Alors 
l'orbite $ \{0\}_1 \times \caO^-$ de $\frmm$ est contenue dans $\caO\cap \frmm$.
\end{lemme}

\dem Lorsqu'on passe aux partitions transposées, celle de $\caU^-$ s'obtient à partir de celle de $\caU$ en enlevant 
deux blocs de taille 1. Pour obtenir $\caO^-$ et $\caO$, il faut calculer les \og $X$-collapses \fg\,  du type $X$ adéquat.
Cette opération ne change que la taille de blocs de mauvaise parité pour $G$ et $G^-$. Si $G$ est un groupe orthogonal,
1 est de bonne parité, et donc  $\mathbf{d}_{\caO^-}$  s'obtient de $\mathbf{d}_{\caO}$ en enlevant 
deux blocs de taille 1.  Si $G$ est un groupe symplectique, et si $d_1>d_2+2$, alors il y déjà au moins un bloc
 de taille $1$ dans   la partition transposée de $\mathbf{d}_{\caU^-}$, et la même conclusion prévaut. \qed  

\begin{rmq} Si $G$ est un groupe symplectique et si $d_1=d_2+2$, alors 
$\mathbf{d}_{\caO}$  s'obtient de $\mathbf{d}_{\caO^-}$ en ajoutant un bloc de taille 2.
\end{rmq}

\begin{cor}\label{corplusdeux} On suppose que $G$ est un groupe orthogonal, ou bien que $d_1>d_2+2$. 
Alors l'application naturelle de $\bar A(\caU^-)$ dans $\bar A(\caU)$ est un isomorphisme.
Toute représentation de 
$\Pi_{BV}(\caU,G)$ est obtenue de manière unique   en prenant le sous-quotient de Langlands  
$\bar X=\bar X(\lambda,\mu)$  de l'induite  $\chi_{\frac{d_1-1}{2},\frac{d_1-1}{2}} \star \bar X^-$
(ici $\chi_{\frac{d_1-1}{2},\frac{d_1-1}{2}}$ est un caractère de $\GL_1(\bbC)$), où
en posant  $\bar X^-=\bar X(\lambda', \mu')$, avec 
\[\lambda^-=(\lambda ^-_1,\ldots, \lambda^-_{n-1}), \quad  \mu^-=(\mu^-_1,\ldots, \mu^-_{n-1}),\] 
on a 
\[\lambda= ( \frac{d_1-1}{2},\lambda^-_1,\ldots, \lambda^-_{n-1}),
\qquad \mu= (\frac{d_1-1}{2},\mu^-,\ldots, \mu^-_{n-1}, \frac{d_1-1}{2}).\]

\end{cor}

\dem On est dans le cas d'une orbite $\caU$ lissement induite à partir de 
$\caU'$ ({\sl cf.} \cite{BV}, \S7). 
Si le symbole de $\caU$ est (\ref{Symbo}), celui de $\caU'$ obtenu en remplaçant $x_{2p}$ 
par $x_{2p}-1$.  Les autres assertion se déduisent aisément en 
comparant les  descriptions    des paquets $\Pi_{BV}(\caU^-,G^-)$ et 
$\Pi_{BV}(\caU,G)$ (qui ont même cardinal). \qed

\section{Un résultat sur les exposants}\label{DetParLan}

On énonce maintenant une propriété des exposants des représentations dans les paquets de Arthur unipotents.
Ceci se trouve déjà (de manière implicite) dans \cite{B1} et de manière totalement explicite dans \cite{barbaschtransparent}. La seule petite difficulté avec cette dernière référence c'est que, sauf erreur de notre part, certaines des représentations considérées ici ne le  sont pas  dans {\sl loc. cit.}
\begin{prop}\label{Propexposants}
Soit $G$ un groupe classique complexe de rang $n$ et $\caU$ une orbite nilpotente dans ${}^ L\frg$. 
Soit $\bar X \in \Pi_{BV}(\caU, G)$.  Alors les exposants
de $\bar X$ sont des caractères de $\bbC^\times$ de  la forme 
$\chi_{a,b}$ avec $a-b\in \{-\frac{1}{2},0,\frac{1}{2}\}$.
\end{prop}

\dem Soit $\mathbf{d}_{\caU}=[d_1,d_2,\ldots,d_k]$ la partition attachée à l'orbite $\caU$.
La démonstration se fait par récurrence sur le rang $n$, en utilisant les résultats de la section
\ref{LemmeRed}. On distingue selon les cas. Traitons par le cas où  $G$ est un groupe symplectique (cas {$\mathbf C_n$}).
Si $d_1=d_2$, soit $\caU^-$ l'orbite associée à la partition obtenue en enlevant les deux premiers blocs
$d_1$ et $d_2$. On est dans la situation du paragraphe \ref{ajout}, plus précisément dans le cas (\ref{newsymb3}). 
En particulier, un élément $\bar X$ de $\Pi_{BV}(\caU,G)$ est obtenu comme étant l'induite irréductible
$\bar X^-\star \Triv_{d_1}$ (voir la proposition \ref{pac1pac2}). La proposition pour $\bar X$ se déduit alors du
résultat appliqué à $\bar X^-$ par hypothèse de récurrence. Si $d_1>d_2+2$, on utilise cette fois le 
corollaire \ref{corplusdeux}, avec pour $\caU^-$ l'orbite associée à la partition obtenue en remplaçant 
$d_1$ par $d_1-2$. Dans le cas où $d_1=d_2+2$, on utilise cette fois le corollaire \ref{corinduc}, 
 avec pour $\caU^-$ l'orbite associée à la partition obtenue en enlevant les deux premiers blocs
$d_1$ et $d_2$ et l'hypothèse de récurrence appliquée aux représentations $\bar X_1$ et $\bar X_{-1}$ 
de ce corollaire. Le cas des groupes orthogonaux se traite de la m\^eme façon. \qed

\section{Identification des paquets de Barbasch-Vogan et d'Arthur}

Le résultat principal de cette section est de montrer que les constructions de Barbasch-Vogan et d'Arthur 
coïncident pour les groupes classiques complexes.

\begin{thm}\label{egalitepaq}
Soient $G$ un groupe classique complexe (non métaplectique) de rang $n$,  $\caU$ une orbite
nilpotente spéciale et paire de ${}^ L\frg$ et  $\psi_G$ le paramètre d'Arthur associé à  $\caU$. 
Soient $\eta$ un caractère de $A({{\mathcal U}})=A(\psi_G)$,   $X_\eta$ la représentation  
 associée par Arthur ({\sl cf.}  \S \ref{pacAr} et en particulier \ref{mult1})  
et $ X_\eta^{BV}$ celle associée par Barbasch-Vogan  ({\sl cf.}  \S \ref{descrip}) si $\eta$ se factorise
par $\overline{A}({\mathcal U})$.
 Si $\eta$ ne se factorise pas par $\overline{A}({\mathcal U})$,  on pose $ X^{BV}_\eta=0$.
On a alors $X_\eta= X^{BV}_\eta$.
En particulier
$\Pi_{BV}({}^L\caU, G)=\Pi(\psi,G)$.
\end{thm}

\begin{rmq} \label{AOAU}
Avant de prouver ce théorème il faut prévenir le lecteur que \cite{BV} utilise le groupe des caractères de 
$\overline{A}(\caO)$ et non pas $\overline{A}(\caU)$ (ici ${\mathcal O}$ est l'orbite duale de ${\mathcal U}$) pour 
paramétrer les représentations attachées à $\psi$. On passe de l'un à l'autre par un isomorphisme entre ces deux 
groupes (cf. l'introduction du paragraphe 12 de \cite{BV}). Et on utilise alors le corollaire 12.4 de \cite{BV} pour 
avoir les formules de transfert que l'on exprime, via l'isomorphisme utilisé dans \cite{BV}, en terme de caractères 
de $\overline{A}({\mathcal U})$. L'isomorphisme entre ces groupes de caractères n'est pas complètement évident, il 
traduit la tensorisation par le caractère signe dans l'ensemble des représentations du groupe de Weyl de $G$ qui est 
aussi celui de $^LG$ (cf. \cite{BV} 5.30 et suivant).
\end{rmq}

\begin{rmq}
Une orbite spéciale paire est soit de bonne parité, soit tous les blocs de la partition associée
sont de mauvaise parité. Cette dernière possibilité peut arriver dans les cas $\mathbf{B}_n$ et $\mathbf{D}_n$ 
et cette situation est décrite dans la section \ref{descrip} :  $A({\mathcal U})$
est trivial et la représentation associée par Barbasch-Vogan est une induite irréductible de la représentation
triviale d'un sous-groupe de Levi isomorphe à un produit de $\GL$.
Grâce au théorème \ref{redbp} et son corollaire, il en est de même de la représentation associée par Arthur.  
On peut donc supposer que  $\caU$ est une  orbite de bonne parité, et c'est ce que l'on fera dans ce qui suit.
\end{rmq}

\subsection{Le cas des groupes orthogonaux impairs}
On commence par le cas des groupes orthogonaux impairs qui est légèrement plus simple puis on généralisera.
La démonstration se fait par récurrence sur la longueur de la suite de Jordan-Hölder 
de $\psi=\Std_G\circ \psi_G$ en tant
que représentation de $W_{\mathbb C}\times \SL_2({\mathbb C})$. On
initialise aisément la récurrence si $\psi$ est irréductible: dans ce
cas ${\mathcal U}$ est l'orbite principale,  $A({\mathcal U })$ est
réduit à un élément et $X_1$ (la représentation associée par Arthur au
caractère nécessairement trivial de $A({\mathcal U})$) est la
représentation triviale. Il est facile de vérifier que $ X^{BV}_1$ est
aussi la représentation triviale.

On suppose maintenant que le résultat est établi pour tous les paramètres $\psi'_G$ tel que
$\psi'=\Std_G\circ \psi'_G$ ait une longueur strictement plus petite que celle de $\psi$.
 On utilise la caractérisation des $X_\eta$  par des propriétés de transfert 
endoscopique \cite{Art13} 2.1.1.
Soit  ${\bf H}=(H,s,\xi:{}^LH\rightarrow {}^LG,\ldots)$ une donnée endoscopique
elliptique avec $s$ dans le centralisateur de $\psi_G$. Le paramètre $\psi_G$ se factorise donc en 
$\xi\circ \psi_H$, où $\psi_H$ est un paramètre d'Arthur pour le groupe endoscopique $H$.
On suppose que $s\neq 1$ et en particulier le
groupe endoscopique $H$ est un produit de groupes $\SO_{2a+1}({\mathbb C})\times
\SO_{2b+1}({\mathbb C})$. On peut par hypothèse de récurrence  appliquer le théorème 
à $H$ et $\psi_H$. On sait que la représentation virtuelle $\sum_{\eta \in \widehat{A(\caU)}}\eta(ss_\psi) X_\eta$ 
est un transfert  de la représentation $X_{\psi_H}^{st}$ de $H$. Le même résultat vaut avec les $X_\eta$
remplacé par les $ X^{BV}_\eta$ d'après  \cite{BV}, \S 12.4. Comme les
représentations $X_\eta$ sont soit irréductibles et non équivalentes deux à deux, soit nulles, on en
déduit que l'ensemble des représentations $X_\eta$ coïncide avec
l'ensemble des représentations $X^{BV}_\eta$. En particulier la
somme des représentations dans chacun des ensembles ne dépend pas de
l'ensemble choisi, ce qui nous donne le cas $s=1$ dans l'égalité qui suit. On a  donc pour 
tout $s\in A({{\mathcal U}})$
\[ \Pi(s):=\sum_{\eta \in \widehat{A(\caU)}  } \eta(s) X_\eta=\sum_{\eta  \in \widehat{A(\caU)}} \eta(s)
 X^{BV}_\eta.\]
Par inversion de Fourier, on en déduit que  pour tout  $\eta \in \widehat{A(\caU)}$,
 \[X_\eta= |A({\mathcal U})|^{-1} \sum_{s\in A({\mathcal U})}
\eta(s) \Pi(s)= X^{BV}_\eta,\] ce qui termine la démonstration.  \qed

\subsection{Le cas général}
Pour imiter la démonstration ci-dessus, il faut initialiser la récurrence. 
Pour les groupes symplectiques, il n'y a pas de difficulté, c'est comme dans le cas des groupes orthogonaux
 impairs.  Si $\psi$ est irréductible il lui correspond la représentation triviale à la fois dans  \cite{Art13} 
 et dans \cite{BV}. Si $G$ est un groupe orthogonal pair, l'initialisation se fait quand $\psi$ 
 est de longueur deux. Pour des raisons de parité, $\psi$ ne peut pas être irréductible. 
 On suppose donc que $\psi$ est de longueur deux et on note ${\mathcal U}_\psi$ l'orbite nilpotente 
 déterminant $\psi$. On remarque que $\overline{A}({\mathcal U}_\psi)$ est réduit à un élément.
Il n'y a aussi qu'une représentation associée par Arthur à ce paramètre: notons $a\geq a'$ la taille 
des deux blocs de Jordan de l'orbite ${\mathcal U}_\psi$.  Toute représentation dans $\Pi(\psi)$ est,
 d'après \cite{pourhowe} obtenue par une correspondance de Howe. Mais ici, une seule correspondance de Howe est 
 possible,  c'est l'image de la représentation triviale de $\Sp_{a'-1}({\mathbb C})$. 
 Si $a'=1$, la représentation cherchée est tout simplement la représentation triviale.
  Il faut encore identifier ces deux représentations, celle construite par \cite{BV} et celle construite par 
  \cite{Art13}. D'après les constructions de \cite{Art13}, 
  c'est nécessairement la représentation dans le paquet de Langlands à l'intérieur du paquet d'Arthur,
   c'est-à-dire le quotient de Langlands de la représentation induite du caractère
$$\chi_{(a-1)/2,(a-1)/2}\otimes \cdots \otimes \chi_{(a'+1)/2,(a'+1)/2}\otimes \Triv_{a'}$$
d'un sous-groupe parabolique de Levi $(a-a')/2$ facteurs ${\mathbb C}^\times $ fois $\GL_{a'}({\mathbb C})$ (peu 
importe le parabolique choisi). Les réductions déjà faites ici montrent que cette représentation est bien celle 
construite par \cite{BV}.

Maintenant que l'on a initialisé la récurrence, on procède exactement comme dans le cas des groupes orthogonaux 
impairs: les données endoscopiques sont, si $G$ est un groupe symplectique un produit d'un groupe symplectique et 
d'un groupe spécial orthogonal pair tandis que si $G$ est un groupe spécial orthogonal pair, c'est le produit de deux 
groupes spéciaux orthogonaux pairs.\qed

\section{Démonstration du théorème \ref{Grossereduc}}\label{DemGrosseReduc}

\subsection{Un résultat de Vogan}
Le théorème suivant ({\sl cf.} \cite{B1}, Thm. 14.1)     est un résultat  attribué par 
Barbasch  à Vogan.
\begin{thm}\label{Vogan}
 Soient $G_0$ un groupe classique de rang $n_0$ et  $\bar X_0$ un caractère unitaire de $G_0$.
Soient, pour tout $i=1,\ldots, s$, un caractère unitaire $\chi_i\circ \textstyle
\det_{n_i}$ de $\GL_{n_i}(\bbC)$.

Soit $G$ le groupe classique de m\^eme type que $G_0$ et de rang $n=n_0+\sum_{i=1}^s n_i$. 
 Soit $\caO$ l'orbite nilpotente induite  $\left( \star_i (\{0\}_{n_i}) \right)\star \{0\}_{G_0}$
 dans $\frg$.  C'est donc une orbite de Richardson pour le sous-groupe parabolique $P=MN$ de $G$
 dont le facteur de Levi $M$ est isomorphe à $\left(\times_i \GL_{n_i}(\bbC)  \right) \times G_0$.
Si l'orbite $\caO$ est paire et normale, et si l'application moment $\caM: \, T^*(G/P) \rightarrow \caO$
est birationnelle, alors la représentation induite  
\[ \bar X= \left( \star_i \chi_i \circ \textstyle \det_{n_i}\right) \star \bar X_0  \] 
est irréductible, de front d'onde  $\mathrm{WF}(\bar X)=\overline \caO$.
\end{thm}

La liste des orbites nilpotentes  vérifiant ces conditions pour les groupes classiques est dans \cite{B1}, \S 14.3.
Nous allons maintenant, à la suite de Barbasch, montrer comment la démonstration du 
théorème \ref{Grossereduc} se ramène à ce résultat.

\subsection{Réduction au théorème \ref{Vogan}}

On suppose que $G$ est un groupe classique de rang $n$. Soient $\caU$ une orbite nilpotente de bonne parité de 
${}^LG$, $\psi=\psi_{\caU}$ le paramètre d'Arthur correspondant et  $\caO$ l'orbite duale dans $\frg$.
Soient $d_{\caU}$ et $d_\caO$ les partitions de $2n+1$ ou $2n$  associées à ces orbites.

Soit $\bar X\in \Pi(\psi,G)$. Avec les notations du  théorème \ref{Grossereduc}, on veut montrer que 
\begin{equation}\label{Redu} Y= \left(\star_{i=1}^s \chi_{i} \circ \textstyle \det_{a_i}\right) 
 \star  \left(\star_{j=1}^t \Triv_{b_j} \right)\star\bar X
\end{equation} est irréductible (les $b_j$ sont pairs dans les cas {$\mathbf C_n$} et {$\mathbf D_n$}, impairs dans 
le cas {$\mathbf B_n$} et {$\mathbf {Mp}_n$}). Pour cela, on se ramène au théorème  \ref{Vogan},
 en trois étapes qui réduisent le problème
à des cas de plus en plus particuliers, jusqu'à finalement se retrouver sous les bonnes hypothèses.
 On suit  de très près \cite{B1}.
Soit $\bar Y$ un sous-quotient irréductible de $Y$. Il s'agit de montrer que $Y=\bar Y$.

\medskip 

{\bf Première étape}. 
Montrons d'abord que l'on peut supposer que l'orbite $\caO$  vérifie  les hypothèses du théorème \ref{Vogan}, 
c'est-à-dire d'après  \cite{B1}, \S 14.3 que la  partition associée $\mathbf{d}_{\caO}$ ne possède
 que  des blocs de bonne parité
(paire pour les groupes symplectiques, impaire pour les groupes orthogonaux).
 
On raisonne par récurrence sur le nombre de blocs de mauvaise parité de $\caO$. Considérons 
l'orbite $\caU^+$ obtenue en rajoutant à $\mathbf{d}_{\caU}$ deux blocs de taille $T$ de bonne parité,
comme dans la section \ref{ajout} dont on adopte les notations dans ce qui suit. 
On a $\caO^+=\{0\}_{T}\star \caO$. On obtient $d_{\caO^+}$ en ajoutant 2 aux $T$ premiers blocs
 de $d_\caO$ et en prenant le \og $X$-collapse\fg \,  de la partition obtenue. 
 En choisissant bien $T$, on peut faire en sorte
 que $\caO^+$ ait deux blocs de mauvaise parité  de moins que $\caO$.
On a vu (sections \ref{ajout}, \ref{ajoutB}, \ref{ajoutC})  que selon les cas 
$\Triv_{T}\star \bar X=\bar X_1^+\oplus \bar X_{-1}^+$, 
ou bien $\Triv_{T}\star \bar X=\bar X_0^+$, où $ \bar X_1^+$,  $ \bar X_{-1}^+$ et $\bar X_0^+$
sont dans $\Pi(\psi^+,G^+)$. Montrons alors que si 
\begin{equation}\label{Yplus} Y^+= \left(\star_{i=1}^s \chi_{i} \circ \textstyle \det_{a_i}\right) 
 \star  \left(\star_{j=1}^t \Triv_{b_j} \right)\star\bar X^+
 \end{equation}
est irréductible pour tout $\bar X^+\in \Pi(\psi^+, G^+)$, il en est de même de $Y$.
On raisonne par l'absurde en supposant que $Y$ est réductible. 
Commençons par supposer qu'il existe $c\leq t$ tel que 
\[   \left(\star_{j=1}^{c-1}
\Triv_{b_j} \right) \star \bar X\]
soit irréductible
mais 
\[  Z= \left(\star_{j=1}^{c} \Triv_{b_j} \right)\star \bar X\]
soit réductible. Dans le groupe de Grothendieck, on calcule : 
\begin{align*} \Triv_{T}\star Z&=    \Triv_{T}\star \left(\star_{j=1}^{c} \Triv_{b_j} \right)\star  \bar X\\
&=  \left(\star_{j=1}^{c} \Triv_{b_j} \right)\star \Triv_{T}\star  \bar X \\
&=\begin{cases}  \left(\star_{j=1}^{c} \Triv_{b_j} \right) \star  \bar  X_1^+  + \left(\star_{j=1}^{c}
\Triv_{b_j} \right) \star \bar   X_{-1}^+  \\
 \text{ ou bien }  \left(\star_{j=1}^{c} \Triv_{b_j} \right) \star \bar  X_{0}^+
 \end{cases}
 \end{align*}

On  déduit du lemme \ref{crucial} que l'un des sous-quotients irréductibles de $Z$  admet un exposant 
$\chi_{a',b'}$   non-sphérique ($a'\neq b'$) avec $a', b'$ demi-entiers (cas {$\mathbf C_n$} et
$\mathbf D_n$) ou   $a', b'$ entiers (cas {$\mathbf B_n$} et {$\mathbf{Mp}_n$}). Il en est donc de même pour
  l'une des représentations $ \left(\star_{j=1}^{c} \Triv_{b_j} \right) \star \bar X^+ $, 
  avec $\bar X^+=\bar X_1^+$, $\bar X_{1}^+ $ ou $\bar  X_0^+$.
Or les exposants de cette représentation, irréductible car par hypothèse $Y^+$ l'est, 
sont soit des  exposants  $\chi_{a,b}$ de $\bar X^+$, donc avec  $a$ 
et $b$ entiers (cas {$\mathbf C_n$} et {$\mathbf D_n$}) ou  demi-entiers (cas {$\mathbf B_n$} et $\mathbf{Mp}_n$),
 soit des exposants  sphériques.  On aboutit  a une contradiction. 
Ainsi
\[  \left(\star_{j=1}^{t} \Triv_{b_j} \right)\star \bar X\]
est irréductible. Supposons  maintenant qu'il existe $c\leq s-1$ tel que 
\[ \star_{i=1}^{c-1} \left(\chi_{i} \circ \textstyle  \det_{a_i}\right)  \star \left(\star_{j=1}^{t}
\Triv_{b_j} \right) \star\bar X\]
soit irréductible, mais 
\[Z=\star_{i=1}^{c} \left(\chi_{i} \circ \textstyle  \det_{a_i}\right)  \star \left(\star_{j=1}^{t}
\Triv_{b_j} \right) \star\bar X\]
soit réductible. Posons $\chi_i=\chi_{a_i,b_i}$, où les $a_i$, $b_i$ sont des nombres complexes
vérifiant $a_i-b_i \in \bbZ$ et $a_i+b_i\in i\bbR$ (les $\chi_i$ sont des caractères unitaires de $\bbC^\times$). 
On peut aussi  supposer que les $\vert a_i-b_i \vert $ croissent avec $i$.

D'après le lemme \ref{crucial}, $Z$ admet un sous-quotient irréductible ayant un exposant  $\chi_{a',b'}$ avec 
$\vert a'-b'\vert >\vert a_c-b_c \vert$. Par un raisonnement analogue à celui ci-dessus, on 
 déduit que l'une des représentations 
 $\left(\star_{i=1}^{c} \chi_{i} \circ \textstyle \det_{a_i} \right)
  \star\left(\star_{j=1}^{t} \Triv_{b_j} \right) \star \bar X^+$, avec $X^+=X_1^+$, $X_{1}^+ $ 
 ou $X_0^+$ admet un exposant $\chi_{a',b'}$ avec 
$\vert a'-b'\vert >\vert a_c-b_c \vert$ et $a',b'$ demi-entiers 
(cas {$\mathbf C_n$} et {$\mathbf D_n$}) ou entiers (cas  {$\mathbf B_n$} et  {$\mathbf{Mp}_n$}) .
 Or cette représentation est irréductible car $Y^+$  
l'est par hypothèse, et  l'ensemble de ses exposants est l'union
des exposants  $\chi_{a,b}$ de $\bar X^ +$,  qui vérifient $\vert a-b\vert\leq \frac{1}{2}$ d'après 
la proposition  \ref{Propexposants}, et d'exposants de la 
forme   $\chi_{a_i+e,b_i+e}$, $e$ entier ou demi-entier, $i\leq c$.
 On aboutit  encore a une contradiction et l'on conclut que 
$Y$ est irréductible pour tout $\bar X \in \Pi(\psi,G)$ si  $Y^+$
l'est pour tout $\bar X^+\in \Pi(\psi^+, G^+)$. Par récurrence, on peut donc supposer que 
$\caO$ vérifie les hypothèses du théorème de Vogan.

\medskip 

{\bf Deuxième étape}.
Ensuite, nous nous ramenons au cas ou $\mathbf{d}_{\caU}$ n'a pas de \og trous \fg, c'est-à-dire que si 
\[\mathbf{d}_{\caU}=[d_1,\ldots ,d_k] \]
alors pour tout $i=1,\ldots ,k-1$, $d_{i+1}=d_i$ ou $d_i-2$.
La méthode est la même que ci-dessus, si $T<d_1$ n'apparaît pas dans  $\mathbf{d}_{\caU}$, 
 on ajoute deux blocs de taille $T$ pour obtenir une nouvelle partition $\mathbf{d}_{\caU^+}$
 vérifiant les mêmes propriétés que $\mathbf{d}_{\caU}$, mais avec un trou en moins.  
L'orbite duale  $\caO^+$ est alors égale à $\{0\}_{T}\star \caO$ et l'on  obtient $d_{\caO^+}$ 
en ajoutant 2 aux $T$  premiers blocs  de $d_\caO$ (il n'y a pas besoin cette de prendre le \og $X$-collapse\fg\,  
car $\caO^+$ n'a que des blocs de bonne parité et  ajouter 2 aux $T$  premiers blocs conserve cette propriété,
 c'est-à-dire que $\caO^+$ vérifie encore les hypothèses du théorème de Vogan). Le même argument que
ci-dessus montre que 
\[ Y= \left(\star_{i=1}^s \chi_{i} \circ \textstyle \det_{a_i}\right) 
 \star  \left(\star_{j=1}^t \Triv_{b_j} \right)\star\bar X \]
est irréductible pour tout $\bar X \in \Pi(\psi,G)$ si 
\[ Y^+= \left(\star_{i=1}^s \chi_{i} \circ \textstyle \det_{a_i}\right) 
 \star  \left(\star_{j=1}^t \Triv_{b_j} \right)\star\bar X^+\]
l'est pour tout $\bar X^+\in \Pi(\psi^+, G^+)$.
On se ramène donc en un nombre d'étape finie au cas  d'une orbite $\caU$ sans trou telle que 
$\caO$ vérifie les hypothèses du théorème de Vogan. 

Le fait que $\caU$ soit sans trou
implique  que tout élément $\bar X$ de $\Pi(\psi,G)$ est obtenu comme sous-quotient irréductible d'une
induite d'un caractère (pas nécessairement unitaire)  d'un sous-groupe de Levi isomorphe à un produit de $\GL$. 
On peut voir ceci en regardant la description explicite des éléments des paquets faites dans la section
 \ref{descrip}.

\medskip 

{\bf Troisième étape}.
La troisième et dernière étape consiste à se ramener au cas où $\bar X$ est de plus induite irréductible 
d'un caractère unitaire   d'un sous-groupe de Levi isomorphe à un produit de $\GL$, ce qui nous ramène
immédiatement aux  hypothèses du théorème de Vogan. Pour cela, utilisons les résultats
suivants :

\begin{lemme}
Il existe des entiers $i_1, \ldots, i_\ell$, avec $\ell+1\leq k$ et $d_{i_j}=d_{i_j+1}+2$, et une représentation 
$\bar Z$  induite irréductible  d'un caractère unitaire   d'un sous-groupe de Levi isomorphe à un produit de $\GL$
tels que $\bar X$ soit un sous-quotient irréductible de 
\begin{equation}\label{LemmeA}
\left( \star_{j=1}^\ell  \chi_{\frac{1}{2}, \frac{1}{2}}\circ\textstyle \det_{d_{i_j}-1}
  \right) \star 
\bar Z.
\end{equation}
\end{lemme}

\dem On peut préciser que $\bar Z$ est une représentation unipotente associée à l'orbite
obtenue de $\caU$  en enlevant à $d_{\caU}$ les blocs  $d_{i_j}$ et $d_{i_j+1}$ pour $j=1,\ldots, \ell$.

 On raisonne par récurrence sur le nombre de blocs de $\mathbf{d}_{\caU}$. 
Nous allons traiter le cas des groupes symplectiques, les autres cas étant similaires, en utilisant les
 lemmes de réduction de la section \ref{LemmeRed} adéquats. 
 Si $d_1=d_2$, soit  $\caU^-$ l'orbite nilpotente dont la partition associée  $\mathbf{d}_{\caU^-}$ 
 est obtenue en enlevant les  deux premiers blocs ($d_1$ et $d_2$) de $\mathbf{d}_{\caU}$. 
Dans ce cas, il découle des résultats de la section \ref{ajout} que $\bar X$ est une induite irréductible 
 de la forme $\bar X=\Triv_{d_1} \star \bar X^- $, où $\bar X^-\in \Pi(\psi^-, G^-)$, et $\psi^-$ 
est le paramètre d'Arthur correspondant à l'orbite  $\caU^-$. On applique ensuite l'hypothèse de récurrence 
à $\bar X^-$ : 

Il existe des entiers  $i_1, \ldots, i_\ell$, avec $\ell+1\leq k-2$ et $d_{i_j}=d_{i_j+1}+2$, et une représentation 
$\bar Z$  induite irréductible  d'un caractère unitaire   d'un sous-groupe de Levi isomorphe à un produit de $\GL$
tels que $\bar X^-$ soit un sous-quotient irréductible de 
\begin{equation}\label{LemmeB}
\left( \star_{j=1}^\ell  \chi_{\frac{1}{2}, \frac{1}{2}}\circ\textstyle \det_{d_{i_j}-1}  \right)
\star \bar Z^-.
\end{equation}
Alors $\bar X$ est sous-quotient de $\left( \star_{j=1}^\ell  \chi_{\frac{1}{2}, -\frac{1}
{2}}\circ\textstyle \det_{d_{i_j}-1}  \right)\star \bar Z$, où $\bar Z= \Triv_{d_1}\star \bar Z^-$.
D'après le corollaire  \ref{corinduc}, $\bar Z$ est bien irréductible, et l'on obtient l'assertion voulue. 

Si $d_1=d_2+2$, on utilise cette fois le lemme \ref{lemmeinduc} et son  corollaire \ref{corinduc}.
D'après ce corollaire, $\bar X$ est un sous-quotient   irréductible de l'induite 
$\chi\circ\det_{d_1-1} \star \bar X^-$, où $\chi=\chi_{\frac{1}{2}, \frac{1}{2}}$,
ou bien $\chi= \chi_{\frac{1}{2}, -\frac{1}{2}}$. On applique l'hypothèse de récurrence à $\bar X^-$ et on conclut comme ci-dessus.
\qed

\begin{lemme} Avec les notations qui précèdent,  les sous-quotients irréductibles de 
\begin{equation}\label{LemmeC}
\left( \star_{j=1}^\ell  \chi_{-\frac{1}{2}, -\frac{1}{2}}\circ\textstyle \det_{d_{i_j}-1}  \right)\star
\left( \star_{j=1}^\ell  \chi_{\frac{1}{2}, \frac{1}{2}}\circ\textstyle \det_{d_{i_j}-1}  \right)
 \star \bar Z .
\end{equation}
sont des induites irréductibles d'un caractère unitaire   d'un sous-groupe de Levi isomorphe à un produit de $\GL$
\end{lemme}

\dem On a pour les groupes généraux linéaires l'égalité suivante dans le groupe de Grothendieck  ({\sl cf.} \cite{BR}, 
Prop 14.1, formule de composition des bouts de séries complémentaires)
\[     \chi_{-\frac{1}{2}, -\frac{1}{2}}\circ\textstyle \det_{n} \times 
 \chi_{\frac{1}{2}, \frac{1}{2}}\circ\textstyle \det_{n} 
=\Triv_{n+1}\times \Triv_{n-1}+    \chi_{\frac{1}{2},- \frac{1}{2}}\circ\textstyle \det_{n} 
\times \chi_{-\frac{1}{2}, \frac{1}{2}}\circ\textstyle \det_{n}.\]
Tous les caractères apparaissant dans le membre de droite sont unitaires. 
Comme $\bar Z$ est elle-même une induite irréductible d'un caractère unitaire   d'un sous-groupe de Levi isomorphe 
à un produit de $\GL$, le résultat s'en déduit immédiatement. \qed 

\bigskip

Revenons à la démonstration du théorème \ref{Grossereduc}, ou plutôt de sa réduction
au théorème  \ref{Vogan}.  D'après les deux lemmes ci-dessus, 
\begin{equation}\label{finalred}
\left( \star_{j=1}^\ell  \chi_{-\frac{1}{2}, -\frac{1}{2}}\circ\textstyle \det_{d_{i_j}-1}  \right) 
\star \bar X
\end{equation}
a tous ses sous-quotients irréductibles qui sont des sous-quotients irréductibles de 
(\ref{LemmeC}), et sont donc des  induites irréductibles d'un caractère unitaire   d'un sous-groupe de Levi
 isomorphe à un produit de $\GL$.  Par le même argument que dans la deuxième étape, 
  (\ref{Redu}) est  irréductible si la même expression avec $\bar X$ remplacé par un sous-quotient irréductible 
  $\bar X^+$ de  (\ref{finalred}), en remarquant que ces $\bar X^+$ sont des éléments
  du paquet d'Arthur attachée à l'orbite obtenue à partir de $\caU$ en  ajoutant à 
  $d_{\caU}$ les blocs  $d_{i_j}$ et $d_{i_j+1}$ pour $j=1,\ldots, \ell$. 
Ceci termine la démonstration de la réduction du théorème \ref{Grossereduc} au thèorème \ref{Vogan}.\qed

\subsection{Le cas du caractère infinitésimal régulier}
On fixe $\psi_G$ un paramètre d'Arthur relatif au groupe $G$.
\begin{thm} On suppose que le caractère infinitésimal déterminé par $\psi_G$
({\sl cf.}(\ref{GC44})) est régulier. Alors $\Pi(\psi_G,G)$ est réduit à un élément,  
la représentation de Langlands à l'intérieur du paquet d'Arthur.
\end{thm}
\dem On reprend la notation $\psi_{u,bp}$ de (\ref{decomposition}). Il est clair que si $\psi_G$ définit un caractère 
infinitésimal régulier de $G$, a fortiori, $\psi_{u,bp}$ définit lui aussi un caractère infinitésimal régulier pour 
un groupe, a priori plus petit mais de même type que $G$.
D'après le théorème \ref{Grossereduc}, il suffit de prouver le théorème pour $\psi_{u,bp}$. 
Dans ce cas, il suffit de remarquer que soit $\psi_{u,bp}$ est nul soit est le paramètre d'un caractère.
 Dans tous les cas, $\overline{A}(\psi_{u,bp})$ est trivial car   $A(\psi_{u,bp})$ l'est déjà.
  \qed 

\begin{rmq} La conclusion du théorème est vraie sous la seule hypothèse que $\psi_{u,bp}$ définit
 un caractère infinitésimal régulier.
Le théorème s'applique en particulier aux paquets contenant des représentations ayant de la cohomologie
 pour un bon  système de coefficients. Ces paquets sont donc réduits à un élément.
\end{rmq}

 \section{Quelques compléments\label{complement}}
\subsection{Groupes spéciaux orthogonaux versus groupes orthogonaux}
Dans cette section, on suppose que $G=\SO_{2n}({\mathbb C})$. La théorie de l'endoscopie tordue qui transfère des 
représentations de $G$ vers des représentations de $\GL_{2n}({\mathbb C})$, ne permet pas de distinguer
 une représentation irréductible $X$ de G  de son image par les automorphismes provenant de 
 ${\mathbf O }_{2n}({\mathbb C})$ si cette représentation n'est pas invariante sous ces automorphismes. 
 Les paramètres $\psi$,  même complétés par la donnée d'un caractère de $A(\psi)$, ne sont donc pas suffisants. 
 Il y a toutefois un cas très important où le problème ne se pose pas car les représentations attachées à $\psi$ 
 sont invariantes sous l'action de ${\mathbf O }_{2n}({\mathbb C})$. Pour exprimer le résultat on 
 rappelle la décomposition de (\ref{decomposition})
  $$\psi=\rho \oplus \rho^* \oplus \psi_{u,bp}.$$
\begin{prop} On suppose que $\psi_{u,bp}\neq 0$, alors les représentations de G attachées à $\psi$ sont invariantes 
sous l'action de ${\mathbf O}_{2n}({\mathbb C})$.
\end{prop}
\dem Toute représentation irréductible attachée à $\psi$ est une induite irréductible d'une représentation $\tau \otimes 
X_{u,bp}$ où $\tau$ est une représentation convenable (unitaire irréductible) d'un produit de groupes $\GL$ et où 
$X_{u,bp}$ est attachée à $\psi_{u,bp}$. On note $2n_{u,bp}$ la dimension de la représentation de 
$\SL_2({\mathbb C})$ définie par $\psi_{u,bp}$. Il suffit donc de montrer que la représentation $X_{u,bp}$ est 
isomorphe à son image par l'automorphisme extérieur induit par ${\mathbf O}_{2n_{u,bp}}({\mathbb C})$.
 Or ceci est vrai 
car $X_{u,bp}$ s'obtient en considérant l'image par la correspondance de Howe d'une représentation irréductible d'un 
groupe symplectique de rang strictement plus petit que $n_{u,bp}$. \qed

\subsection{Effet de la conjugaison complexe}
 En considérant $G$ comme le groupe des points complexes d'un groupe déployé défini sur ${\mathbb R}$, on note
  $\sigma$ la conjugaison complexe. Si $G=\SO_{2n}({\mathbb C})$, on 
  peut aussi choisir une forme quasi-déployée non déployée de $G$ définie sur ${\mathbb R}$, d'où une conjugaison 
  complexe $\sigma_{nd}$; on sait que si $G=\SO_{2n}({\mathbb C})$, les conjugaisons $\sigma$ et 
  $\sigma_{nd}$ diffèrent par un automorphisme extérieur.
Il est donc clair que si $G=\SO_{2n}({\mathbb C})$, pour déterminer si une représentation 
est $\sigma$ ou $\sigma_{nd}$ invariante il faut pouvoir la distinguer de son image sous 
${\mathbf O}_{2n}({\mathbb C})$, ce que nous n'avons pas fait. 
C'est ce qui explique l'hypothèse que l'on met dans le cas où $G=\SO_{2n}({\mathbb C})$

\begin{thm} On fixe $\psi$ comme précédement et ici on suppose que si
 $G=\SO_{2n}({\mathbb C})$, alors $\psi_{u,bp}\neq 0$. Soit $X$ une représentation associée à $\psi$. 
 Alors $\sigma(X)\simeq X$ pour $X$ une conjugaison complexe définie sur $G$.
\end{thm}
\dem On décompose $\psi$ d'où la partie unipotente de bonne parité $\psi_{u,bp}$. Soit $X$ comme dans l'énoncé.
 On sait 
d'après le théorème \ref{Grossereduc} qu'il existe une représentation irréductible $X_{u,bp}$ associée à 
$\psi_{u,bp}$ et une représentation, $\rho$, induite de caractères unitaires d'un groupe $\GL_m({\mathbb C})$ tel que 
$X$ soit l'induite irréductible de $\rho\otimes X_{u,bp}$. Comme $\rho$ est une induite de caractères unitaires, on a 
$\sigma(\rho) \simeq \rho^* $. Et l'induite de $\rho^*\times X_{u,bp}$ est irréductible isomorphe à $X$ dans tous les 
cas; c'est ici que l'on utilise l'hypothèse faite quand $G=\SO_{2n}({\mathbb C})$. Il suffit donc de montrer que 
$\sigma(X_{u,bp})\simeq X_{u,bp}$. 

Le plus simple est de partir des formules explicites de \cite{BV}. On rappelle les constructions de loc. cite.
 On note ${\mathcal O}$ l'orbite duale de l'orbite unipotente associée à $\psi_{u,bp}$. Ici on est en droit de 
 supposer que $\psi=\psi_{u,bp}$. On note $W$ le groupe de Weyl de $G$, si $G$ est de type $B$ ou $C$. Dans le cas où 
 $G=\SO_{2n}({\mathbb C})$, on considère le groupe de Weyl de 
 ${\mathbf{O} }_{2n}({\mathbb C})$ au lieu de $\SO_{2n}({\mathbb C})$.  
 Pour $x\in \overline{A}(\caO)$, en suivant Lusztig, \cite{BV} associe une représentation $\sigma_x$ de $W$; 
 évidemment, il y a une difficulté dans le cas où $G=\SO_{2n}({\mathbb C})$ que l'on résoud ainsi. 
 Seules nous intéressent les composantes irréductibles de ces représentations ayant des vecteurs invariants sous le 
 stabilisateur dans $W$ de $\lambda_{{\mathcal O}}$ (avec les notations de loc. cite). Or ce stabilisateur contient 
 des éléments de ${\mathbf{O} }_{2n}({\mathbb C})\setminus \SO_{2n}({\mathbb C})$; 
 on peut donc étendre uniquement une telle représentation composante irréductible à $W$ en une représentation 
 irréductible. Ensuite \cite{BV} définit
  $$R_x:=|\mathrm{Stab}_W(\lambda_{\caO})|^{-1}\sum_{w\in W}{\mathrm{trace}}
 (\sigma_x(w)) X(\lambda_{\mathcal O},w\lambda_{{\mathcal O}})$$
et pour $\eta$ un caractère de $\overline{A}({\mathcal O})$, $$X_\eta:=|\overline{A}({\mathcal O})|^{-1}\sum_{x\in 
\overline{A}({\mathcal O})}\eta(x)R_x.$$
L'effet de la conjugaison complexe est de transformer $X(\lambda_{\mathcal O},w\lambda_{{\mathcal O}})$ en 
$X(w\lambda_{\mathcal O},\lambda_{{\mathcal O}})$ c'est-à-dire $X(\lambda_{\mathcal O},w^{-1}\lambda_{{\mathcal 
O}})$. Or pour tout $w\in W$, $w^{-1}$ est conjugué dans $W$ de $w$ d'où l'invariance par conjugaison complexe des 
représentations $R_x$. \qed

\subsection{Front d'onde des représentations dans les paquets d'Arthur\label{frontdonde}}
Soit $\psi$ un morphisme comme ci-dessus et soit $\bar X \in \Pi(\psi,G)$. On note ${\mathcal U}_{\psi}$ l'orbite 
nilpotente de ${}^L\frg$ qui est l'orbite de l'image par la différentielle en l'identite de $\psi$ 
d'un élément nilpotent non trivial de 
 $\frs\frl_2({\mathbb  C})$. 

\begin{cor}   Le front d'onde $WF(\bar X)$ est l'adhérence de l'orbite duale de ${\mathcal U}_{\psi}$
 dans $\frg$.
\end{cor}

\dem On utilise la décomposition de $\psi$ en 
sous-représen\-tations irréductibles comme dans (\ref{decomposition}). A la composante $\psi_{u,bp}$, on associe une 
orbite $\caU'$  dont on note $\caO'$ l'orbite duale. C'est donc une orbite nilpotente spéciale pour un groupe de même 
type que $G$ mais de rang en général plus petit que l'on note $G'$. On pose $\caO_{\psi}$ l'induite de $\caO'$ 
définie par
$$
\caO_{\psi}= *_{(\chi_p,N_p)}\{0\}_{N_p}*\caO'.$$
C'est une orbite spéciale. On montre que l'adhérence de cette orbite est $WF(\bar X)$.

Si $\psi=\psi_{u,bp}$ ce corollaire est un corollaire du fait que les représentations de \cite{Art13} sont celles 
construites par \cite{BV} et pour les représentations de \cite{BV} cela fait pratiquement partie de la définition. Le 
cas général est un corollaire du théorème \ref{Grossereduc}.

Il reste à montrer que cette orbite est la duale de ${\mathcal U}_\psi$. On reprend la décomposition complète de 
$\psi$
$$
\psi=\oplus_{\chi,p} \chi\boxtimes R_p \oplus \chi^{-1}\boxtimes R_p \oplus \psi_{u,bp}.$$
On pose $^LM$ le sous-groupe de Levi de $^LG$ isomorphe à $\times_{\chi,p}\GL_p({\mathbb C})\times ^LG'$. 
On remarque  que ${\mathcal U}_{\psi}$ contient l'orbite de  
$^LM$  qui est le produit des orbites nilpotentes principales sur chaque facteur $\GL_p({\mathbb C})$ et qui est 
${\mathcal U}'$ pour le facteur $^LG'$. On note ${\mathcal U}_{M}$ cette orbite. On note $M$ le sous-groupe de Levi 
de $G$ qui correspond à $^LM$ dans la dualité; c'est un produit avec les mêmes notations que ci-dessus 
$\times_{\chi,p}\GL_p(\bbC)\times G'$. On remarque encore que l'orbite nilpotente de $M$ qui est triviale sur 
tous les facteurs $\GL_p({\mathbb C})$ et vaut $\caO'$ sur $G'$ est la duale de ${\mathcal U}_M$. Pour obtenir le 
corollaire, il n'y a plus qu'à appliquer un résultat de Spaltenstein que l'on rappelle ci-dessous. \qed

\subsection{Rappel d'un résultat de Spaltenstein\label{spaltenstein}}
En \cite{spaltenstein}, \S III.11.7  la dualité et l'induction pour les orbites unipotentes sont reliées. 
Soit $^LM$ un sous-groupe de Levi de $^LG$ et soit ${\mathcal U}_M$ une orbite unipotente de $^LM$. On considère 
$\caU_M$ comme une orbite unipotente de $^LG$ par l'inclusion $i_{M,G}$. On note $M$ le sous-groupe de Levi de $G$ 
correspond à $^LM$ dans la dualité et $j_{M,G}$ l'induction. On note $d^*_G$ la dualité entre orbites unipotentes de 
$^LG$ vers $G$ et $d^*_M$ son analogue pour $^LM$, en suivant essentiellement les notations de \cite{spaltenstein} à 
ceci près qu'ici la dualité change le groupe (sauf pour $G=\SO_{2n}{\mathbb C})$ qui ici devient donc plus simple) 
alors que Spaltenstein utilise plutôt l'expression de la dualité sans changer de groupe (cf. \cite{spaltenstein} 
III.10.3 pour le rapport).

\begin{thm}(Spaltenstein \cite{spaltenstein} III.11.7) L'orbite duale de $i_{M,G}({\mathcal U}_M)$ est l'orbite 
induite de l'orbite duale de ${\mathcal U}_M$, c'est-à-dire que l'on a l'égalité
$$ j_{M,G} \circ d^*_M ({\mathcal U}_M)= d^*_G \circ i_{M,G}({\mathcal U}_M)
$$
\end{thm}
Il est dommage que ce ne soit pas exactement l'énoncé de Spaltenstein mais on laisse au lecteur le soin de s'y 
ramener. Pour l'application qui nous intéresse (corollaire \ref{frontdonde}), 
on remarque qu'il suffit de traiter, pas à pas, le 
cas d'un parabolique maximal de la forme $\GL_m({\mathbb C})\times G'$ et d'une orbite de la forme le produit de 
l'orbite principal sur le facteur $\GL$ est une orbite spéciale sur le facteur $G'$. Alors cela est un calcul sur les 
partitions duales analogue à ceux que l'on a fait dans le paragraphe \ref{mauvaiseparite}.

\subsection{Intersection entre les paquets d'Arthur\label{intersection}}
A tout morphisme $\psi$  comme dans tout ce travail, on va associer  un morphisme explicite 
 $\psi_{sp}$,  vérifiant
\begin{thm}
(i) Pour tout morphisme $\psi$, $\Pi(\psi,G)\subset \Pi(\psi_{sp},G)$

(ii) Soit $\psi,\psi'$ tel que $\Pi(\psi,G)\cap \Pi(\psi',G)\neq \emptyset$ alors  $\psi_{sp}=\psi'_{sp}$.
\end{thm}
Cela donne déjà des renseignements assez précis mais on pourrait aller plus loin et décrire 
vraiment les intersections.  Nous ne le faisons pas ici.

\medskip

\dem Pour montrer le théorème il faut bien sûr décrire $\psi_{sp}$.
 On  fixe $\psi$ et on note $\psi_{bp}$  la somme 
des sous-représentations de $W_\bbC\times \SL_2({\mathbb C})$ incluses dans $\psi$  de la forme 
$\Triv\boxtimes R_a$  avec $a$ de bonne parité, ou bien  de la forme $\chi_{1/2,-1/2}^{\pm 1}\boxtimes R_a$ 
avec $a$  de la mauvaise parité. On obtient  ainsi  une décomposition de $\psi$ sous la forme 
$\psi=\psi'\oplus \psi_{bp}$.
 On note $m_{bp}$ la dimension de la représentation $\psi_{bp}$. En restreignant $\psi_{bp}$ à $\SL_2({\mathbb C})$ 
 on obtient une orbite unipotente pour un groupe de même type que $^LG$ mais de rang $\lfloor m_{bp}/2
 \rfloor$ que l'on note   ${\mathcal U}_0$. On note ${\mathcal U}_{sp}$ l'orbite biduale de ${\mathcal U}_0$. 
Par définition $\psi_{bp,sp}$ est la représentation de $W_{\bbC}\times \SL_2({\mathbb C})$ de dimension $m_{bp}$  
qui définit ${\mathcal U}_{sp}$ par  restriction à $\SL_2(\bbC)$ et qui sur ${\mathbb C}^\times$ 
est triviale  sur les  composantes isotypiques des représentations de $\SL_2({\mathbb C})$
 correspondant à une représentation irréductible  dont la dimension est de bonne parité et qui
  pour les autres composantes isotypiques est une somme de caractères  
$\chi_{1/2,-1/2}$ et  $\chi_{-1/2,1/2}$, intervenant en nombres égaux. 
On pose alors: $\psi_{sp}=\psi'\oplus \psi_{bp,sp}$.
 On obtient bien ainsi un morphisme de $W_\bbC\times \SL_2({\mathbb C})$ dans $^LG$. 

On a la caractérisation intrinsèque de ${\mathcal U}_{sp}$ vue comme une orbite nilpotente: cette orbite est l'orbite 
duale du front d'onde de n'importe quelle représentation associée à $\psi_{bp}$: en effet on a calculé le front 
d'onde des représentations attachées à $\psi_{bp}$ en faisant $\psi=\psi_{bp}$ en \ref{frontdonde}. On vient donc de 
vérifier que $\psi_{bp,sp}$ est uniquement déterminé par n'importe quel élément de $\Pi(\psi_{bp})$.

Pour mieux comprendre la situation, décrivons ${\mathcal U}_{sp}$ du point de vue combinatoire
 ce qui permettra de caractériser uniquement  $\psi_{bp,sp}$ en fonction de $\psi_{bp}$. 
 On commence par l'orbite ${\mathcal U}_{u,bp}$ qui est définie par la 
restriction de $\psi_{u,bp}$ à $\SL_2({\mathbb C})$. On a sa partition, ${\mathcal E}_{u,bp}$  formée d'entiers tous 
de même parité, la bonne. La partition définissant ${\mathcal U}_0$ ({\sl cf.} ci-dessus) 
s'obtient en ajoutant à cette 
partition des entiers de mauvaise parité chacun un nombre pair de fois. On note ${\mathcal E}_{mp}$ ces entiers avec 
leur multiplicité. Soit $a\in {\mathcal E}_{mp}$. On dit que $a$ est à échanger si quand on ajoute deux copies de $a$ 
à la partition de ${\mathcal U}_{u,bp}$ on obtient une orbite qui n'est pas spéciale. On note ${\mathcal E}^-_{mp}$ 
l'ensemble ${\mathcal E}_{mp}$ dont on a retiré exactement deux fois tout entier $a$ qui est à échanger. Pour être 
clair si un entier $a$ à échanger intervient avec multiplicité $2r$ dans ${\mathcal E}_{mp}$, il intervient avec 
multiplicité $2(r-1)$ dans ${\mathcal E}_{mp}^-$. La partition de ${\mathcal U}_{sp}$ est alors l'union de ${\mathcal 
E}_{mp}^-$ avec ${\mathcal E}_{u,bp}$ et de l'ensemble formé des éléments $(a+1,a-1)$ où $a$ parcourt l'ensemble des 
entiers à échanger.
Avec cela, on vérifie facilement que $\psi_{sp}$ définit le même caractère infinitésimal que $\psi_{bp}$. En 
particulier il a la même propriété d'intégralité ou de demi-intégralité que le caractère infinitésimal de la 
représentation triviale.
On remarque pour la suite que cette propriété d'intégralité ou de demi-intégralité du caractère infinitésimal ainsi  
que le fait que  la restriction  de $\psi_{bp,sp}$ à $W_\bbC$ est par définition une somme de caractères 
$\chi_{1/2,-1/2}^{\pm 1}$ et de caractère triviaux fait que $\psi_{bp,sp}$ est uniquement déterminé par l'orbite 
unipotente ${\mathcal U}_{sp}$ (qui détermine la restriction de $\psi_{bp,sp}$ à $\SL_2({\mathbb C})$).

On démontre maintenant la propriété intermédiaire suivante. On a l'inclusion  $\Pi(\psi_{bp})\subset 
\Pi(\psi_{bp,sp})$. En effet soit $X$ un élément de $\Pi(\psi_{bp,sp})$; $X$ a comme caractère infinitésimal le 
caractère infinitésimal défini par $\psi_{bp,sp}$. On a vu que le front d'onde de $X$ est la fermeture d'une orbite 
nilpotente spéciale dont la duale est l'orbite associée à $\psi_{bp,sp}$. Comme la bidualité est 
l'identité dans l'ensemble des orbites spéciales (Spalstenstein, voir \cite{CMcG}), le front d'onde
 de $X$ est la fermeture de l'orbite duale de 
${\mathcal U}_{sp}$. On en conclut que $X\in \Pi(\psi_{bp,sp})$ en utilisant la définition même de \cite{BV} à 
condition que ${\mathcal U}_{sp}$ soit paire. C'est le cas quand ${\mathcal E}^-_{mp}$, défini ci-dessus, est 
l'ensemble vide. On se ramène facilement à ce cas puisque si ${\mathcal E}^-_{mp}$ ne fait plus que produire une 
induction $ (*_a\chi_{1/2,-1/2}\circ \det_{a} )* X'$ où $a$ parcourt l'ensemble des éléments de 
${\mathcal E}^-_{mp}$ avec une multiplicité moitié  et où $X'$ est convenable. Cela termine la preuve de l'assertion 
intermédiaire.  Au passage cela démontre (i) en toute généralité gr\^ace encore à \ref{Grossereduc}.

 Montrons  (ii) dans le cas particulier où $\psi$ et $\psi'$ vérifient $\psi=\psi_{bp}$ et $\psi'=\psi'_{bp}$. Dans 
 ce cas (ii) résulte du fait que $\psi_{bp,sp}$ est uniquement déterminé par un élément de $\Pi(\psi_{bp},G)
 \cap  \Pi(\psi'_{bp},G)$ si cet ensemble est non vide, comme on l'a vu ci-dessus.

Montrons (ii) en toute généralité. On fixe $\psi,\psi'$ comme dans l'énoncé de (ii) et on suppose que 
$\Pi(\psi,G)\cap \Pi(\psi',G)\neq \emptyset$ et on fixe $X$ dans cette intersection. On a donc défini $\psi_{bp}$ et 
$\psi'_{bp}$. On note $m$ la dimension de la représentation $\psi$ qui est aussi la dimension de la représentation 
$\psi'$ et on note $m_{bp}$ et $m'_{bp}$ les dimensions des représentations $\psi_{bp}$ et $\psi'_{bp}$. Et on 
raisonne par récurrence sur $\max(m-m_{bp}, m-m'_{bp})$. Si ce nombre est nul, on vient de démontrer (ii) et 
cela initialise la récurrence. Supposons que ce nombre ne soit pas nul. On fixe 
$p\in \frac{1}{2}{\mathbb Z}_{\geq 0}$ et 
$a\in {\mathbb N}$ tel que la représentation $\chi_{p,-p}\boxtimes R_a$ intervienne dans $\psi$ ou dans $\psi'$ et 
pas ni dans $\psi_{bp}$ ni dans $\psi'_{bp}$; ceci est tout à fait possible car une représentation irréductible 
intervenant dans $\psi$ et non dans $\psi_{bp}$ n'intervient par la définition même de $\psi_{bp}$ pas non plus dans 
$\psi'_{bp}$. Et on impose en plus à $p$    d'être maximum avec cette propriété. Par symétrie on suppose que 
$\chi_{p,-p}\boxtimes R_a$ intervient dans $\psi$ et on montre d'abord que cette représentation intervient 
nécesssairement dans $\psi'$. En effet, on note $\psi^-$ le morphisme qui se déduit de $\psi$ en enlevant la 
représentation $\chi_{p,-p}\boxtimes R_a$ et sa contragrédiente.  On a montré en \ref{Grossereduc} qu'il existe 
$X^-\in \Pi(\psi^-)$ tel que $X$ soit l'induite irréductible $\chi_{p,-p}\circ \det_a \star X^-$. 
Cela donne les paramètres de Langlands de $X$ en fonction de ceux de $X^-$. On  a aussi décrit $X$ comme 
induite à partir d'un élément de $\Pi(\psi'_{bp})$ en utilisant le fait que $X\in \Pi(\psi',G)$ et on voit, en 
regardant les paramètres de Langlands que cela force le fait que la représentation $\chi_{p,-p}\boxtimes R_a$ et sa 
contragrédiente interviennent dans $\psi'$. On note alors $\psi^{'-}$ l'analogue de $\psi^-$ et on obtient aussi 
$X^{'-}$ comme $X^-$. En regardant les paramètres de Langlands, on vérifie que $X^-\simeq X^{'-}$. D'où $\Pi(\psi^-) 
\cap \Pi(\psi^{'-})\neq \emptyset$. Et on obtient alors (ii), en appliquant l'hypothèse de récurrence à ces 
morphismes.\qed

\bigskip
\bibliographystyle{smfalpha}
\bibliography{MR2}

\providecommand{\bysame}{\leavevmode ---\ }
\providecommand{\og}{``}
\providecommand{\fg}{''}
\providecommand{\smfandname}{et}
\providecommand{\smfedsname}{\'eds.}
\providecommand{\smfedname}{\'ed.}
\providecommand{\smfmastersthesisname}{M\'emoire}
\providecommand{\smfphdthesisname}{Th\`ese}
\begin{thebibliography}{ABV92}

\bibitem[ABV92]{ABV}
{\scshape J.~Adams, D.~Barbasch {\normalfont \smfandname} D.~A. Vogan, Jr.} --
  \emph{The {L}anglands classification and irreducible characters for real
  reductive groups}, Progress in Mathematics, vol. 104, Birkh\"auser Boston,
  Inc., Boston, MA, 1992.

\bibitem[Art84]{Art84}
{\scshape J.~Arthur} -- {\og On some problems suggested by the trace
  formula\fg}, Lie group representations, {II} ({C}ollege {P}ark, {M}d.,
  1982/1983), Lecture Notes in Math., vol. 1041, Springer, Berlin, 1984,
  p.~1--49.

\bibitem[Art89]{Art89}
\bysame , {\og Unipotent automorphic representations: conjectures\fg},
  \emph{Ast\'erisque} (1989), no.~171-172, p.~13--71, Orbites unipotentes et
  repr{\'e}sentations, II.

\bibitem[Art13]{Art13}
\bysame , \emph{The endoscopic classification of representations}, American
  Mathematical Society Colloquium Publications, vol.~61, American Mathematical
  Society, Providence, RI, 2013, Orthogonal and symplectic groups.

\bibitem[Bar]{barbaschtransparent}
{\scshape D.~Barbasch} -- {\og Unipotent representations and theta
  correspondence\fg}, Notes d'un exposé à Dubrovnik, disponibles à \url{
  http://www.math.cornell.edu/~barbasch/}.

\bibitem[Bar89]{B1}
{\scshape D.~Barbasch} -- {\og The unitary dual for complex classical {L}ie
  groups\fg}, \emph{Invent. Math.} \textbf{96} (1989), no.~1, p.~103--176.

\bibitem[Bar03]{Baruch}
{\scshape E.~M. Baruch} -- {\og A proof of {K}irillov's conjecture\fg},
  \emph{Ann. of Math. (2)} \textbf{158} (2003), no.~1, p.~207--252.

\bibitem[BR10]{BR}
{\scshape A.~I. Badulescu {\normalfont \smfandname} D.~Renard} -- {\og Unitary
  dual of {${\rm GL}(n)$} at {A}rchimedean places and global
  {J}acquet-{L}anglands correspondence\fg}, \emph{Compos. Math.} \textbf{146}
  (2010), no.~5, p.~1115--1164.

\bibitem[BV85]{BV}
{\scshape D.~Barbasch {\normalfont \smfandname} D.~A. Vogan, Jr.} -- {\og
  Unipotent representations of complex semisimple groups\fg}, \emph{Ann. of
  Math. (2)} \textbf{121} (1985), no.~1, p.~41--110.

\bibitem[CM93]{CMcG}
{\scshape D.~H. Collingwood {\normalfont \smfandname} W.~M. McGovern} --
  \emph{Nilpotent orbits in semisimple {L}ie algebras}, Van Nostrand Reinhold
  Mathematics Series, Van Nostrand Reinhold Co., New York, 1993.

\bibitem[GF]{GG}
{\scshape W.~T. Gan {\normalfont \smfandname} G.~Fan} -- {\og The
  {L}anglands-{W}eissman program for {B}rylinski-{D}eligne extensions\fg},
  prépublication, \url{http://arxiv.org/abs/1409.4039 }.

\bibitem[KV95]{KV}
{\scshape A.~W. Knapp {\normalfont \smfandname} D.~A. Vogan, Jr.} --
  \emph{Cohomological induction and unitary representations}, Princeton
  Mathematical Series, vol.~45, Princeton University Press, Princeton, NJ,
  1995.

\bibitem[Lan79]{Lang79}
{\scshape R.~P. Langlands} -- {\og Automorphic representations, {S}himura
  varieties, and motives. {E}in {M}\"archen\fg}, Automorphic forms,
  representations and {$L$}-functions ({P}roc. {S}ympos. {P}ure {M}ath.,
  {O}regon {S}tate {U}niv., {C}orvallis, {O}re., 1977), {P}art 2, Proc. Sympos.
  Pure Math., XXXIII, Amer. Math. Soc., Providence, R.I., 1979, p.~205--246.

\bibitem[LS79]{LSp}
{\scshape G.~Lusztig {\normalfont \smfandname} N.~Spaltenstein} -- {\og Induced
  unipotent classes\fg}, \emph{J. London Math. Soc. (2)} \textbf{19} (1979),
  no.~1, p.~41--52.

\bibitem[M{\oe}g]{pourhowe}
{\scshape C.~M{\oe}glin} -- {\og Paquets d'arthur sp\'eciaux unipotents aux
  places archim\'ediennes et correspondance de {H}owe\fg}, prépublication,
  \url{http://webusers.imj-prg.fr/~colette.moeglin/pourhowe.pdf }.

\bibitem[M{\oe}g96]{duke}
{\scshape C.~M{\oe}glin} -- {\og Repr\'esentations quadratiques unipotentes des
  groupes classiques {$p$}-adiques\fg}, \emph{Duke Math. J.} \textbf{84}
  (1996), no.~2, p.~267--332.

\bibitem[Spa82]{spaltenstein}
{\scshape N.~Spaltenstein} -- \emph{Classes unipotentes et sous-groupes de
  {B}orel}, Lecture Notes in Mathematics, vol. 946, Springer-Verlag, Berlin-New
  York, 1982.

\bibitem[Tad09]{Tadic}
{\scshape M.~Tadi{\'c}} -- {\og $\mathrm{GL}(n,\mathbb{C})\, \hat{}$ and
  $\mathrm{GL}(n,\mathbb{R})\, \hat{}$\fg}, Automorphic forms and
  {$L$}-functions {II}. {L}ocal aspects, Contemporary Mathematics, Amer. Math.
  Soc., Providence, R.I., 2009, p.~285--313.

\bibitem[Vog86]{VogGL}
{\scshape D.~A. Vogan, Jr.} -- {\og The unitary dual of {${\rm GL}(n)$} over an
  {A}rchimedean field\fg}, \emph{Invent. Math.} \textbf{83} (1986), no.~3,
  p.~449--505.

\bibitem[Wei]{Weiss}
{\scshape M.~Weissman} -- {\og {$L$}-groups and parameters for covering
  groups\fg}, prépublication, \url{ http://arxiv.org/abs/1507.01042 }.

\bibitem[Zhe74]{Zhe}
{\scshape D.~P. Zhelobenko} -- \emph{Garmonicheskii analiz na poluprostykh
  kompleksnykh gruppakh {L}i}, Izdat. ``Nauka'', Moscow, 1974, Sovremennye
  Problemy Matematiki. [Current Problems in Mathematics].

\end{thebibliography}

\end{document}